\documentclass[AMA,Times1COL]{WileyNJDv5}

\usepackage{cleveref}
\usepackage{subcaption}
\usepackage{algpseudocode}
\usepackage{amsmath,amsthm,amstext,amssymb}
\usepackage{upgreek}
\usepackage{ulem}
\usepackage{enumerate}
\usepackage[shortlabels]{enumitem}
\usepackage{mathrsfs}
\usepackage{inconsolata}
\usepackage{fancyhdr}
\usepackage{booktabs}
\usepackage{caption}
\usepackage{subcaption}
\usepackage{placeins}
\usepackage{makecell}
\usepackage{listofitems}
\usepackage{comment}
\usepackage{multirow}
\usepackage{adjustbox}

\algtext*{EndWhile}
\algtext*{EndFor}
\algtext*{EndIf}

\newtheorem*{theorem*}{Theorem}
\newtheorem*{definition*}{Definition}
\newtheorem*{remark*}{Remark}
\newtheorem*{corollary*}{Corollary}

\usepackage{stmaryrd}
\SetSymbolFont{stmry}{bold}{U}{stmry}{m}{n}
\DeclareFontFamily{U}{mathx}{}
\DeclareFontShape{U}{mathx}{m}{n}{<-> mathx10}{}
\DeclareSymbolFont{mathx}{U}{mathx}{m}{n}
\DeclareMathAccent{\widehat}{0}{mathx}{"70}
\DeclareMathAccent{\widecheck}{0}{mathx}{"71}

\newcommand{\fboxc}[2]{{\fboxrule=0pt\fcolorbox{#1}{#1}{{#2}}}}
\definecolor{uniSlightblue}{HTML}{00BEFF}

\definecolor{uniSlightblue}{HTML}{00BEFF}

\let\T\relax
\newcommand{\T}[1]{{#1}^{\sf  T}}
\newcommand{\normg}[1]{\left\lVert #1 \right\rVert}
\newcommand{\norm}[1]{| #1 |}

\renewcommand{\div}[1]{{\rm div }\left( #1 \right)}

\newcommand{\sty}[1]{{\boldsymbol{#1}}}

\newcommand{\styy}[1]{{\mathbb{#1}}}

\newcommand{\fe}{\sty{ e}}
\newcommand{\ff}{\sty{ f}}
\newcommand{\fg}{\sty{ g}}

\newcommand{\fq}{\sty{ q}}
\newcommand{\fr}{\sty{ r}}
\newcommand{\fs}{\sty{ s}}

\newcommand{\fu}{\sty{ u}}

\newcommand{\fw}{\sty{ w}}
\newcommand{\fx}{\sty{ x}}
\newcommand{\fy}{\sty{ y}}

\newcommand{\fB}{\sty{ B}}

\newcommand{\fG}{\sty{ G}}

\newcommand{\fI}{\sty{ I}}

\newcommand{\zero}{\sty{ 0}}

\newcommand{\ffC}{\styy{ C}}
\newcommand{\ffD}{\styy{ D}}

\newcommand{\ffI}{\styy{ I}}

\newcommand{\ffN}{\styy{ N}}

\newcommand{\ffP}{\styy{ P}}

\newcommand{\ffR}{\styy{ R}}

\newcommand{\ffZ}{\styy{ Z}}
\newcommand{\fkappa}{\mbox{\boldmath $\kappa$}}
\newcommand{\fsigma}{\mbox{\boldmath $\sigma$}}
\newcommand{\fdelta}{\mbox{\boldmath $\delta$}}
\newcommand{\feps}{\mbox{\boldmath $\varepsilon $}}

\newcommand{\cA}{{\cal A}}
\newcommand{\cB}{{\cal B}}

\newcommand{\cF}{{\cal F}}

\newcommand{\cK}{{\cal K}}
\newcommand{\cL}{{\cal L}}
\newcommand{\cM}{{\cal M}}

\newcommand{\cO}{{\cal O}}

\newcommand{\cR}{{\cal R}}
\newcommand{\cS}{{\cal S}}
\newcommand{\cT}{{\cal T}}
\newcommand{\cU}{{\cal U}}

\newcommand{\pd}[2]{\displaystyle\frac{\displaystyle\partial #1}{\displaystyle\partial #2}}

\newcommand{\ol}[1]{\overline{#1}}
\renewcommand{\ul}[1]{\underline{#1}}
\newcommand{\ull}[1]{\ul{\ul{#1}}}
\newcommand{\ulll}[1]{\ul{\ull{#1}}}

\definecolor{uniSlightblue}{HTML}{00BEFF}
\definecolor{uniSblue}{HTML}{004191}
\definecolor{uniSlblue}{HTML}{00BEFF}
\definecolor{uniSgray}{RGB}{62, 68, 76}

\newcommand{\Her}[1]{{#1}^{\sf H}}
\newcommand{\Adj}[1]{{#1}^{\ast}}
\newcommand{\dual}{{}^\ast}
\newcommand{\dInt}[1]{{\, \mathrm{d} #1 }}
\newcommand{\vecOp}[1]{\mathrm{vec} \left( #1 \right)}
\newcommand{\conjOp}[1]{\mathrm{conj} \left( #1 \right)}

\newcommand{\field}[1]{ \boxed{ #1 } }
\newcommand{\diag}[1]{\mathrm{diag} \left( #1 \right)}
\newcommand{\eig}[1]{\mathrm{eig} \left( #1 \right)}

\newcommand{\orderOf}[1]{\cO \left( #1 \right)}
\newcommand{\her}{{}^\mathsf{H}}
\newcommand{\cond}{\mathrm{cond}_2}

\newcommand{\diffOp}{\mathsf{L}}
\newcommand{\diffOpParam}{{\diffOp_{\param}}}
\newcommand{\diffOpParamSample}{{\diffOp_{\param\sampleIdx}}}
\newcommand{\solOp}{\mathsf{G}}
\newcommand{\conv}{\ast}
\newcommand{\symOp}[1]{{\mathrm{sym} \left( #1 \right)}}
\newcommand{\bilinForm}{a}
\newcommand{\linForm}{l}

\newcommand{\realPart}[1]{\mathrm{Re}\left( #1 \right)}
\newcommand{\identity}{\mathrm{id}}
\newcommand{\mandelBasis}[1]{\fB^{( #1 )}}
\newcommand{\emptySet}{\varnothing}
\newcommand{\discrIntv}[2]{{\left[ #1 .. #2 \right]}}
\newcommand{\nDims}{d}
\newcommand{\nComp}{c}
\newcommand{\nChannels}{c}

\newcommand{\nWeights}{{N_{\rm w}}}
\newcommand{\nSamples}{{N_{\rm s}}}

\newcommand{\nElem}{{N_{\rm elem}}}
\newcommand{\nNodes}{{N_{\rm nodes}}}
\newcommand{\nDof}{{N_{\rm dof}}}
\newcommand{\nLayers}{{L}}
\newcommand{\nEpochs}{{N_\mathrm{epochs}}}

\newcommand{\rveLen}{L_0}
\newcommand{\unitVec}{\fe}
\newcommand{\unitVecI}{\unitVec_{(i)}}
\newcommand{\CG}{\emph{CG}}
\newcommand{\GMRES}{\emph{GMRES}}
\newcommand{\UNO}{\emph{UNO}}
\newcommand{\UNOCG}{\emph{UNO-CG}}
\newcommand{\UNOCGnaive}{\emph{UNO-CG (naive)}}

\newcommand{\JacCG}{\emph{Jac-CG}}
\newcommand{\Jac}{\emph{Jac}}
\newcommand{\FANS}{{\emph{FANS}}}

\newcommand{\ParPDE}{parametric PDE}
\newcommand{\ParPDEs}{parametric PDEs}
\newcommand{\domain}{\Omega}
\newcommand{\boundary}{\partial\domain}
\newcommand{\dimI}{i}
\newcommand{\dimJ}{j}
\newcommand{\edgeIM}{{\Gamma^{-}_\dimI}}
\newcommand{\edgeIP}{{\Gamma^{+}_\dimI}}
\newcommand{\edgeJM}{{\Gamma^{-}_\dimJ}}
\newcommand{\edgeJP}{{\Gamma^{+}_\dimJ}}
\newcommand{\edgesI}{{\Gamma^{\pm}_\dimI}}
\newcommand{\edges}{{\Gamma^{\pm}}}
\newcommand{\onEdgeM}[1]{\left. #1 \right|_{\Gamma^-}}
\newcommand{\onEdgeP}[1]{\left. #1 \right|_{\Gamma^+}}
\newcommand{\onEdgeIM}[1]{\left. #1 \right|_{\edgeIM}}
\newcommand{\onEdgeIP}[1]{\left. #1 \right|_{\edgeIP}}
\newcommand{\onEdgeJM}[1]{\left. #1 \right|_{\edgeJM}}
\newcommand{\onEdgeJP}[1]{\left. #1 \right|_{\edgeJP}}
\newcommand{\onBoundary}[1]{\left. #1 \right|_{\boundary}}
\newcommand{\paramSymbol}{\mu}
\newcommand{\param}{\paramSymbol}
\newcommand{\sol}{\fu}
\newcommand{\solD}{\ul{u}}
\newcommand{\resSymbol}{r}
\newcommand{\res}{\fr}
\newcommand{\resD}{\ul{r}}
\newcommand{\precResSymbol}{s}
\newcommand{\precRes}{\fs}
\newcommand{\precResD}{\ul{s}}
\newcommand{\test}{\fw}

\newcommand{\rhs}{\ff}
\newcommand{\rhsD}{\ul{f}}
\newcommand{\solParam}{{\sol_{\param}}}
\newcommand{\solDParam}{{\solD_{\param}}}

\newcommand{\solDParamSample}{{\solD\sampleIdx_{\param}}}
\newcommand{\rhsParam}{{\rhs_{\param}}}
\newcommand{\rhsDParam}{{\rhsD_{\param}}}

\newcommand{\rhsDParamSample}{{\rhsD\sampleIdx_{\param}}}
\newcommand{\convRate}{C}
\newcommand{\cgRate}{\convRate_{\CG}}

\newcommand{\cgIters}{{N_\mathrm{iter}}}

\newcommand{\cgItersBound}{{N_\mathrm{iter,est}}}

\newcommand{\Greens}{\fG}
\newcommand{\GreensParam}{\Greens_{\param}}
\newcommand{\GreensConv}{\overline{\Greens}}
\newcommand{\GreensConvParam}{\GreensConv_{\param}}
\newcommand{\femBasis}{\boldsymbol{\varphi}}
\newcommand{\eigMin}{\lambda_{\rm min}}
\newcommand{\eigMax}{\lambda_{\rm max}}
\newcommand{\refmat}{{\mathrm{r}}}
\newcommand{\pert}{\widecheck}
\newcommand{\temperature}{\vartheta}
\newcommand{\tempFluct}{\widetilde{\temperature}}
\newcommand{\tempMacro}{\ol{\temperature}}
\newcommand{\tempGrad}{\fg}

\newcommand{\tempGradMacro}{\ol{\tempGrad}}
\newcommand{\flux}{\fq}

\newcommand{\fluxUnit}{\mathrm{[W/(mK)]}}
\newcommand{\heatCond}{\kappa}
\newcommand{\heatCondTensor}{\fkappa}
\newcommand{\displacement}{\fu}
\newcommand{\dispFluct}{\widetilde{\displacement}}

\newcommand{\strain}{\feps}
\newcommand{\strainMacro}{\ol{\strain}}
\newcommand{\strainComp}{\varepsilon}
\newcommand{\strainMandel}{\ul{\strainComp}}
\newcommand{\strainMandelMacro}{\ol{\ul{\strainComp}}}
\newcommand{\stress}{\fsigma}

\newcommand{\stressComp}{\sigma}
\newcommand{\stressMandel}{\ul{\stressComp}}
\newcommand{\stiffTensor}{\ffD}

\newcommand{\stressUnit}{\mathrm{[GPa]}}
\newcommand{\noUnit}{\mathrm{[-]}}
\newcommand{\lameLambda}{\lambda^{\rm Lamé}}
\newcommand{\lameMu}{\mu^{\rm Lamé}}
\newcommand{\sampleI}{j}

\newcommand{\sampleIdx}{{}^{(j)}}

\newcommand{\layerI}{l}
\newcommand{\layerIFirst}{0}
\newcommand{\layerILast}{\nLayers}
\newcommand{\layerIdx}{{}^{\{\layerI\}}}
\newcommand{\layerIdxFirst}{{}^{\{\layerIFirst\}}}
\newcommand{\layerIdxLast}{{}^{\{\layerILast\}}}
\newcommand{\layerIdxNext}{{}^{\{\layerI+1\}}}
\newcommand{\layerIdxPrev}{{}^{\{\layerI-1\}}}
\newcommand{\nChannelsLayer}{\nChannels_\layerI}

\newcommand{\nChannelsLayerPrev}{\nChannels_{\layerI-1}}
\newcommand{\nChannelsLayerFirst}{\nChannels_{\layerIFirst}}

\newcommand{\paramSpace}{\cM}
\newcommand{\linopSpace}{\mathfrak{L}}
\newcommand{\solSpace}{\cU}
\newcommand{\solSpaceD}{\solSpace^{\rm h}}
\newcommand{\SquareIntSpace}{L^2}
\newcommand{\SobolevSpace}{H^1}
\newcommand{\symSpace}[1]{{Sym \left( #1 \right)}}
\newcommand{\spdSpace}[1]{{Sym_+ \left( #1 \right)}}
\newcommand{\PrecIn}{\fr}
\newcommand{\PrecInD}{\ul{r}}
\newcommand{\PrecInSet}{\cR}
\newcommand{\PrecInDSet}{\PrecInSet^{\rm h}}
\newcommand{\PrecOut}{\fs}
\newcommand{\PrecOutD}{\ul{s}}
\newcommand{\PrecOutSet}{\cS}
\newcommand{\PrecInDFwd}{\widehat{\PrecInD}}

\newcommand{\PrecInSample}{\PrecIn \sampleIdx}
\newcommand{\PrecInDSample}{\PrecInD \sampleIdx}
\newcommand{\PrecInDXSample}{\PrecInD_{\mathrm{x}}^{(j)}}
\newcommand{\PrecInDYSample}{\PrecInD_{\mathrm{y}}^{(j)}}
\newcommand{\PrecInDZSample}{\PrecInD_{\mathrm{z}}^{(j)}}
\newcommand{\PrecInDFwdSample}{\PrecInDFwd \sampleIdx}
\newcommand{\PrecInDFwdXSample}{\PrecInDFwd^{(j)}_{{\rm x}}}
\newcommand{\PrecInDFwdYSample}{\PrecInDFwd^{(j)}_{{\rm y}}}
\newcommand{\PrecInDFwdZSample}{\PrecInDFwd^{(j)}_{{\rm z}}}

\newcommand{\PrecInDAdjXSample}{\Adj{\left(\PrecInDFwdXSample\right)}}
\newcommand{\PrecInDAdjYSample}{\Adj{\left(\PrecInDFwdYSample\right)}}
\newcommand{\PrecInDAdjZSample}{\Adj{\left(\PrecInDFwdZSample\right)}}
\newcommand{\PrecOutDFwdXSample}{\PrecOutDFwd^{(j)}_{{\rm x}}}
\newcommand{\PrecOutDFwdYSample}{\PrecOutDFwd^{(j)}_{{\rm y}}}
\newcommand{\PrecOutDFwdZSample}{\PrecOutDFwd^{(j)}_{{\rm z}}}

\newcommand{\PrecOutDFwd}{\widehat{\PrecOutD}}

\newcommand{\PrecOutSample}{\PrecOut \sampleIdx}
\newcommand{\PrecOutDSample}{\PrecOutD \sampleIdx}

\newcommand{\PrecOutDFwdSample}{\PrecOutDFwd \sampleIdx}
\newcommand{\PrecOutDSet}{\PrecOutSet^{\rm h}}
\newcommand{\Prec}{\mathsf{P}}
\newcommand{\PrecDSymbol}{P}
\newcommand{\PrecD}{\ull{\PrecDSymbol}}

\newcommand{\PrecSparse}{\ull{Q}}
\newcommand{\StiffMatSymbol}{A}
\newcommand{\StiffMat}{\ull{\StiffMatSymbol}}
\newcommand{\StiffMatParam}{\StiffMat_{\param}}
\newcommand{\StiffMatParamSample}{\StiffMat_{\param\sampleIdx}}
\newcommand{\IdentD}{\ull{I}}
\newcommand{\weight}{\theta}
\newcommand{\weights}{{\ul{\weight}}}

\newcommand{\featureDelta}{\delta}
\newcommand{\featureDeltaSample}{\delta \sampleIdx}
\newcommand{\featureA}{\ul{\alpha}}
\newcommand{\featureB}{\ul{\beta}}
\newcommand{\featureI}{m}
\newcommand{\featureIdx}{{}_{\featureI}}
\newcommand{\PrecEntries}{\ul{v}}

\newcommand{\localParam}{\psi}
\newcommand{\globalParam}{\Psi}
\newcommand{\nPrecC}{C}
\newcommand{\FundSol}{\Phi}
\newcommand{\FundSolMat}{\ull{\FundSol}}
\newcommand{\FundSolFourier}{\widehat{\FundSol}}
\newcommand{\FundSolMatFourier}{\ull{\FundSolFourier}}
\newcommand{\cgSearchDir}{\ul{d}}
\newcommand{\cgMatvecResult}{\ul{p}}
\newcommand{\cgStepWidth}{\alpha}
\newcommand{\cgDelta}{\gamma}
\newcommand{\cgError}{\ul{e}}
\newcommand{\cgI}{m}
\newcommand{\cgIdx}{{}^{(\cgI)}}
\newcommand{\cgIdxFirst}{{}^{(0)}}

\newcommand{\cgTol}{\epsilon}
\newcommand{\diracDelta}{\fdelta}
\newcommand{\fnoConvKernel}{\overline{\ull{K}}\layerIdx}
\newcommand{\fnoConvKernelFourier}{\widehat{\ull{K}}\layerIdx}
\newcommand{\fnocgConvKernelSymbol}{\widetilde{K}}
\newcommand{\fnocgConvKernelFourier}{\ull{\fnocgConvKernelSymbol}}
\newcommand{\fnocgFundSolSymbol}{\widetilde{\FundSol}}
\newcommand{\fnocgFundSolMat}{\ull{\fnocgFundSolSymbol}}
\newcommand{\fnocgBypass}{\widetilde{\ull{\fnoBypassSymbol}}}
\newcommand{\Loss}{\cL}

\newcommand{\lossSample}{l \sampleIdx}
\newcommand{\optimStep}{\Delta \weights}
\newcommand{\optimGrad}{\ul{g}}
\newcommand{\optimHess}{\ull{H}}
\newcommand{\trafo}{\cT}
\newcommand{\trafoFFT}{\cF}
\newcommand{\trafoFFTInv}{\trafoFFT^{-1}}
\newcommand{\trafoDST}{\cT_{\mathsf{ST}}}
\newcommand{\trafoInv}{\trafo^{-1}}

\newcommand{\trafoD}{\ull{T}}
\newcommand{\trafoDInv}{\trafoD^{-1}}
\newcommand{\trafoDHer}{\trafoD \her}

\newcommand{\modeSymbol}{k}
\newcommand{\mode}{\ul{\modeSymbol}}
\newcommand{\modeComp}{\modeSymbol}
\newcommand{\modeI}{i}
\newcommand{\modeIdx}{{}^{\langle \modeI \rangle}}
\newcommand{\modeIdxBypass}{{}^{\langle 0 \rangle}}
\newcommand{\modeIdxFirst}{{}^{\langle 1 \rangle}}
\newcommand{\modeIdxLast}{{}^{\langle \nModes \rangle}}
\newcommand{\nModes}{{k_{\mathrm{max}}}}
\newcommand{\modeSet}{\cK}
\newcommand{\modeSetConj}{\modeSet_{\rm c}}

\newcommand{\modeParam}{M}
\newcommand{\PIso}{\ffP^\mathrm{iso}}
\newcommand{\Isym}{\ffI^\mathrm{s}}

\newcommand{\fnoWeightSymbol}{R}
\newcommand{\fnoWeightTensor}{\ulll{\fnoWeightSymbol}\layerIdx}
\newcommand{\fnoAct}{\sigma\layerIdx}
\newcommand{\fnoFeatures}{\ul{w}}

\newcommand{\fnoKernelSymbol}{K}
\newcommand{\fnoKernel}{\ull{\fnoKernelSymbol} \layerIdx}
\newcommand{\fnoBypassSymbol}{W}
\newcommand{\fnoBypass}{\ull{\fnoBypassSymbol} \layerIdx}
\newcommand{\fnoProjOp}{\mathsf{V}_{\mathrm{lift}}}
\newcommand{\fnoBackprojOp}{\mathsf{V}_{\mathrm{proj}}}

\newcommand{\learningRate}{\eta}

\articletype{Research article}

\startpage{1}

\begin{document}

\title{Accelerating Conjugate Gradient Solvers for Homogenization Problems with Unitary Neural Operators}

\author[1]{Julius Herb}

\author[1]{Felix Fritzen*}

\authormark{HERB \textsc{et al.}}
\titlemark{Accelerating Conjugate Gradient Solvers for Homogenization Problems with Unitary Neural Operators}

\address[1]{\orgdiv{Institute of Applied Mechanics, Data Analytics in Engineering}, \orgname{University of Stuttgart}, \orgaddress{\state{Universitätsstr. 32, 70569 Stuttgart}, \country{Germany}}}

\corres{* Felix Fritzen, Institute of Applied Mechanics, Universitätsstr. 32, 70569 Stuttgart, Germany. \email{fritzen@mib.uni-stuttgart.de}}

\colorlet{revision}{black}
\colorlet{REVISION}{revision}


\abstract[Abstract]{Rapid and reliable solvers for parametric partial differential equations (PDEs) are needed in many scientific and engineering disciplines.
For example, there is a growing demand for composites and architected materials with heterogeneous microstructures.
Designing such materials and predicting their behavior in practical applications requires solving homogenization problems---typically governed by PDEs---for a wide range of material parameters and microstructures.
While classical numerical solvers offer reliable and accurate solutions supported by a solid theoretical foundation, their high computational costs and slow convergence remain limiting factors.
As a result, scientific machine learning is emerging as a promising alternative, aiming to rapidly approximate solutions using surrogate models.
However, such approaches often lack guaranteed accuracy and physical consistency.
This raises the question of whether it is possible to develop hybrid approaches that combine the advantages of both data-driven methods and classical solvers.
To address this, we introduce UNO-CG, a hybrid solver that accelerates conjugate gradient (CG) solvers using specially designed machine-learned preconditioners, while ensuring convergence by construction.
As a preconditioner, we propose Unitary Neural Operators (UNOs) as a modification of the established Fourier Neural Operators.
Our method can be interpreted as a data-driven discovery of Green's functions, which are then used much like expert knowledge to accelerate iterative solvers.
We evaluate UNO-CG on various homogenization problems involving materials with heterogeneous microstructures and millions of degrees of freedom.
Our results demonstrate that UNO-CG enables a substantial reduction in the number of CG iterations and is competitive with handcrafted preconditioners for homogenization problems that involve expert knowledge.
Moreover, UNO-CG maintains strong performance across a variety of boundary conditions, where many specialized solvers are not applicable, highlighting its versatility and robustness, which is supported by our extensive numerical study.}

\keywords{composite materials, computational homogenization, partial differential equations, conjugate gradient, machine learning, neural operators, preconditioning}

\jnlcitation{\cname{
\author{J. Herb}, and
\author{F. Fritzen}}.
\ctitle{Accelerating Conjugate Gradient Solvers for Homogenization Problems with Unitary Neural Operators} \cjournal{\it } \cvol{}}

\maketitle

\renewcommand\thefootnote{}
\footnotetext{\textbf{Abbreviations:} CG, conjugate gradient; UNO, Unitary Neural Operator; FNO, Fourier Neural Operator; PDE, partial differential equation; FEM, finite element method; FANS, Fourier-Accelerated Nodal Solvers; GMRES, generalized minimal residual method; FFT, fast Fourier transform; ML, machine learning; BC, boundary conditions.}

\renewcommand\thefootnote{\fnsymbol{footnote}}
\setcounter{footnote}{1}

\section{Introduction}
\label{sec:introduction}

\subsection{Motivation by application}
\label{ssec:introduction-motivation}

While it is assumed in many engineering applications that components are characterized by a homogeneous microstructure and can be described by explicit and closed-form material laws, this is not always the case in reality.
Real materials often exhibit a heterogeneous microstructure, which can significantly affect the material's properties.
Pronounced examples of this are metal-matrix composites (MMCs\cite{Chawla2013}), which consist of a metallic matrix reinforced by fibers or particles of another, often ceramic, material. MMCs can show exceptional material properties, making them highly demanded materials for challenging applications, such as in aerospace engineering, the automotive industry, and the biomedical sector \cite{Pooja2025}.

However, determining the effective material behavior of components with heterogeneous microstructures remains challenging. 
Multi-scale simulations are typically required to capture the nontrivial relation between macroscale behavior and microscale heterogeneities.
One established approach to address this challenge is the FE\textsuperscript{2} (FE square) method \cite{Schroder2014}, where a finite element simulation on the macroscale is coupled with microscale finite element simulations at each Gauss point and iteration of the macroscopic analysis. Despite modern computational resources, the FE\textsuperscript{2} approach remains prohibitively expensive from a computational perspective for many real-world applications.

Microscale simulations in this context represent homogenization problems \cite{Yvonnet2019}.
They are, e.g., crucial for the design of composite and architected materials \cite{Dirrenberger2019}, which may also involve topology optimization \cite{Osanov2016}.
The microscale simulations are coupled to the macroscale simulation through the macroscopic loading acting as input, and the resulting effective stress and stiffness information is handed back to the macroscale simulation.
This sought-after input-output relation depends heavily on the microstructure, a set of material parameters, and the macroscopic loadings, leading to a prohibitively high-dimensional input space.

From a mathematical point of view, these homogenization problems can be formulated as parametric partial differential equations (PDEs), with repeated solution queries involving different material parameters, microstructural geometry, and loading conditions. Over the past decades, numerous methods have been developed to accelerate the solution of such problems, which will be summarized in the following.

\subsection{Numerical iterative solvers for parametric PDEs}
\label{ssec:introduction-iterative-solvers}
{\color{revision}
The classical approach to solving parametric PDEs is to discretize them to an algebraic system of equations, for example, using the finite element method (FEM) or finite difference method.
The resulting linear or linearized system then consists of a stiffness matrix and a right-hand side vector that have to be assembled for a given parameter.
Direct solvers can yield a solution by factorizing stiffness matrix, e.g., using a LU or QR decomposition \cite{Stoer2002}.
The general drawback of this approach in a many-query scenario is that, for each new parameter, the stiffness matrix has to be reassembled and decomposed, which is computationally prohibitively expensive.
For this reason, direct solvers are typically infeasible for large-scale homogenization problems.
}

Unlike direct solvers, iterative solvers approximate the solution iteratively with steadily improving accuracy, typically performing one or more matrix-vector products with the stiffness matrix in each iteration. Notably, direct access to the stiffness matrix is often not necessary, enabling the use of matrix-free methods and avoiding the computationally expensive assembly of the stiffness matrix.
It can be proven that many iterative methods converge to the exact solution, and error bounds can be derived, e.g., for the conjugate gradient method (\CG{}\cite{Hestenes1952}) or the generalized minimal residual method ($\GMRES$ \cite{Saad1986}).
If the problem at hand is poorly conditioned, e.g., induced by high contrast in the local material properties or fine discretizations, numerical iterative solvers can suffer from slow convergence.
Therefore, preconditioners are often applied to improve the convergence behavior. In order to guarantee convergence for the preconditioned solver, these preconditioners have to fulfill certain requirements. At the same time, the efficient computational evaluation of the preconditioner is of paramount importance for the overall performance of the solver.

For homogenization problems, special iterative solution methods exist that are based on alternative discretizations such as the generalized methods of cells \cite{Aboudi2004} or the finite cell method \cite{Duester2008}.
{\color{revision}Particularly computationally efficient are iterative solvers based on the Fast Fourier Transform (FFT),} such as the Moulinec-Suquet scheme \cite{Moulinec1998} and many other methods that build up on it, as summarized in \cite{Schneider2021}.
While many of these specialized solvers do not fit directly into the framework of established iterative solvers such as the \CG{} method, some of these can act as preconditioners. For example, Fourier-Accelerated Nodal Solvers (FANS \cite{Leuschner2018}) combine a FEM discretization with FFT-based approaches and are applicable as preconditioners for Krylov subspace methods.

\subsection{Machine-learned surrogates for parametric PDEs}
\label{ssec:introduction-inexact}

While classical solvers for parametric PDEs can reliably provide the solution with almost arbitrary accuracy and are equipped with a rich mathematical theory, including convergence guarantees and error bounds, they suffer from high computational cost.
Unfortunately, many application scenarios, such as multi-scale simulations or materials design, require the solution for many different parameters, i.e., they are used in a many-query context.
For example, homogenization problems have to be solved for a large number of different microstructures. 
Often, such parametric variations lead to nontrivial changes of the solution, requiring unacceptable resources for the execution of simulations.

Due to these challenges, there has been an increase in research activities in the field of scientific machine learning (ML)\cite{Cuomo2022,Brunton2023,Watson2025}.
Therein, the goal is to approximate the results of classical simulations by machine-learned surrogates, which often allow a significantly lower computational effort in the evaluation.
The aim of scientific ML is to uncover the hidden relation between input parameters and quantities of interest, e.g., emerging from the solution of parametric PDEs.
Ideally, the ML prediction is a good approximation of the exact solution while being much faster to evaluate than classical solvers once the training is completed.
Unfortunately, there are seldom guarantees that the predicted solutions are accurate or physically meaningful.
Especially for parameter extrapolation outside of the training data set, the accuracy can be poor, while there are typically no error estimates to indicate so, at the same time.
Furthermore, machine-learned surrogates are often black-box models with limited interpretability of the learned weights of the model.

One class of machine-learned surrogates for PDEs is based on physics-informed neural networks (PINNs \cite{Raissi2019}). These learn the solution of a PDE on a given discretization. Unfortunately, PINNs are often unable to generalize parametric PDEs for parameters outside the training data without the need for additional retraining. This property limits their applicability for parametric homogenization problems.
More interesting for this application are neural operators, which aim at approximating a so-called solution operator\cite{Kovachki2023} for parametric PDEs and, thereby, learn a mapping between function spaces instead of relying on one fixed discretization.
Neural operators have been originally proposed in \cite{Li2020,Li2020MultipoleGN} and later formalized in the framework of parametric PDEs in \cite{Kovachki2023}.
Especially, Fourier Neural Operators (FNOs \cite{Li_fno_2021}) turned out to be powerful.
Extensions of these have been explored in \cite{li_physics-informed-no_2023,Zhao_incremental_2022}.
As an alternative to neural operators, DeepONets are also used as machine-learned surrogate models for parametric PDEs, but are based on a different architecture involving several branch nets and trunk nets\cite{Lu2021}.

\subsection{Machine learning enhanced hybrid solvers for parametric PDEs}
\label{ssec:introduction-hybrid-solvers}

On the one hand, machine-learned surrogates discussed in \cref{ssec:introduction-inexact} can rapidly provide approximate solutions for parametric PDEs, but without prior guarantees.  On the other hand, iterative solvers addressed in \cref{ssec:introduction-iterative-solvers} converge to the solution of discretized parametric PDEs with controllable accuracy, but often at the expense of many iterations.

In the present study, we focus on fusing iterative solvers with machine-learned surrogates to combine their undeniable advantages, namely the convergence guarantees of iterative schemes and the computational efficiency of machine learning. The resulting methods are referred to as \emph{hybrid solvers}.

Although this concept is still relatively new and unexplored, there are already some approaches in this direction in the literature.
These can be categorized into direct and indirect preconditioning approaches that are, e.g., compared in \cite{Kopanicakova2025}.
While in \emph{direct preconditioning} a machine-learned model itself is applied as a preconditioner in each iteration, in \emph{indirect preconditioning} only established algebraic preconditioners are constructed with the help of machine learning.
Often, algebraic multigrid (AMG) preconditioners are used for indirect preconditioning, as for example in \cite{Kopanicakova2025, Li2025}.
However, this usually involves a substantial additional computational overhead.

On the other hand, direct preconditioning is often more challenging to realize, as preconditioners for many iterative solvers have to fulfill certain properties to guarantee convergence.
For example, the \CG{} method requires the preconditioner to be a symmetric and positive definite operator.
Instead, direct machine-learned preconditioners are proposed for relaxation methods in \cite{Zhang2024}, and for the Generalized Minimal Residual ($\GMRES$) method in \cite{Xiang2022, Kopanicakova2025}, which is a Krylov subspace method that is also applicable to matrices that are neither symmetric nor positive definite \cite{Saad1986}.
On a different note, machine-learned preconditioners based on neural Green's operators have been recently presented, giving rise to an interesting interpretation as machine-learned Green's functions \cite{Melchers2024}.
However, this approach is also only applicable for the $\GMRES$ method.

Despite the reported progress on hybrid solvers mainly for the \GMRES{} method, the benefits of the oftentimes symmetric and positive definite systems found in most finite element discretizations cannot be exploited by the aforementioned hybrid solvers.
Contrary to that, the \CG{} method is far superior to the $\GMRES$ method in terms of both performance and memory consumption: it benefits from the symmetry, avoids the storage of the Krylov subspace, and has fewer hyperparameters.
Moreover, the \CG{} method is guaranteed to converge if the preconditioner is unconditionally symmetric and positive definite, a property that has not been hardwired into machine learning enhanced hybrid solvers to our knowledge.

The target of the present work is to extend previous efforts to enhance the performance of \CG{} solvers by using machine-learned preconditioners inspired by neural operators in a direct preconditioning approach with convergence guarantees.
A striking similarity of the resulting \UNOCG{} scheme with the established FFT-based solver \FANS{}\cite{Leuschner2018,Keshav2022} is uncovered, which draws a novel connection between very successful existing solvers and our new hybrid solver and gives rise to physical interpretation of the machine-learned preconditioner as machine-learned Green's functions similar to \cite{Melchers2024}.

\begin{figure}[ht]
	\centering
    \includegraphics[scale=1.0]{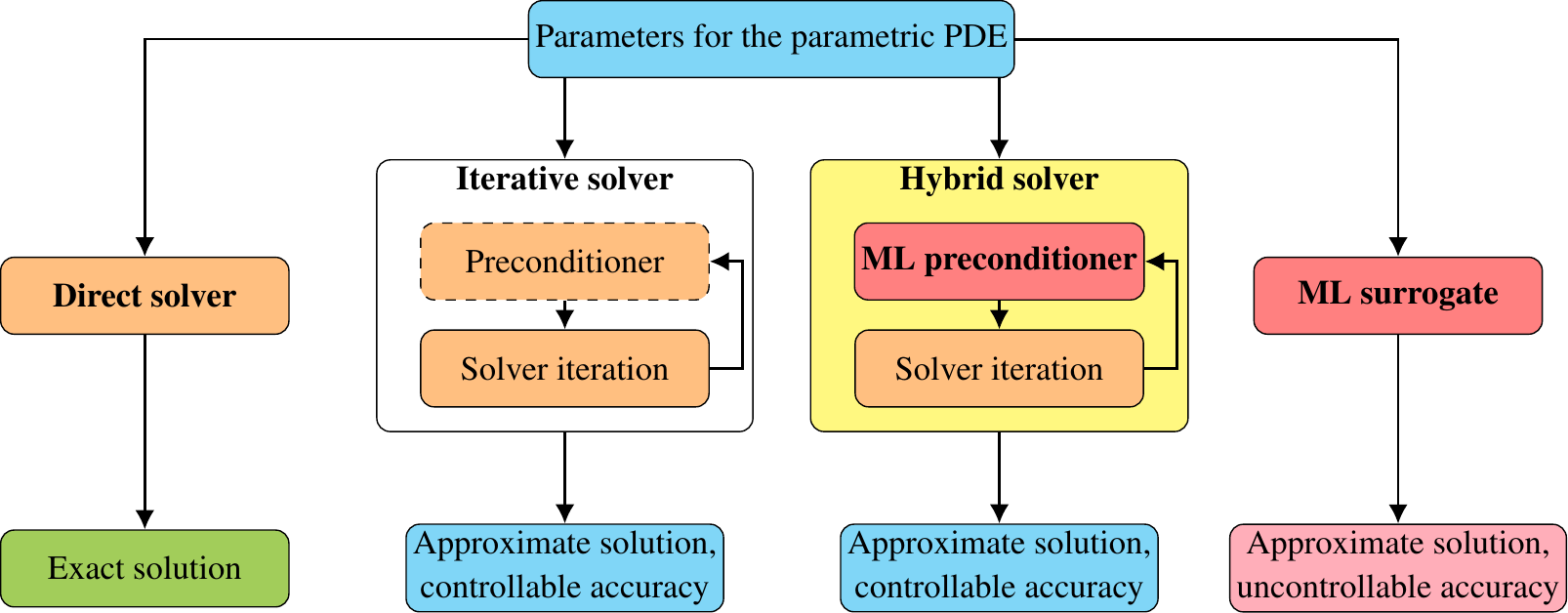}
	\caption{Different techniques for solving parametric PDEs, including a direct solver, iterative solver (\cref{ssec:introduction-iterative-solvers}), machine-learned surrogate (\cref{ssec:introduction-inexact}), and a hybrid solver with a machine-learned preconditioner (\cref{ssec:introduction-hybrid-solvers}).}
	\label{fig:solvers-overview}
\end{figure}

An overview of the different approaches described here for solving parametric PDEs, including a hybrid solver featuring a machine-learned preconditioner, is given in \cref{fig:solvers-overview}.

\subsection{Notation}
\label{ssec:notation}

All scalar-valued quantities such as~$\alpha \in \ffR$,~$\modeParam \in \ffN$ are written in regular font.
First-order tensors~$\flux$ and second-order tensors~$\heatCondTensor$ are written in bold font.
The~$\dimI$-th first-order unit tensor is denoted by~$\unitVecI \in \ffR^{n}$ for $1 \leq \dimI \leq n$.
Fourth-order tensors are denoted by~$\stiffTensor$.
Algebraic vectors are represented by an underline~$\ul{v} \in \ffR^n$ with~$n \in \ffN$, whereas matrices are characterized by two underlines, e.g.,~$\StiffMat \in \ffR^{n_1 \times n_2}$.
Arrays with three axes are written as~$\ulll{R} \in \ffR^{n_1 \times n_2 \times n_3}$, where $n_1, n_2, n_3 \in \ffN$.

Continuous vector-fields~$\res: \domain \to \ffR^{\nComp}$ with~$\nComp \in \ffN$ components (or nodal DOF) on a domain~$\domain \subsetneq \ffR^{\nDims}$ with~$\nDims \in \ffN$ can be discretized on a~$\nDims$-dimensional regular grid with~$n = N_1 \cdots N_\nDims$ nodes, where~$N_i \in \ffN$.
In this case, we denote the discrete field using boxes~$\field{ \res }: \ffN_0^\nDims \to \ffR^\nComp$.
Here, we use zero-based indexing for the discrete fields, i.e.,
\begin{align}
    \field{ \res }(\mode) = \field{ \res }(\modeComp_1, \dots, \modeComp_\nDims) \in \ffR^\nComp \,, && \text{for} && 0 \leq \modeComp_i \leq N_i - 1 \,, &&
    1 \leq i \leq \nDims \,, &&
    \mode \in \ffN_0^\nDims \,.
\end{align}
Throughout this article, fields (e.g., $\res$) are often represented both on a discrete grid~$\field{ \res }$ and in a vectorized representation~$\resD = \vecOp{\field{ \res }} \in \ffR^{\nComp N_1 \cdots N_\nDims}$ containing all degrees of freedom (DOF).
The mapping between these is given by indexing the discrete grid with discrete positions~$\mode\modeIdx \in \ffN^{\nDims}_0$ such that $\resSymbol_{\nComp(i-1)+j} = \field {\resSymbol_j} \left( \mode\modeIdx \right)$ for~$1 \leq \modeI \leq n$ and~$1 \leq j \leq \nComp$.
Further, the point-wise product or Hadamard product~$\odot$ between vectors~$\ul{a}$ and~$\ul{b}$ is defined as~$\left(\ul{a} \odot \ul{b}\right)_j = a_j b_j$.
The bullet point symbol~$\bullet$ is used as a placeholder for mathematical objects.
The identity map is written as $\mathrm{id}$.
The symmetric part of a matrix is denoted by $\symOp{ \bullet } = \frac{1}{2} \left( \bullet + \T{\bullet} \right)$, the space of symmetric matrices by $\symSpace{\ffR^{\nDims \times \nDims}}$, and the space of symmetric positive definite matrices by~$\spdSpace{\ffR^{\nDims \times \nDims}}$.
While the complex conjugate of a complex number, vector, or matrix is written as~$\Adj{\bullet}$, the conjugate transpose or Hermitian transpose is written as~$\Her{\bullet} = \Adj{\left(\T{\bullet}\right)}$.
Calligraphic symbols, e.g.,~$\cA$ are used for sets and spaces.
Operators that map between function spaces are written as~$\diffOp$, $\Prec$, for example.
The space of linear operators that map from a space~$\cA$ to a space~$\cB$ is denoted as~$\linopSpace(\cA;\cB)$.

\subsection{Problem setting}
\label{ssec:parPDE}

To formulate all covered methods in a common framework, we introduce a
linear parameter-dependent PDE (\ParPDE) that is described by a linear differential operator~$\diffOpParam$, as established in numerical analysis and also in scientific machine learning \cite{Kovachki2023}.
This operator may depend on parameters~$\param \in \paramSpace$,
where the set of admissible parameters is denoted as~$\paramSpace$.
The differential operator acts on the solution~$\solParam: \domain \to \ffR^\nComp$ of the parameter-dependent PDE with~$\nComp \in \ffN$ components and maps it to a right-hand side~$\rhsParam: \domain \to \ffR^\nComp$ that is given.
Both functions are defined on a bounded domain~$\domain \subsetneq \ffR^\nDims$ in~$\nDims \in \ffN$ dimensions.
In the following, the function space containing all parameter-dependent solutions~$\solParam$ with parameters~$\param \in \paramSpace$ is referred to as~$\solSpace$.
Usually, these function spaces are chosen to be spaces of square-integrable functions~$\SquareIntSpace(\bullet;\, \bullet)$ or subspaces thereof. For example, a possible function space for the solution is typically given by the Sobolev space~$\solSpace = \SobolevSpace(\domain;\, \ffR^\nComp) \subsetneq \SquareIntSpace(\domain;\, \ffR^\nComp)$ of weakly differentiable functions.
Based on that, the parameter-dependent differential operator~$\diffOpParam: \paramSpace \to \linopSpace\left( \solSpace; \solSpace\dual \right)$ is defined as a mapping from~$\paramSpace$ to the space of linear operators that themselves map the solution space~$\solSpace$ to its dual space~$\solSpace\dual$.
Then, the general formulation of a \ParPDE{} reads
\begin{align} \label{eq:parametric-PDE}
    \diffOpParam \, \solParam(\fx) = \rhsParam(\fx) \,, &&\fx \in \domain \,.
\end{align}
It is complemented by a set of suitable boundary conditions (BC) such as, e.g., homogeneous Dirichlet BC, i.e., ~$\onBoundary{ \solParam } = \sty{0}$, Neumann BC,
or periodic BC. Also, combinations of these are possible by partitioning the boundary $\partial \domain$.

Some \ParPDEs{} can be solved analytically using (in general matrix-valued) Green's functions~$\GreensParam: \domain \times \domain \to \ffR^{\nComp \times \nComp}$ that also depend on the parameters~$\param$.
Often,~$\GreensParam$ is also denoted as the PDE's fundamental solution.
An admissible Green's function fulfills the respective PDE's boundary conditions and at the same time satisfies the equation
\begin{align} \label{eq:fundamental-sol}
	\diffOpParam \, \GreensParam(\fx_0, \bullet) = \diracDelta_{\fx_0} \,, && \forall \; \fx_0 \in \domain \,,
\end{align}
where~$\diracDelta_{\fx_0}$ is the Dirac delta function (or distribution) that is centered at the position~$\fx_0 \in \domain$.
This Green's function can be interpreted similarly to an inversion of the differential operator~$\diffOpParam$, and with its help, the solution of the parametric PDE in \eqref{eq:parametric-PDE} can be obtained using the integral operator
\begin{equation} \label{eq:pde-integral-operator}
	\solParam(\fx) = \int_\domain \GreensParam(\fx, \fy) \, \rhsParam(\fy) \dInt{\fy} \,,
\end{equation}
since for a linear differential operator~$\diffOpParam$ and by applying \eqref{eq:fundamental-sol} it holds
\begin{align}
	\diffOpParam \, \solParam(\fx) = \diffOpParam \int_\domain \GreensParam(\fx, \fy) \, \rhsParam(\fy) \dInt{\fy}
	= \int_\domain \diracDelta_{\fx}(\fy) \, \rhsParam(\fy) \dInt{\fy} = \rhsParam(\fx) \,.
\end{align}
In the special case where the differential operator~$\diffOpParam$ is translation invariant \cite{Olver2014}, the Green's function can be taken to be a convolution kernel, i.e., there exists a function~$\GreensConvParam: \domain \to \ffR^{\nComp \times \nComp}$ such that~$\GreensParam(\fx, \fy) = \GreensConvParam(\fx - \fy)$, and thus \eqref{eq:pde-integral-operator} can be formulated as the convolution
\begin{equation} \label{eq:pde-convolution}
	\solParam(\fx) = \GreensConvParam \conv \, \rhsParam = \int_\domain \GreensConvParam(\fx - \fy) \, \rhsParam(\fy) \dInt{\fy} \,.
\end{equation}
However, the approach of directly solving \ParPDEs{} using a convolution with a Green's function as in~\eqref{eq:pde-convolution} is typically only straightforward for trivial cases, for instance, if the parameters~$\param$ are constant over the domain, and if periodic BC are considered at the same time.

Many homogenization problems, as motivated in \cref{ssec:introduction-motivation}, can be described in this framework of linear \ParPDEs{} if a linear material law is considered.
For many of the PDEs underlying homogenization problems, a Green's function can be derived for spatially homogeneous parameters, i.e., a homogeneous material with constant material parameters.
However, the difficulty lies in the heterogeneity of the material parameters due to microstructures with different phases.
Solving these parametric PDEs is particularly challenging when there is a high contrast between the phases.
Thus, it is not possible to apply~\eqref{eq:pde-convolution} directly for the presented homogenization problems, where the material parameters are heterogeneous, i.e., differ drastically in the different phases.
Nevertheless, we will see that the concept of Green's functions will turn out to be useful in an indirect manner for solving homogenization problems multiple times throughout this work, using classical solvers and also using our proposed hybrid solver.

\section{Machine-learned surrogates for parametric PDEs}
\label{sec:mlpde}

\subsection{Physics-Informed Neural Networks (PINNs)}
\label{ssec:pinn}

A first approach to solving \ParPDEs{} using machine learning is given by Physics-Informed Neural Networks (PINNs).
Originally proposed in \cite{Raissi2019}, PINNs aim at approximating the PDE's solution for fixed parameters~$\param \in \paramSpace$ instead of approximating a solution operator that maps from the parameters~$\param \in \paramSpace$ to the solution~$\solParam \in \solSpace$ of the parameter-dependent PDE.
Thus, PINNs only have limited generalization capabilities for parametric PDEs, in general.
As soon as the parameters~$\param$, the right-hand side~$\rhsParam$, or the domain~$\domain$ changes, PINNs usually have to be retrained, making their application unsuitable for \ParPDEs{} as introduced in \cref{ssec:parPDE} and, thus, for all homogenization problems covered in \cref{sec:problems}.

\subsection{Neural Operators}
\label{ssec:neuralop}

Recently, neural operators (NOs) \cite{Li2020,Kovachki2023} started to emerge as a new class of machine-learned surrogates that widely circumvent the limitations of Physics-Informed Neural Networks as they try to learn the mapping from the parameters~$\param \in \paramSpace$ to the corresponding solution of the parametric PDE~$\solParam$ directly instead of learning the solution~$\solParam$ for a single instance of the \ParPDEs{} with fixed parameters~$\param \in \paramSpace$.
Originally, neural operators were proposed in \cite{Li2020}.
An extensive theoretical framework, including an overview of different realizations of neural operators and featuring a universal approximation theorem for nonlinear differential operators, is available in \cite{Kovachki2023}.
Because of the generalization capabilities of neural operators that include extrapolation to previously unseen parameters~$\param \in \paramSpace$ during inference after the model is trained, neural operators can be considered a more powerful tool than PINNs in this regard.

In contrast to many classical approaches in scientific machine learning, neural operators do not just aim at learning a mapping between finite-dimensional vector spaces, for instance, a mapping from the discretized parameters on a grid to the discretized solution field on the same grid.
Instead, neural operators try to learn a mapping from the space of the parameters~$\paramSpace$ to the function space of the PDE's solution~$\solSpace$. Thereby, trained neural operators can be potentially invariant to the discretization of the PDE \cite{Kovachki2023}.
For a parametric PDE as defined in \eqref{eq:parametric-PDE}, a solution operator~$\solOp$ that maps the parameters~$\param \in \paramSpace$ to the corresponding solution~$\solParam \in \solSpace$ of the parametric PDE can be defined under the assumption that an inverse operator~$\diffOp^{-1}_{\param}$ of the differential operator $\diffOpParam$ exists as
\begin{align}
	\solOp: \paramSpace \to \solSpace \,, && \param \mapsto \diffOp^{-1}_{\param} \, \rhsParam = \solParam \,.
\end{align}
In fact, this is the operator that neural operators try to approximate for linear parametric PDEs \cite{Kovachki2023}. By that, neural operators have the possibility to generalize to new parameters~$\param \in \paramSpace$ not part of the training data.

\subsubsection{Architecture}
\label{sssec:architecture}

First, neural operators apply a learnable lifting operator~$\fnoProjOp$ from the parameters~$\param \in \paramSpace$ and potentially further arguments to a vector field~$\fnoFeatures\layerIdxFirst: \domain \to \ffR^{\nChannelsLayerFirst}$,~$\nChannelsLayerFirst \in \ffN$, consisting of artificial features.
This lifting operator could be realized using an ANN that is applied point-wise to the input field. The operator is defined as
\begin{align} \label{eq:neuralop-projection}
	\fnoProjOp: \paramSpace \times \dots \to \SquareIntSpace(\domain;\,\ffR^{\nChannelsLayerFirst}) \,, &&\fnoFeatures\layerIdxFirst(\fx) = \fnoProjOp \left( \param, \dots \right) \,.
\end{align}
Neural operators aim at approximating the operator~$\solOp: \paramSpace \to \solSpace$ by learning matrix-valued integral kernels~$\fnoKernel: \domain \times \domain \to \ffR^{\nChannelsLayer \times \nChannelsLayerPrev}$ in multiple neural operator (NO) layers that map from one feature field~$\fnoFeatures\layerIdxPrev(\fx)$ with~$\nChannelsLayerPrev \in \ffN$ channels to another feature field~$\fnoFeatures\layerIdx(\fx)$ with~$\nChannelsLayer \in \ffN$ channels.
The superscript~$\bullet\layerIdx$ expresses the association with the~$\layerI$-th layer of the model, where~$1 \leq \layerI \leq \nLayers$ , and~$\nLayers \in \ffN$.
In addition, a linear operator in the form of the matrix~$\fnoBypass \in \ffR^{\nChannelsLayer \times \nChannelsLayerPrev}$ based on learnable weights is applied as a bypass (or shortcut) in parallel to the integral operator in each layer.
Finally, each layer features a nonlinear activation function~$\fnoAct: \ffR \to \ffR$ that is applied point-wise at each spatial position~$\fx$ to each feature channel.
Hence, the full mathematical description of the $\layerI$-th neural operator layer reads
\begin{equation} \label{eq:neuralop-integral-operator}
	\fnoFeatures\layerIdx(\fx) = \fnoAct \left( \fnoBypass \, \fnoFeatures\layerIdxPrev(\fx) + \int_{\domain} \fnoKernel(\fx, \fy) \, \fnoFeatures\layerIdxPrev(\fy) \,\dInt{\fy} \right) \,.
\end{equation}
For~$\nChannelsLayerPrev=\nChannelsLayer=\nComp$, i.e., if the number of input and output channels match the number of physical components, the integral kernel~$\fnoKernel$ is similar to the Green's function that is used in the integral operator~\eqref{eq:pde-integral-operator} to obtain the solution of linear PDEs.
However, the difference to an integral operator involving a classical Green's function is that the integral operator in~\eqref{eq:neuralop-integral-operator} does not necessarily act on the right-hand side of the PDE~$\rhsParam$ but on an arbitrary function that is the result of the previous NO layer.
Moreover, the number of channels~$\nChannels_\layerI$ in each layer is usually chosen to be higher than the number of components of the PDE's solution~$\nComp$.
Herein lies the approximation power of neural operators, similar to conventional ANNs that utilize hidden layers---with each having more artificial features than physical input or output features.

\begin{figure}[ht]
	\centering
    \includegraphics[scale=1.0]{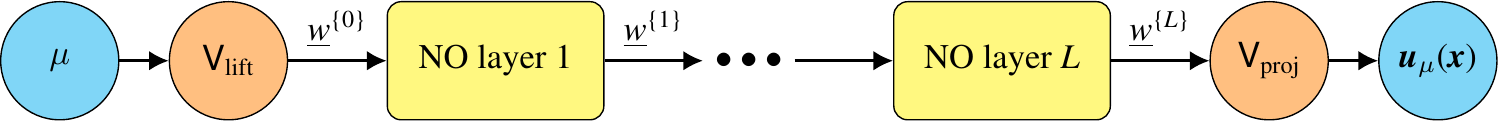}
	\caption{Architecture of a neural operator model with a lifting operator~$\fnoProjOp$, several neural operator (NO) layers, and a projection operator~$\fnoBackprojOp$. Figure inspired by \cite{Li_fno_2021}.}
	\label{fig:neuralop-architecture}
\end{figure}

Finally, after~$\layerILast$ neural operator (NO) layers, a projection of the resulting vector field~$\fnoFeatures\layerIdxLast$ consisting of artificial features to an approximation of the solution is performed. This is described by the operator
\begin{align} \label{eq:neuralop-backprojection}
	\fnoBackprojOp: \SquareIntSpace(\domain;\,\ffR^{c_L}) \to \solSpace \,, &&
	\solParam(\fx) = \fnoBackprojOp\left( \fnoFeatures\layerIdxLast(\fx) \right) \,, && \fx \in \domain \,,
\end{align}
which---as the lifting operator~$\fnoProjOp$---should be a point-wise operation in order to guarantee that the learned neural operator is invariant to the resolution used for the discretization of the parameters~$\param \in \paramSpace$ and the solution~$\solParam \in \solSpace$.
An overview of the general structure is given in \cref{fig:neuralop-architecture}.
{\color{revision}
While neural operators typically consist of multiple layers with nonlinear activations, linear architectures with a single neural operator layer and identity activation~$\sigma=\mathrm{id}$ are also used as special case for linear PDEs \cite{Kovachki2023,Li2020}.}

In the general framework of neural operators, each layer's integral kernel is a function~$\fnoKernel \in \SquareIntSpace(\domain \times \domain; \ffR^{\nChannelsLayer \times \nChannelsLayerPrev})$.
In order to learn an approximation for such a function, it has to be somehow parametrized.
For that, different approaches are proposed that are summarized in \cite{Kovachki2023} while the most popular ones are Graph Neural Operators \cite{Li2020} and Fourier Neural Operators (FNOs \cite{Li_fno_2021}).

\subsubsection{Fourier Neural Operators}
\label{sssec:fno}

\begin{figure}[ht]
	\centering
    \includegraphics[scale=1.0]{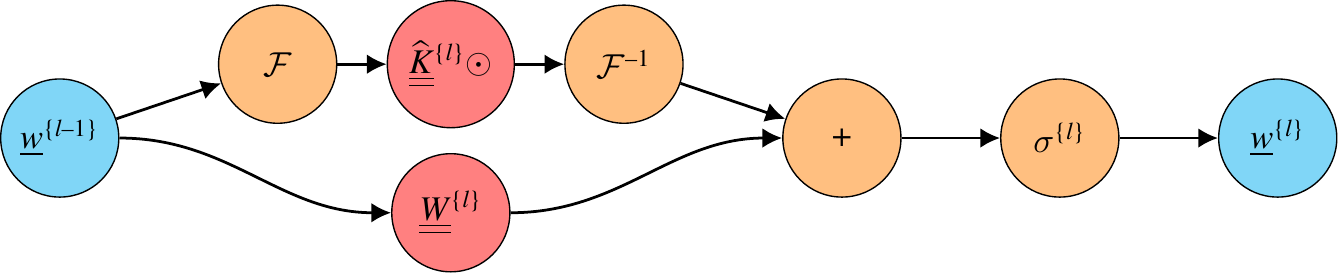}
	\caption{Architecture of a Fourier layer based on the Fourier transform $\trafoFFT$ that features a learnable kernel~$\fnoKernel$ in Fourier space and a learnable bypass~$\fnoBypass$. Figure inspired by \cite{Li_fno_2021}.}
	\label{fig:fno-architecture}
\end{figure}

In Fourier Neural Operators (FNOs), the integral kernel~$\fnoKernel$ is parameterized in the Fourier space.
To this end, it is assumed that a function~$\fnoConvKernel \in \SquareIntSpace(\domain; \ffR^{\nChannelsLayer \times \nChannelsLayerPrev})$ exists, such that~$\fnoKernel(\fx, \fy) = \fnoConvKernel(\fx - \fy)$.
Based on that, the application of the integral operator can be reformulated as the convolution operation
\begin{align} \label{eq:fno-integral-operator}
	\int_\domain \fnoKernel(\fx, \fy) \, \fnoFeatures\layerIdxPrev(\fy) \,\dInt{\fy} = \int_\domain \fnoConvKernel(\fx - \fy) \, \fnoFeatures\layerIdxPrev(\fy) \,\dInt{\fy} = \left( \fnoConvKernel \conv \fnoFeatures\layerIdxPrev \right)(\fx) \,.
\end{align}
By virtue of the convolution theorem, this can be computed {\color{revision}at $\cO(n \log n)$ complexity} using the FFT for $1 \leq i \leq \nChannelsLayer$ via
\begin{align} \label{eq:fno-convolution}
	\left( \fnoConvKernel \conv \fnoFeatures\layerIdxPrev \right)_i =
	\sum_{j=1}^{\nChannelsLayerPrev} \overline{K}\layerIdx_{ij} \conv v\layerIdxPrev_j = \trafoFFTInv \left(\sum_{j=1}^{\nChannelsLayerPrev} \widehat{K}\layerIdx_{ij} \odot \trafoFFT \left( v\layerIdxPrev_j \right) \right) \,, &&
	\widehat{K}\layerIdx_{ij} = \trafoFFT\left( \overline{K}\layerIdx_{ij} \right) \,.
\end{align}
Overall, a single neural operator layer in FNOs (also called Fourier layer) is given by
\begin{equation} \label{eq:fno-operator}
	\fnoFeatures\layerIdx(\fx) = \fnoAct \left( \fnoBypass \, \fnoFeatures\layerIdxPrev(\fx) + \trafoFFTInv \left( \fnoConvKernelFourier \odot \trafoFFT \left( \fnoFeatures\layerIdxPrev(\fx) \right) \right) \right) \,,
\end{equation}
where we use the operator $\odot$ as shorthand notation for the matrix-valued Hadamard product with~$\fnoConvKernelFourier$ from~\eqref{eq:fno-convolution}.
The resulting architecture of a Fourier layer is visualized in \cref{fig:fno-architecture}.

The ansatz for FNOs is to parametrize~$\fnoConvKernel \in \SquareIntSpace(\domain; \ffR^{\nChannelsLayer \times \nChannelsLayerPrev})$ as a periodic function that can be represented by a truncated~$\nDims$-dimensional Fourier series using a finite set of discrete frequency modes~$\mode\modeIdx \in \ffZ^{\nDims}$.
The common approach is to select a number of $\nModes$ low-frequency modes as part of the finite set of modes
\begin{align}
	\modeSet = \left\{ \mode\modeIdx: \; 1 \leq \modeI \leq \nModes \right\} \subsetneq \ffZ^\nDims \,,
	&& \nModes \in \ffN \,,
\end{align}
that should be included in the parametrization of $\fnoConvKernel$.
Note that we choose~$\ffZ^\nDims$ instead of~$\ffN^\nDims_0$ to also include the corresponding negative frequencies in Fourier space.
The discrete function~$\fnoConvKernelFourier: \ffZ^\nDims \to \ffC^{\nChannelsLayer \times \nChannelsLayerPrev}$ in Fourier space from~\eqref{eq:fno-operator} is then parametrized using the third-order weight tensor
\begin{align}
    \fnoWeightTensor \in \ffC^{\nModes \times \nChannelsLayer \times \nChannelsLayerPrev} \,, &&
    1 \leq \layerI \leq \nLayers \,, &&
    \nModes = \left| \modeSet \right| \,.
\end{align}
In most applications, the solution $\solParam$ of the parametric PDE is real-valued instead of complex-valued.
This is also the case for the elastic and thermal homogenization problems covered in \cref{sec:problems}.
To enforce~$\fnoFeatures\layerIdxNext(\fx)$ to be real-valued, the parametrized kernel in Fourier space~$\fnoConvKernelFourier$ has to fulfill the property of conjugate symmetry, that is,
\begin{align} \label{eq:fno-complex-conjugate}
	\fnoConvKernelFourier \left( \conjOp{\mode} \right) = \left( \fnoConvKernelFourier \left( \mode \right) \right)^\ast \,, &&
	\conjOp{\mode} = \T{\begin{bmatrix} \modeComp_1 & \cdots & \modeComp_{\nDims - 1} & {-\modeComp_\nDims} \end{bmatrix}} \,, &&
	\forall \mode \in \ffZ^\nDims \,.
\end{align}
This can be achieved by learning weights for the integral kernel~$\fnoConvKernelFourier$ only on a half-space by choosing the set of frequency modes such that~$\modeSet \subsetneq \ffZ^{\nDims - 1} \times \ffN_0$ instead of~$\modeSet \subsetneq \ffZ^\nDims$.
The set of conjugate frequency modes corresponding to $\modeSet$ is defined as
\begin{align} \label{eq:modesetconj-fno}
	\modeSetConj = \left\{ \conjOp{\mode\modeIdx} : \mode\modeIdx \in \modeSet \right\} \setminus \modeSet \,.
\end{align}
Then, the function $\fnoConvKernelFourier$ is defined based on the learnable weights $\fnoWeightTensor \in \ffC^{\nModes \times \nChannelsLayer \times \nChannelsLayerPrev}$ for~$1 \leq \layerI \leq \nLayers$ as
\begin{align} \label{eq:fno-def-param}
	\fnoConvKernelFourier \left( \mode \right) =
	\begin{cases}
		\fnoWeightSymbol\layerIdx_{i \bullet \bullet} & \text{if} \quad \exists i : \, \mode = \mode\modeIdx \in \modeSet \,,\\
		\left( \fnoWeightSymbol\layerIdx_{i \bullet \bullet} \right)^* & \text{if} \quad \exists i : \, \mode = \conjOp{\mode\modeIdx} \in \modeSetConj \,,\\
		\ull{0} & \text{else} \,,
	\end{cases}
	&& \forall \mode \in \ffZ^\nDims \,.
\end{align}
Effectively, this means that the attributions of all frequency modes that are not in~$\modeSet$ or~$\modeSetConj$ are set to zero.
For simplicity, a hyperparameter~$\modeParam \in \ffN$ is used to determine the desired number of modes in the sets~$\modeSet$ and~$\modeSetConj$ in a convenient way.
For the mode selection, we introduce a notation of integer intervals as~$\discrIntv{a}{b} = \left[a,b\right] \cap \ffZ$.
In particular, the set $\discrIntv{a}{b}$ contains~$a, b \in \ffZ$.
Then, the sets~$\modeSet$ and~$\modeSetConj$ can be, for example\footnote{Note that other choices are also possible, e.g., $\modeSet = \discrIntv{-m_1}{m_1} \times \cdots \times \discrIntv{-m_{\nDims-1}}{m_{\nDims-1}} \times \discrIntv{0}{m_{\nDims}}$ and $\modeSetConj = \discrIntv{{-m_1}}{m_1} \times \cdots \times \discrIntv{-m_{\nDims-1}}{m_{\nDims-1}} \times \discrIntv{-m_{\nDims}}{{-1}}$ with $m_1,\dots,m_{\nDims} \in \ffN$.}, defined as
\begin{align}
	\modeSet = \discrIntv{{-M}}{M}^{\nDims - 1} \times \discrIntv{0}{M} \,, &&
	\modeSetConj = \discrIntv{{-M}}{M}^{\nDims - 1} \times \discrIntv{{-M}}{{-1}} \,.
\end{align}

The property of complex conjugateness in~\eqref{eq:fno-complex-conjugate} can be automatically fulfilled if the real FFT of a real-valued input~$\fnoFeatures\layerIdx$ is used in the Fourier layer and only the learnable weights related to the parametrized modes~$\modeSet$ are applied in the Fourier space.
As it is common in software implementations\cite{FFTW05}, we assume that a~$\nDims$-dimensional real FFT is computed by performing a one-dimensional real FFT along the last (i.e.,~$\nDims$-th dimension) and full FFTs along all other dimensions (i.e.,~$1 \leq \dimI \leq \nDims-1$).

\section{Iterative solvers for parametric PDEs}
\label{sec:solvers}

The machine-learned surrogates presented in \cref{sec:mlpde} can hypothetically predict the solution of parametric PDEs with high performance.
However, these are difficult to interpret black-box models, and there are no guarantees for the accuracy of the predictions, which makes their use unacceptable in many scientific applications.
This is especially true for the homogenization problems introduced in \cref{sec:problems}, as their use in multi-scale simulations requires highly accurate and physically consistent solutions of the underlying parametric PDEs.

In contrast, iterative solvers for PDEs follow a different paradigm.
In each iteration, the residual of the discretized PDE is evaluated to correct the approximated solution with different algorithms.
If implemented cleverly, for example, based on sparse linear algebra or even as matrix-free methods, a single iteration can be executed significantly faster than direct solvers and often even faster than machine-learned surrogates.
If only a few iterations are required to reach the desired accuracy, {\color{revision}the resulting overall runtime enables the fast solution of parametric PDEs}.
At the same time, convergence can often be proven, and there are a priori error bounds for the approximated solution in each iteration.
In this work, we only cover linear PDEs and iterative solvers for the corresponding linear algebraic systems.
However, many of them can be generalized to nonlinear problems, e.g., using a line search method.

\subsection{Discretization}
\label{ssec:discretization}

Here, we focus on the finite element method (FEM) often used in computational mechanics, while using other discretization methods, e.g., the finite difference method or the finite volume method, is also possible for \ParPDEs{}.
In the end, our presented framework of hybrid solvers with machine-learned preconditioners is independent of the discretization scheme used, and it can be applied to different discretizations.
We restrict ourselves to a regular grid for various reasons:
\begin{itemize}
	\item Parameters for \ParPDEs{} are often defined on regular grids. For example, microstructures in homogenization are often given by images that stem from CT scans or synthetic methods.
	\item Fast solvers for homogenization problems and other \ParPDEs{} can take advantage of the Fast Fourier Transform (FFT) and, hence, require a discretization based on regular grids.
	\item Most approaches in scientific machine learning, including FNOs, operate only on regular grids.
	\item A discretization with regular grids enables several technical performance optimizations since all grid cells are the same up to translation, e.g., the use of matrix-free implementations.
\end{itemize}

\begin{figure}[ht]
    \begin{subfigure}[b]{0.32\textwidth}
        \begin{center}
    \includegraphics[scale=1.0]{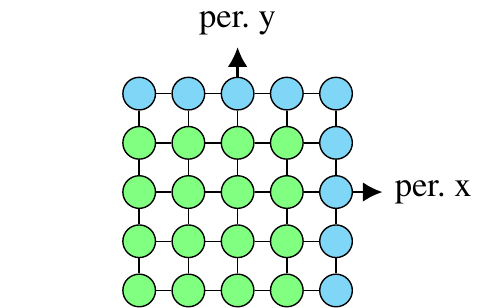}
        \end{center}
        \vspace{-1.0\baselineskip}
        \caption{Periodic BC.}
    \end{subfigure}
\begin{subfigure}[b]{0.32\textwidth}
        \begin{center}
\includegraphics[scale=1.0]{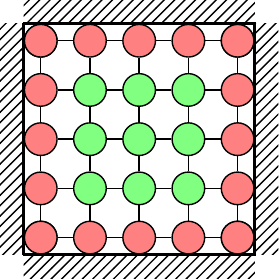}
        \end{center}
        \vspace{-1.0\baselineskip}
        \caption{Dirichlet BC.}
    \end{subfigure}
\begin{subfigure}[b]{0.32\textwidth}
        \begin{center}
\includegraphics[scale=1.0]{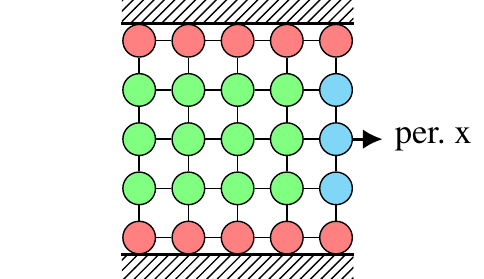}
        \end{center}
        \vspace{-1.0\baselineskip}
        \caption{Mixed BC (per. $\mathrm{x}$, Dir. $\mathrm{y}$).}
    \end{subfigure}
	\caption{Figures of a regular grid in $\nDims=2$ dimensions with free nodes (green), periodic nodes (blue), and fixed nodes (red) for periodic BC (left), Dirichlet BC (middle), and mixed BC (right).}
	\label{fig:regular-grid}
\end{figure}

The considered~$\nDims$-dimensional regular grid consists of~$\nElem = N_1 \cdots N_\nDims$ elements in form of rectangles ($\nDims=2$) or bricks ($\nDims=3$).
Thus, the grid contains in total~$\nNodes = (N_1 + 1) \cdots (N_\nDims + 1)$ nodes.
The degrees of freedom (DOF) of the discretization are associated with a subset of these {\color{revision}nodes}---called free nodes.
The number of free nodes will be denoted with~$n$ and depends on the boundary conditions.
When using periodic BC on all edges, there are~$n_{\rm per} = N_1 \cdots N_\nDims$ free nodes.
However, when using Dirichlet BC on all edges, only~$n_{\rm Dir} = (N_1 - 1) \cdots (N_\nDims - 1)$ free nodes remain.
Depending on the problem, one or multiple nodal DOF~$\nComp$ are assigned per node.
Hence, the discretized problem includes~$\nDof = \nComp n$ DOF in total.
As an example, regular grids with different boundary conditions are shown in \cref{fig:regular-grid}.

In the following, the \ParPDE{} from~\eqref{eq:parametric-PDE} is discretized using the Galerkin FEM, which involves deriving a variational form for the \ParPDE{}, which leads to
\begin{align}
\label{eq:variational-form}
    \bilinForm (\solParam, \test; \param) = \linForm (\test; \param) && \forall \test \in \solSpace \,, && \param \in \paramSpace \,,
\end{align}
with a parameter-dependent bilinear form~$\bilinForm(\bullet, \bullet; \param)$ and a parameter-dependent linear functional~$\linForm(\bullet; \param)$.
In the FEM, the domain~$\domain$ is divided into a finite number~$\nElem \in \ffN$ of elements~$\domain^{e}$, where~$1 \leq e \leq \nElem$. Based on the resulting mesh, an approximate solution~$\sol^{\rm h}_{\param} \in \solSpaceD$ in a finite-dimensional function space~$\solSpaceD \subsetneq \solSpace$ is determined, such that the variational form
\begin{align} \label{eq:fem-weak-form-discr}
	\bilinForm \left( \sol^{\rm h}_{\param}, \test^{\rm h}; \param \right) = \linForm \left( \test^{\rm h}; \param \right) &&
	\forall \test^{\rm h} \in \solSpaceD \,,
\end{align}
{\color{revision}is fulfilled for all test functions~$\test^{\rm h} \in \solSpaceD$.}
The discrete solution space $\solSpaceD$ is spanned by the FE ansatz functions $\femBasis_i$, leading to the algebraic system ($1 \leq i, j \leq cn$)
\begin{align} \label{eq:fem-linear-system}
	\StiffMatParam \, \solDParam = \rhsDParam \,, && \param \in \paramSpace  \,, &&
 	\left(\StiffMatParam\right)_{ij} = \bilinForm \left( \femBasis_j, \femBasis_i; \param \right) \,, &&
 	\left( \rhsDParam \right)_i = \linForm(\femBasis_i; \param) \,.
\end{align}
Herein,~$\StiffMatParam\in \ffR^{\nComp n \times \nComp n}$ denotes the FEM stiffness matrix and~$\rhsDParam \in \ffR^{\nComp n}$ the right-hand-side vector, which both depend on the parameter $\param$.
The vector $\solDParam \in \ffR^{\nComp n}$ contains the DOF of the solution $\solParam$ in the discrete solution space $\solSpaceD$.
For elliptic problems, the bilinear form~$\bilinForm$ is symmetric and coercive, implying symmetry and positive definiteness of~$\StiffMatParam$.

\subsection{Conjugate Gradient (\texorpdfstring{\CG{}}{CG}) method}
\label{ssec:cg}

The conjugate gradient (\CG{}) method is a well-known iterative solver to compute the solution of a linear system
\begin{align} \label{eq:cg-system}
	\StiffMat\, \solD = \rhsD \,, && \text{with} && \StiffMat &\in \ffR^{\nComp n \times \nComp n} \,, & \rhsD & \in \ffR^{\nComp n} \,, & \solD & \in \ffR^{\nComp n} \,,
\end{align}
where~$\StiffMat$ is a symmetric and positive definite matrix, often emerging from a FEM discretization.
The \CG{} method belongs to the class of Krylov subspace methods, where approximate solutions~$\solD\cgIdx \in \ffR^{\nComp n}$ are constructed iteratively based on certain optimality principles starting from an initial guess~$\solD\cgIdxFirst \in \ffR^{\nComp n}$.
Each iteration~$\cgI=1, 2, \dots$ involves the computation of a search direction~$\cgSearchDir \in \ffR^{\nComp n}$ and the corresponding optimal step width~$\cgStepWidth \in \ffR$.
For that, the action of one matrix-vector product with the matrix~$\StiffMat$ has to be evaluated per iteration at a computational effort in~$\orderOf{n}$ {\color{revision}provided that~$\StiffMat$ is sparse}.

One of the problems with standard iterative solvers such as the \CG{} method is that they suffer from slow convergence for poorly conditioned systems, e.g., induced by finely resolved discretizations, i.e., they require many iterations to converge.
In order to lower the number of required iterations, the \CG{} method is usually used in combination with a problem-specific preconditioner.
Here, we consider the approach of left-preconditioning\cite{Barrett1994} using a preconditioner matrix~$\PrecD \in \ffR^{\nComp n \times \nComp n}$.
This leads to the preconditioned linear system
\begin{equation} \label{eq:prec-linear-system}
	\PrecD \, \StiffMat\, \solD = \PrecD \; \rhsD \,.
\end{equation}
Using a suitable preconditioner can lead to a significantly faster convergence.
In the case of the \CG{} method, the preconditioner matrix~$\PrecD$ has to be symmetric and positive definite. Then the \CG{} solver applied to \eqref{eq:prec-linear-system} converges to the solution of the original system.
Note that the trivial preconditioner $\PrecD = \IdentD$ recovers the \CG{} method without preconditioner (unpreconditioned \CG{}).

A simple algebraic preconditioner is the Jacobi preconditioner defined as
\begin{align}
	\PrecD_\Jac \in \ffR^{\nComp n \times \nComp n} \,, && \left( \PrecD_\Jac \right)_{ij} =
    \begin{cases}
		1 \, / \, A_{ij} & \text{if} \quad i = j \,,\\
		0 & \text{else} \,,
	\end{cases} && 1 \leq i, j \leq \nComp n \,.
\end{align}
In the context of the current work $\StiffMat=\StiffMatParam$, $\solD=\solDParam$, $\rhsD=\rhsDParam$ holds. 
This implies that the preconditioner~$\color{revision}\PrecD_\Jac$ depends explicitly on the parameter-dependent stiffness matrix~$\StiffMatParam$ {\color{revision}that itself is different for each parameter realization~$\param \in \paramSpace$. 
In the special case of using a Jacobi preconditioner for problems on regular grids, the full stiffness matrix~$\StiffMatParam$ does not have to be assembled, but more sophisticated algebraic preconditioners (e.g., incomplete Cholesky \cite{Saad2003}) require direct access to the assembled matrix.}
The complete algorithmic description of the preconditioned conjugate gradient (\CG{}) method is available in \cref{alg:pcg}.

\begin{algorithm}
    \caption{Preconditioned conjugate gradient (\CG{}) method}
\label{alg:pcg}
    \hspace*{\algorithmicindent} \textbf{Data:} $\StiffMat\in \ffR^{\nComp n \times \nComp n}$, $\rhsDParam \in \ffR^{\nComp n}$, $\solD\cgIdxFirst_{\param} \in \ffR^{\nComp n}$, $\PrecD \in \ffR^{\nComp n \times \nComp n}$, $\cgTol > 0$ \\
    \hspace*{\algorithmicindent} \textbf{Result:} $\solDParam \in \ffR^{\nComp n}$
    \begin{algorithmic}[1]
\State{$\resD \gets \rhsDParam - \StiffMat\, \solD\cgIdxFirst_{\param}$} \algorithmiccomment{Initial residual}
\State{$\cgSearchDir \gets \ul{0} \in \ffR^{\nComp n}$; $\cgDelta_1 \gets 1 \in \ffR$} \algorithmiccomment{Initialization}
	\While{ $\normg{\resD}_\infty > \cgTol $ }
		\State{$\precResD \gets \PrecD \, \resD$} \algorithmiccomment{Apply preconditioner}
		\State{$\cgDelta_0 \gets \cgDelta_1; \; \cgDelta_1 \gets \resD \cdot \precResD$}
		\State{$\cgSearchDir \gets \precResD + \frac{\cgDelta_1 }{\cgDelta_0} \cgSearchDir$} \algorithmiccomment{Update search direction}
		\State{$\cgMatvecResult \gets \StiffMat\, \cgSearchDir$} \algorithmiccomment{Matrix-vector product}
		\State{$\cgStepWidth \gets \cgDelta_1 \, / \left(\cgSearchDir \cdot \cgMatvecResult \right)$} \algorithmiccomment{Compute optimal step width}
		\State{$\resD \gets \resD - \cgStepWidth \cgMatvecResult$} \algorithmiccomment{Update residual}
		\State{$\solDParam \gets \solDParam + \cgStepWidth \cgSearchDir$} \algorithmiccomment{Update solution approximation}
    \EndWhile
    \end{algorithmic}
\end{algorithm}

It can be proven\cite{Hackbusch1994} that when performed in exact arithmetics, the unpreconditioned \CG{} method (i.e., \cref{alg:pcg} with $\PrecD = \ull{I}$) converges to the solution $\solDParam$ after at most~$\nComp n$ iterations.
When the error in the~$\cgI$-th iteration is measured in the energy norm as
\begin{align} \label{eq:energy-norm}
    \normg{\cgError\cgIdx}_{\StiffMat} = \sqrt{ \T{\cgError\cgIdx} \StiffMat\, \cgError\cgIdx } \,, &&
	\cgError\cgIdx = \solD\cgIdx_{\param} - \solDParam \,,
\end{align}
an error bound can be derived\cite{Hackbusch1994} that reads
\begin{align} \label{eq:cg-convergence}
	\normg{\cgError\cgIdx}_{\StiffMat} \leq 2 \cgRate^\cgI \normg{\cgError\cgIdxFirst}_{\StiffMat}
    \,, &&
    \cgRate = \frac{\sqrt{\cond(\StiffMat)} - 1}{\sqrt{\cond(\StiffMat)} + 1} \,, &&
    \cond(\StiffMat) = \frac{ \norm{\eigMax(\StiffMat)} }{ \norm{\eigMin(\StiffMat)} } \,,
\end{align}
where~$\cgRate$ is the convergence rate of the \CG{} method and~$\cond(\StiffMat)$ refers to the condition number of the stiffness matrix~$\StiffMat$, which is defined based on its extremal eigenvalues~$\eigMax$ and~$\eigMin$.

If the preconditioner matrix~$\PrecD$ is symmetric and positive definite, a Cholesky decomposition~$\PrecD = \ull{L} \, \T{\ull{L}}$ exists and the preconditioned \CG{} method is equivalent to the unpreconditioned \CG{} method for the modified matrix $\widetilde{\StiffMat}_{\param} = \ull{L} \, \StiffMat \, \T{\ull{L}}$ that is also symmetric and positive definite, and has the same eigenvalues as~$\PrecD \, \StiffMat$ as shown in \cite{Golub2013}.
Thus, the convergence theory for unpreconditioned \CG{} can also be applied to preconditioned \CG{} as long as~$\PrecD$ is symmetric and positive definite\cite{Barrett1994}.

\subsection{Solvers tailored to homogenization problems}
\label{ssec:solvers-homogenization}

To overcome some of the challenges of solvers for parametric PDEs that arise when solving homogenization problems in a many-query context, different solution schemes specifically tailored to homogenization problems have been developed in the last decades.
Some of the most important methods are outlined in the following.

\subsubsection{Lippmann-Schwinger equation}

The foundation for many solvers tailored to homogenization problems is laid by the Lippmann-Schwinger decomposition that originates from computational quantum mechanics~\cite{Lippmann1950}.
There, it is convenient to utilize a decomposition of the Hamiltonian with an unknown eigenfunction into a suitable free Hamiltonian with a known eigenfunction that can be solved using a corresponding Green's function, and into a remaining perturbation potential.
The Lippmann-Schwinger ansatz is also the basis for many solvers tailored to homogenization problems.
While these solvers are usually directly formulated for homogenization problems as introduced in \cref{sec:problems}, we will formulate them in the general framework of \ParPDEs{} from \cref{ssec:parPDE}.
To this end, we decompose the differential operator of the linear \ParPDE{} \eqref{eq:parametric-PDE} into a spatially homogeneous differential operator~$\diffOp_{\refmat}$ which is induced by a chosen set of reference parameters and a remaining differential operator for the perturbation~$\pert{\diffOp}_{\param}$ as
\begin{align} \label{eq:Lippmann-decomposition}
	\diffOpParam = \diffOp_{\refmat} + \pert{\diffOp}_{\param} \,, &&
	\pert{\diffOp}_{\param} = \diffOpParam - \diffOp_{\refmat} \,.
\end{align}
Plugging this ansatz into~\eqref{eq:parametric-PDE} yields
\begin{align} \label{eq:Lippmann-ansatz}
	\diffOp_{\refmat} \, \solParam
	= \rhsParam - \pert{\diffOp}_{\param} \, \solParam
	= \rhsParam + \left(\diffOp_{\refmat} - \diffOpParam\right) \solParam \,.
\end{align}
Since $\diffOp_{\refmat}$ is based on spatially homogeneous parameters, it is often possible to derive a corresponding Green's function $\Greens_{\refmat}$ such that~$\diffOp_{\refmat} \Greens_{\refmat} = \delta_{\fx}$.
Applying this to~\eqref{eq:parametric-PDE} leads to the fixed-point integral equation
\begin{align} \label{eq:Lippmann-convolution}
	\solParam = \int_{\domain} \Greens_{\refmat}(\fx - \fy) \, \pert{\rhs}_{\param}(\fx) \dInt{\fx} =
	\Greens_{\refmat} \conv \pert{\rhs}_{\param} \,, &&
	\pert{\rhs}_{\param} = \rhsParam - \pert{\diffOp}_{\param} \, \solParam \,,
\end{align}
that is also known as the Lippmann-Schwinger equation. The problem \cref{eq:Lippmann-convolution} can be used to iterate for $\solParam$ using the fixed-point iteration
\begin{align} \label{eq:Lippmann-iteration}
	\sol^{(\cgI)}_{\param} = \Greens_{\refmat} \conv \left( \rhsParam - \pert{\diffOp}_{\param} \solParam^{(\cgI-1)} \right) = \sol^{(\cgI-1)}_{\param} -  \Greens_{\refmat} \conv \left( \diffOpParam \, \sol^{(\cgI-1)}_{\param} - \rhsParam \right) \,,
\end{align}
where~$\cgI \in \ffN$ is the iteration number, and $\sol^{(0)}_{\param}$ an initial guess.

\subsubsection{Collocation method of Moulinec-Suquet}
\label{sssec:Moulinec-Suquet}

In the pioneering work~\cite{Moulinec1994}, one of the first solvers for homogenization problems that employ the Fast Fourier Transform (FFT) is presented.
Nowadays, this is referred to as the \textit{collocation method} of Moulinec-Suquet.
The main idea is to perform the convolution with an analytically derived Green's function~$\Greens_{\refmat}$ as in~\eqref{eq:Lippmann-convolution} {\color{revision}in Fourier space at $\cO(n \log n)$ cost using the FFT}.
While solvers for homogenization problems based on the FEM typically act on the primary variables such as displacements or temperature, the collocation method of Moulinec-Suquet acts on secondary variables that depend on these, e.g., on strains and temperature gradients.
Reformulating the parametric PDE in these variables simplifies the solution process and gives rise to a Green's function~$\Greens_{\refmat}$ for a spatially homogeneous reference material with parameters~$\param_{\refmat}$ that can be easily represented in Fourier space.
Reformulating the convolution with the Green's function $\Greens_{\refmat}$ in~\eqref{eq:Lippmann-convolution} using the convolution theorem leads to
\begin{align}
    \sol^{(\cgI)}_{\param} = \sol^{(\cgI-1)}_{\param} - \trafoFFTInv \left( \widehat{\Greens}_{\refmat} \odot \trafoFFT \left( \diffOp_{\param} \sol^{(\cgI-1)}_{\param} - \rhsParam \right) \right) \,,
	&& \text{with} &&
	\widehat{\Greens}_{\refmat} = \trafoFFT \left( \Greens_{\refmat} \right) \,.
\end{align}

Numerous extensions of the Moulinec-Suquet scheme (also: \textit{basic scheme}) are available nowadays.
For instance, a generalization to nonlinear composites is presented in \cite{Moulinec1998}.
In \cite{Zeman2010}, a significant acceleration of the basic scheme is achieved by replacing the fixed-point iteration by an iterative scheme that is similar to the preconditioned \emph{CG} method in \Cref{alg:pcg}.
A comparison of different variants based on the basic scheme is available in \cite{Moulinec2014}, and a variety of FFT-based solvers are summarized in \cite{Lucarini2021}.

However, the collocation method of Moulinec-Suquet and similar solution schemes that are purely FFT-based instead of employing a FEM discretization may suffer from various disadvantages.
For instance, checker-boarding or aliasing effects can be observed that originate from the formulation of the Green's function in Fourier space.
Also, Gibb's phenomenon introduces artifacts due to the mathematical nature of the Fourier series itself \cite{Leuschner2018,Carslaw1925}, mainly triggered by jumps in the material parameters.
Moreover, the extensive mathematical theory behind the FEM does not apply.

\subsubsection{Space Lippmann Schwinger (SLS) scheme}
\label{sssec:sls}

The Space Lippmann Schwinger scheme (SLS \cite{Yvonnet2012}) is another specialized solver for homogenization problems based on the Lippmann-Schwinger equation in~\eqref{eq:Lippmann-iteration} just like the collocation method of Moulinec-Suquet.
However, it avoids the artifacts due to Gibb's phenomenon by performing the Lippmann-Schwinger update rule \eqref{eq:Lippmann-iteration} in each iteration based on a FEM discretization instead of directly in Fourier space.
Here, the convolution with $\Greens_{\refmat}$ is computed in the real space with the help of truncated transformation tensors that are obtained using a FEM discretization \cite{Yvonnet2012}.
One main disadvantage of the Space Lippmann Schwinger scheme is that it is {\color{revision}computationally} less efficient compared to the collocation method of Moulinec-Suquet, as it does not rely on FFT-acceleration for the convolution operation {\color{revision}and thus exhibits a higher per-iteration cost}.

\subsubsection{Fourier-Accelerated Nodal Solvers (\FANS)}
\label{sssec:fans}

\begin{figure}[ht]
	\centering
    \includegraphics[scale=1.0]{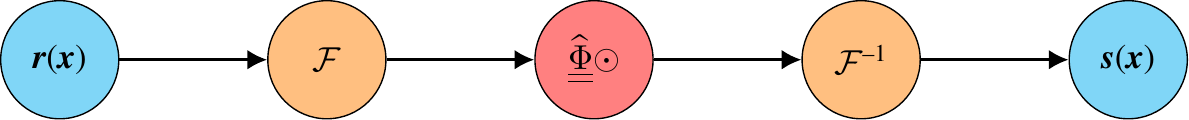}
	\caption{Structure of the \FANS{} preconditioner based on the Fourier transform~$\trafoFFT$ and a precomputed fundamental solution in Fourier space~$\FundSolMatFourier$. Note the analogy to \cref{fig:fno-architecture}.}
	\label{fig:FANS-architecture}
\end{figure}

In recent years, attempts have been made to merge {\color{revision}highly computationally efficient} FFT-based solvers like the collocation method of Moulinec-Suquet from \cref{sssec:Moulinec-Suquet} into the established framework of the FEM, e.g., via the development of Fourier-Accelerated Nodal Solvers (\FANS \cite{Leuschner2018}).
As the name suggests, the primary variables for \FANS{} are the nodal values of a FEM discretization instead of secondary variables as in the collocation method of Moulinec-Suquet.
An alternative approach to \FANS{} that also connects the FEM with Fourier-based solvers was simultaneously developed in \cite{Schneider2017}.
The goal of both approaches is to obtain solvers that reproduce the same results as a FEM discretization, i.e., avoiding the artifacts that purely FFT-based solvers usually introduce, but with an efficiency that is on par with classical FFT-based solvers.
\FANS{} also use a Lippman-Schwinger ansatz as in~\eqref{eq:Lippmann-ansatz}.
However, this is rewritten directly in the discretized setting as
\begin{align} \label{eq:FANS-Lippmann}
	\StiffMat_{\refmat} \, \solDParam = \rhsDParam + \left( \StiffMat_{\refmat} - \StiffMat \right) \solDParam \,,
\end{align}
where~$\StiffMat_{\refmat}$ is the FEM stiffness matrix arising from the discretization of the differential operator~$\diffOp_{\refmat}$ corresponding to a reference material with spatially homogeneous parameters~$\param_{\refmat}$, and~$\StiffMat$ is the FEM stiffness matrix for the parametric differential operator~$\diffOp_{\param}$.
Under the assumption that~$\StiffMat_{\refmat}$ is invertible, the ansatz in~\eqref{eq:FANS-Lippmann} is equivalent to the preconditioned linear system
\begin{align} \label{eq:FANS-preconditioned}
	\StiffMat_{\refmat}^{-1} \StiffMat\, \solDParam = \StiffMat_{\refmat}^{-1} \rhsDParam \,.
\end{align}
Comparing \cref{eq:FANS-preconditioned} and the preconditioned linear system \cref{eq:prec-linear-system}, the role of the preconditioner is taken by $\StiffMat_{\refmat}^{-1}$, which is symmetric and positive definite for suitable reference parameters~$\param_{\refmat}$.
Hence, \eqref{eq:FANS-preconditioned} suggests the use of iterative solvers such as the \CG{} method, as long as $\cond ( \StiffMat_{\refmat}^{-1} \, \StiffMat ) \ll \cond( \StiffMat)$.
Fortunately, this is the case if $\param_{\refmat}$ is chosen properly with respect to the set of admissible parameters $\paramSpace$ for the parametric PDEs.
The idea of preconditioners based on a reference material such as \FANS{}  is discussed in detail in \cite{Ladecky2023}, where it is also extended to other discretizations.
For spatially homogeneous reference parameters~$\param_{\refmat}$ and if, e.g., the FEM on a regular grid is used as discretization approach, the matrix~$\StiffMat_{\refmat}$ has a block circulant structure as discussed in \cite{Ladecky2023}.
Due to the special structure of $\StiffMat_{\refmat}^{-1}$, the preconditioner can be applied efficiently using a matrix-valued fundamental solution~$\FundSolMat: \domain \to \ffR^{\nComp \times \nComp}$ that can be interpreted similarly to a Green's function in the discretized setting.
For that, we introduce an underlying preconditioner operator~$\Prec$, that maps a function~$\res$ associated with the residual vector in the discrete setting~$\resD = \vecOp{\field{\res}}$ to a function~$\precRes$ associated with the preconditioned residual vector~$\precResD = \vecOp{\field{\precRes}}$.
In \FANS{}, the application of the preconditioning operator~$\precRes = \Prec \res$ is given by the convolution of the residual field~$\res$ with this fundamental solution
\begin{align} \label{eq:FANS-prec-convolution}
    \precResSymbol_i(\fx)
    = \sum_{j=1}^\nComp \FundSol_{ij} \conv \resSymbol_j
    = \sum_{j=1}^\nComp \int_\domain \FundSol_{ij}(\fx - \fy) \resSymbol_j(\fy) \dInt{\fy} \,, &&
    1 \leq i \leq \nComp \,.
\end{align}
Again, the convolution is evaluated in the Fourier space by virtue of the convolution theorem as
\begin{align} \label{eq:FANS-prec-operator}
    \precResSymbol_i(\fx) = \trafoFFTInv \left( \sum_{j=1}^\nComp \FundSolFourier_{ij} \odot \trafoFFT \left( \resSymbol_j(\fx) \right) \right)
    \,, &&
    \FundSolFourier_{ij} = \trafoFFT \left( \FundSol_{ij} \right) \,, &&
    1 \leq i, j \leq \nComp \,,
\end{align}
In the discrete setting, the application of the \FANS{} preconditioner, i.e., $\precResD = \StiffMat_{\refmat}^{-1} \resD$, can be computed via
\begin{align} \label{eq:FANS-prec-fields}
    \precResD = \vecOp{ \field{\precRes} } \,, &&
    \resD = \vecOp{ \field{\res} } \,, &&
    \field{ \precResSymbol_i } = \trafoFFTInv \left( \sum_{j=1}^\nComp \field{ \FundSolFourier_{ij} } \odot \trafoFFT \left( \field{ \resSymbol_j } \right) \right) \,, && 1 \leq i \leq \nComp \,.
\end{align}
at a time complexity of~$\orderOf{n \log n}$ when relying on established FFT algorithms.
In this formulation, the \FANS{} preconditioner maps the residual vector~$\resD$ that contains the values of the discrete field~$\field{ \res }$ to the preconditioned residual vector~$\precResD$ containing the values of~$\field{ \precRes }$, see also \cref{fig:FANS-architecture}.

\begin{figure}[ht]
    \includegraphics[scale=1.0]{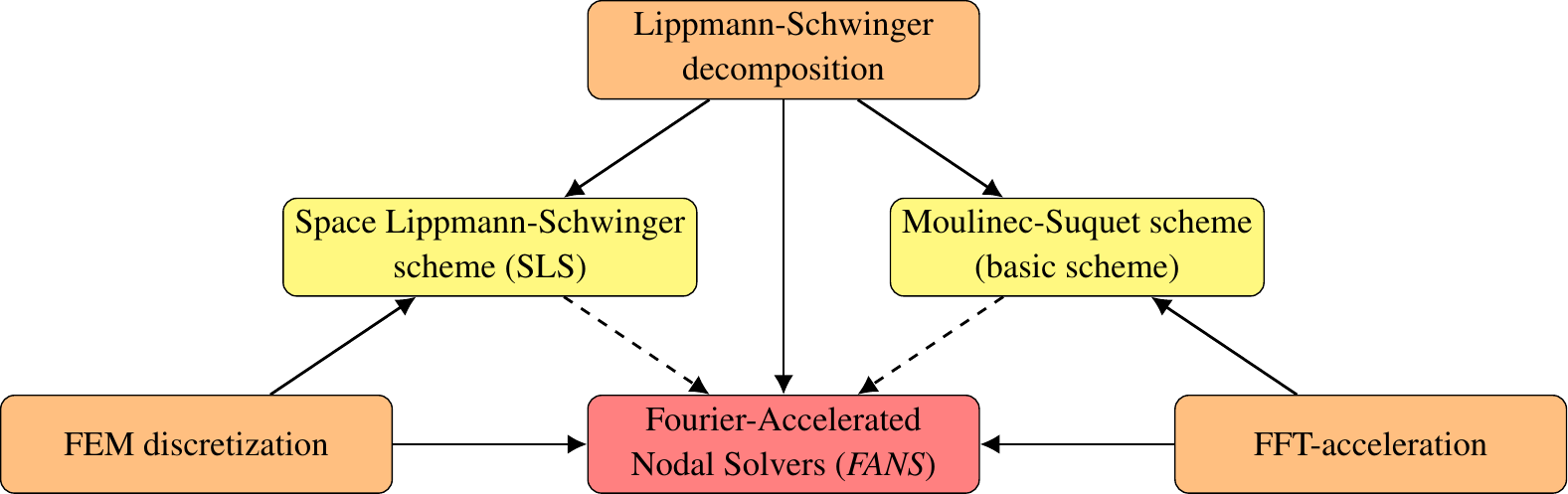}
	\caption{\normalsize Relations of \FANS{} (red) to other solvers (yellow) for homogenization problems.}
	\label{fig:FANS-relations}
\end{figure}

The matrix-valued fundamental solution~$\FundSolMatFourier$ in Fourier space can be assembled on a given discrete grid either by solving an auxiliary linear system using the matrix~$\StiffMat_{\refmat}$ {\color{revision}or using an equivalent matrix-free operator. Alternatively, the fundamental solution can also be constructed directly in Fourier space based on the FEM gradient operators as explained in~\cite{Leuschner2018}}.
Hence, assembling the fundamental solution for \FANS{} requires some sort of expert knowledge in the form of the suitable reference parameters~$\param_{\refmat}$ and the FEM gradient operators used in the discretization.
The relations of \FANS{} with other solvers tailored to homogenization problems are shown in \cref{fig:FANS-relations}.
An extension of \FANS{} is available in \cite{Keshav2022}.
Further solvers designed for homogenization problems can be found in, e.g., \cite{Lucarini2021,zeman_finite_2017,Schneider2019,vondrejc_fft-based_2014,Risthaus2024}.

\section{Machine learning enhanced solvers for parametric PDEs}
\label{sec:mlsolvers}

\subsection{Machine-learned preconditioners}
\label{ssec:mlprec}

On the one hand, machine learning could be used to provide an initial guess for an iterative solver, such that it converges faster.
By that, a machine-learned surrogate model is used as a prior for the solver.
This is typically easy to realize as there are few requirements on the initial guess, but the reduction in computing time is often only moderate, as the solver still accounts for most of the computational effort, unless the solution is directly forecasted almost perfectly.
On the other hand, ML models can be integrated on a deeper level into iterative solvers to leverage their performance.
A common approach in numerical analysis that is backed by an extensive theoretical framework is the use of preconditioners, which are invoked in each iteration of the iterative solver.
It has been shown many times that by using algebraic-based preconditioners, e.g., based on an incomplete Cholesky factorization (IC\cite{Saad2003}), or physical-based preconditioners such as \FANS{}~\cite{Leuschner2018}, the number of iterations required can be considerably reduced, almost independent of the quality of the initial guess.
This raises the question of whether machine-learned models can also take over the role as preconditioners, 
and to what extent convergence can be guaranteed.
The latter point is of crucial relevance, as the lack of error control is one major disadvantage over conventional machine-learned surrogates that forecast the solution directly.

\begin{figure}[ht]
	\centering
    \includegraphics[scale=1.0]{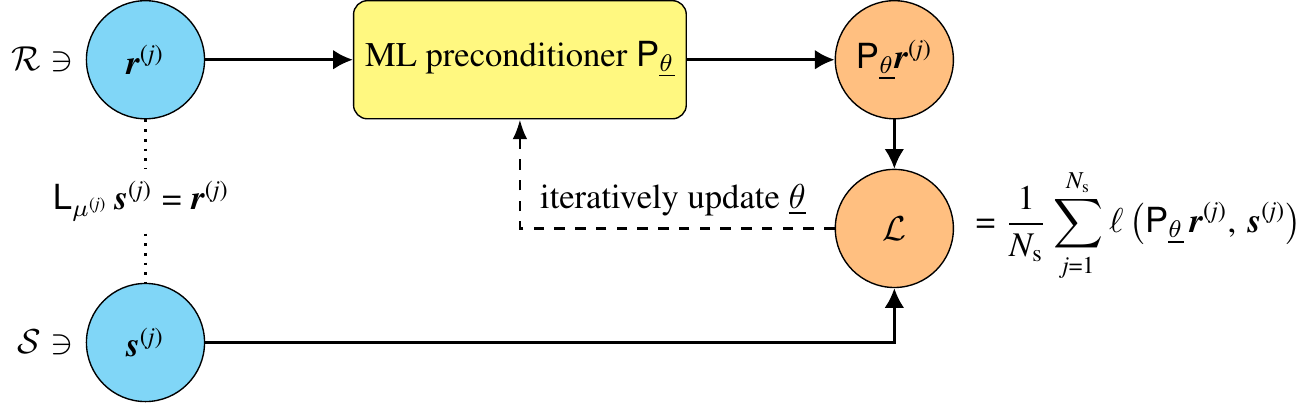}
        \vspace{-0.5\baselineskip}
	\caption{Training a machine-learned preconditioner~$\Prec_\weights$ with learnable weights~$\weights$ based on available training data sets~$\PrecInSet$ and~$\PrecOutSet$ for the parametric PDE with differential operator~$\diffOpParam$ at hand.}
	\label{fig:mlprec-training}
\end{figure}

Let~$\Prec_\weights: \SquareIntSpace(\Omega,\, \ffR^\nComp) \to \SquareIntSpace(\Omega,\, \ffR^\nComp)$ be the operator of a machine-learned surrogate such as a neural operator that is based on some vector of learnable weights~$\weights \in \ffR^\nWeights$, where~$\nWeights \in \ffN$, and that shall be used as a preconditioner.
We consider training data to calibrate the model in the form of finite sets of functions
\begin{align} \label{eq:prec-train-data}
    \PrecInSet = \left\{ \PrecInSample: 1 \leq j \leq \nSamples \right\} \subsetneq \SquareIntSpace(\Omega,\, \ffR^\nComp) \,, &&
    \PrecOutSet = \left\{ \PrecOutSample: 1 \leq j \leq \nSamples \right\} \subsetneq \cU \subseteq \SquareIntSpace(\Omega,\, \ffR^\nComp) \,, && \nSamples \in \ffN \,,
\end{align}
that satisfy the \ParPDE{} for some given arbitrary parameters~$\param\sampleIdx$, that is,
\begin{align}
	\diffOpParamSample \PrecOutSample(\fx) = \PrecInSample(\fx) \,, &&
	\param\sampleIdx \in \paramSpace \,, &&
	1 \leq \sampleI \leq \nSamples \,.
\end{align}
Hence, the data generation consists of solving the parametric PDE for many different parameters~$\param\sampleIdx$ and right-hand sides~$\PrecInSample$.
Given training data~$\PrecInSet$ and~$\PrecOutSet$, the preconditioner~$\Prec_\weights$ can be trained by minimizing a loss function of the form
\begin{align} \label{eq:prec-loss}
    \!\min_\weights \quad \Loss \,, && \text{where} && 
    \Loss = \frac{1}{\nSamples}
    \sum_{\sampleI = 1}^\nSamples \lossSample \,, &&
    \lossSample = \ell \left( \Prec_\weights \, \PrecInSample ,\, \PrecOutSample \right) \,,
\end{align}
with a suitable sample-wise loss function~$\ell: \SquareIntSpace\left( \domain;\, \ffR^\nComp \right) \times \SquareIntSpace\left( \domain;\, \ffR^\nComp \right) \to \ffR$ that measures the deviation of its arguments (prediction and reference fields).
By solving the optimization problem~\eqref{eq:prec-loss}, we aim for a machine-learned preconditioner that is a good approximation of the inverse of the parametric differential operator~$\diffOpParam$ over the set of parameters~$\paramSpace$, that is,
\begin{align} \label{eq:mlprec-objective}
	\Prec_\weights \, \diffOpParamSample \approx \identity \,, &&
	\forall \param\sampleIdx \in \paramSpace \,.
\end{align}
Note that in this framework, the machine-learned preconditioner~$\Prec_{\weights}$ does not depend on the specific parameters~$\param\sampleIdx \in \paramSpace$ of the \ParPDE{}.
This is also the case for many classical preconditioners, including \FANS{} (\cref{sssec:fans}).
Since the machine-learned preconditioner~$\Prec_\weights$ is not aware of the present parameters and only has a limited design space, the objective in \eqref{eq:mlprec-objective} obviously cannot be fulfilled exactly.
The ansatz of training a machine-learned preconditioner by solving the optimization problem in~\eqref{eq:prec-loss} based on the training data defined in~\eqref{eq:prec-train-data} is visualized in \cref{fig:mlprec-training}.

A remaining crucial aspect is how to guarantee the convergence for hybrid solvers that feature a machine-learned preconditioner.
For efficient Krylov subspace methods such as the \CG{} solver, there are strict requirements on the preconditioner that are not easy to fulfill for machine-learned models.
A machine-learned surrogate, for example, could fulfill the objective in~\eqref{eq:mlprec-objective} very well, but as long as the convergence of the resulting hybrid solver can not be proven, all guarantees of the iterative numerical scheme are lost, representing a severe drawback.

\subsection{Unitary Neural Operators as a preconditioner for \texorpdfstring{\CG{}}{CG} (\texorpdfstring{\UNOCG{}}{UNO-CG})}
\label{ssec:fnocg}

While most approaches for machine-learned preconditioners in the literature, as summarized in \cref{ssec:introduction-hybrid-solvers}, only deal with numerical iterative solvers having no special requirements on the preconditioner itself, such as (flexible)~$\GMRES$, in the following, a method is presented to accelerate a conjugate gradient (\CG{}) solver with a special machine-learned preconditioner.
To guarantee the convergence of a preconditioned \CG{} solver, the preconditioner~$\Prec_\weights$ must satisfy:

\begin{enumerate}[label=\textbf{[A\arabic*{}]}]
\item \label{prop:lin} \textbf{Linearity}:   $\Prec_\weights \left( \alpha \fr^{(i)} + \fr^{(j)} \right) = \alpha \, \Prec_\weights \, \fr^{(i)} + \Prec_\weights \, \fr^{(j)}, \; \forall \alpha \in \ffR, \; \forall \fr^{(i)}, \fr^{(j)}$
\item \label{prop:sym} \textbf{Symmetry}: $\fr^{(i)} \cdot \Prec_\weights \, \fr^{(j)} = \fr^{(j)} \cdot \Prec_\weights \, \fr^{(i)} \; \forall \fr^{(i)}, \fr^{(j)}$
\item \label{prop:pos} \textbf{Positive definiteness}: $\fr \cdot \Prec_\weights \, \fr > 0 \; \forall \fr \neq \zero \,.$
\end{enumerate}

Furthermore, the condition of the problem should be improved (considerably) to lower the number of needed solver iterations. In the following, we carefully design a machine learning model meeting \ref{prop:lin}--\ref{prop:pos}, and, thus, with convergence guarantees.

\subsubsection{Architecture}
\label{sssec:fnocg-architecture}
\begin{figure}[ht]
	\centering
    \includegraphics[scale=1.0]{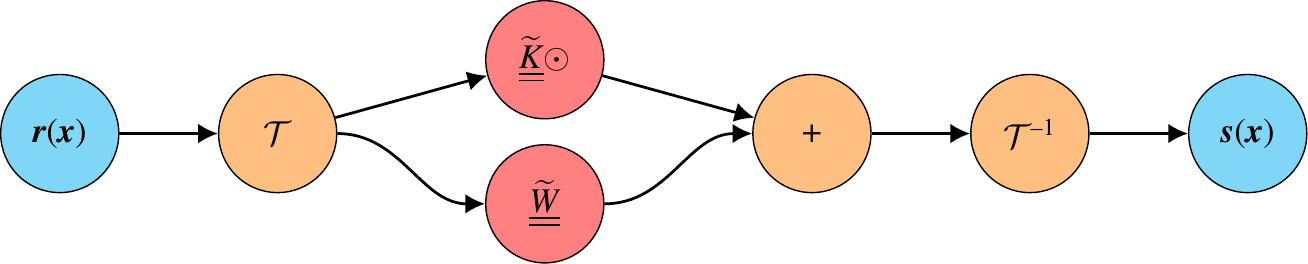}
	\caption{Architecture of the \UNO{} preconditioner used in \UNOCG{} based on a unitary transform~$\trafo$ that features a learnable kernel~$\fnocgConvKernelFourier$ and a learnable inner bypass~$\fnocgBypass$. Note the similarity to \cref{fig:fno-architecture} and \cref{fig:FANS-architecture}.}
	\label{fig:fnocg-architecture}
\end{figure}

We propose an architecture inspired by a single Fourier layer from FNOs as introduced in \cref{sssec:fno}. Through specific modifications, the properties \ref{prop:lin}--\ref{prop:pos} are enforced.
First, \emph{linearity} \ref{prop:lin} is achieved by restricting the nonlinear activation to the identity mapping~$\fnoAct{} = \mathrm{id}$.
Additionally, the forward and inverse Fourier transform~$\trafoFFT, \trafoFFTInv$ are replaced by a more general unitary transform~$\trafo$ and its inverse~$\trafoInv{}=\trafo^\ast$, i.e., the inverse transform is equal to the adjoint transform.
This relaxation represents a significant generalization compared to classical Fourier layers: Unitary transforms comprise, e.g., the Discrete Fourier Transform (DFT), but also the discrete sine and cosine transform (DST, DCT). Most importantly, implementations of these algorithms with complexity~$\orderOf{n \,\log n}$ exist, yielding comparable performance for all of them\cite{FFTW05}.
Second, the number of output channels~$\nChannelsLayer$ is fixed to the number of nodal DOF~$\nComp$, which is one of several requirements to ensure symmetry~\ref{prop:sym}.
Furthermore, the bypass is now applied in the transformed space, i.e., an \emph{inner bypass} is used.
This results in the operator
\begin{align}\label{eq:fnocg-operator}
	\precRes(\fx) &= \trafoInv \left( \fnocgBypass \, \trafo \left( \res(\fx) \right) + \fnocgConvKernelFourier \odot \trafo \left( \res(\fx) \right) \right) \,,
\end{align}
with a bypass matrix~$\fnocgBypass \in \ffR^{\nComp \times \nComp}$ and a matrix-valued function~$\fnocgConvKernelFourier: \ffR^{\nDims} \to \ffR^{\nComp \times \nComp}$.
Both~$\fnocgBypass$ and~$\fnocgConvKernelFourier$ must be learned from training data.
We refer to this modified Fourier layer as \UNO{} preconditioner, whose architecture is shown in \cref{fig:fnocg-architecture}.

\subsubsection{Parametrization}
\label{sssec:fnocg-parametrization}

Exactly like in FNOs as introduced in \cref{sssec:fno}, our proposed \UNO{} preconditioner only acts on a finite number of (usually low-frequency) modes from a set~$\modeSet \in \ffZ^{\nDims}$.
This is realized by parametrizing~$\fnocgConvKernelFourier$ only on these discrete modes~$\mode\modeIdx \in \modeSet$ via learnable weights.
While for FNOs,~$\fnoBypass$ is a general real-valued matrix and $\fnoConvKernelFourier$ is allowed to be complex-valued,  in the \UNO{} preconditioner we restrict~$\fnocgBypass$ and~$\fnocgConvKernelFourier$ to be symmetric and real-valued to guarantee the symmetry~\ref{prop:sym}.
In order to ensure \ref{prop:pos}, both~$\fnocgBypass$ and~$\fnocgConvKernelFourier$ are required to be positive definite.
For that, we employ a local parametrization for symmetric and positive definite~$\nComp \times \nComp$ matrices $\localParam: \ffR^{\nPrecC} \to \spdSpace{\ffR^{\nComp \times \nComp}}$
that constructs these matrices from vectors containing~$\nPrecC = \frac{\nComp(\nComp+1)}{2}$ values.
We propose an approach to assign weights~$\weights^{\langle \modeI \rangle} \in \ffR^{\nPrecC}$ to the symmetric bypass matrix~$\fnocgBypass \in \ffR^{\nComp \times \nComp}$ and to the symmetric matrix-valued kernel~$\fnocgConvKernelFourier\left(\mode^{\langle \modeI \rangle}\right) \in \ffR^{\nComp \times \nComp}$ for each mode~$\mode^{\langle \modeI \rangle} \in \modeSet$ with~$1 \leq \modeI \leq \nModes$ via
\begin{align} \label{eq:fnocg-parametrization}
    \fnocgBypass = \localParam \left( \weights\modeIdxBypass \right) \,, &&
    \fnocgConvKernelFourier(\mode) = \begin{cases}
        \localParam \left( \weights\modeIdx \right) & \text{if} \quad \exists \modeI: \mode = \mode\modeIdx \in \modeSet \,, \\
        \localParam \left( \weights\modeIdx \right) & \text{if} \quad \exists \modeI: \mode = \conjOp{\mode\modeIdx} \in \modeSetConj \,, \\
        \ull{0} & \text{else} \,,
    \end{cases} &&
    \mode \in \ffZ^{\nDims} \,,
\end{align}
where a set of conjugate modes~$\modeSetConj$ similar to FNOs (see \cref{sssec:fno}) is considered.
The set~$\modeSetConj$ shall be chosen such that the output $\precRes$ of the \UNO{} preconditioner is always real-valued.
For the Fourier transform~$\trafo=\trafoFFT$ this is given by~\eqref{eq:modesetconj-fno}. For the sine transform~$\trafo=\trafoDST$, an empty set~$\modeSetConj=\emptySet$ is sufficient, for example.
Further, the weights~$\weights\modeIdx \in \ffR^{\nPrecC}$ are collected in the vector
\begin{align}
	\weights = \T{ 
     \begin{bmatrix}
 	 \T{\left(\weights\modeIdxBypass\right)} & \T{\left(\weights\modeIdxFirst\right)} & \cdots & \T{\left(\weights\modeIdxLast\right)}
 	 \end{bmatrix} }
    \in \ffR^{\nWeights} \,, &&
	\nWeights = \nPrecC \left( \nModes + 1 \right) \,.
\end{align}

For the local parametrization~$\localParam$, we utilize the Cholesky-like decomposition of a symmetric~$\nComp \times \nComp$ matrix based on a triangular matrix having only~$\nPrecC = \frac{\nComp(\nComp+1)}{2}$ nonzero entries.
This approach is explained below for the most relevant special cases.

\paragraph*{Special case of one nodal DOF ($c=1$):}
In this case, the local parametrization~$\localParam: \ffR \to \ffR$ can simply be defined as $\localParam \left( \weight\modeIdx \right) = \left( \weight\modeIdx \right)^2$ for $0 \leq \modeI \leq \nModes$.
\paragraph*{Special case of two nodal DOF ($c=2$):}
Here, the local parametrization~$\localParam: \ffR^{3} \to \spdSpace{\ffR^{2 \times 2}}$ is given by
\begin{align} \label{eq:local-param-2}
	\localParam \left( \weights\modeIdx \right) = \localParam \left( \begin{bmatrix} \weight\modeIdx_1 \\ \weight\modeIdx_2 \\ \weight\modeIdx_3 \end{bmatrix} \right) = \begin{bmatrix} \weight\modeIdx_1 & 0 \\ \weight\modeIdx_3 & \weight\modeIdx_2 \end{bmatrix} \T{\begin{bmatrix} \weight\modeIdx_1 & 0 \\ \weight\modeIdx_3 & \weight\modeIdx_2 \end{bmatrix}} = \begin{bmatrix} \weight\modeIdx_1 \cdot \weight\modeIdx_1 & \weight\modeIdx_1 \cdot \weight\modeIdx_3 \\ \weight\modeIdx_1 \cdot \weight\modeIdx_3 & \weight\modeIdx_2 \cdot \weight\modeIdx_2 + \weight\modeIdx_3 \cdot \weight\modeIdx_3 \end{bmatrix} \,, &&
	0 \leq \modeI \leq \nModes \,,
\end{align}
which is based on a Cholesky-like (note that the diagonal of the lower/upper triangular matrix is not enforced to be positive here, but it is in the Cholesky factorization) decomposition for each mode.
\paragraph*{Special case of three nodal DOF ($c=3$):}
Analogously, the local parametrization in the case of three nodal DOF is defined in \cref{sec:local-parametrization-c-3}.

\vskip1em

\begin{figure}[ht!]
    \centering
    \begin{subfigure}[b]{0.32\textwidth}
        \centering
        \includegraphics[scale=0.69]{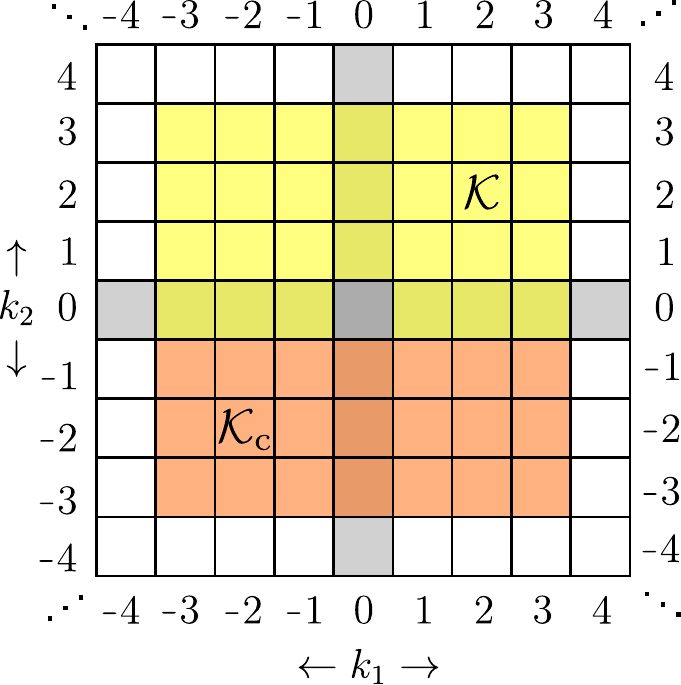}
        \vspace{-0.2\baselineskip}
        \caption{Periodic BC.}
        \label{fig:fnocg-mode-selection-per2d}
    \end{subfigure}
    \hfill
    \begin{subfigure}[b]{0.32\textwidth}
        \centering
        \includegraphics[scale=0.69]{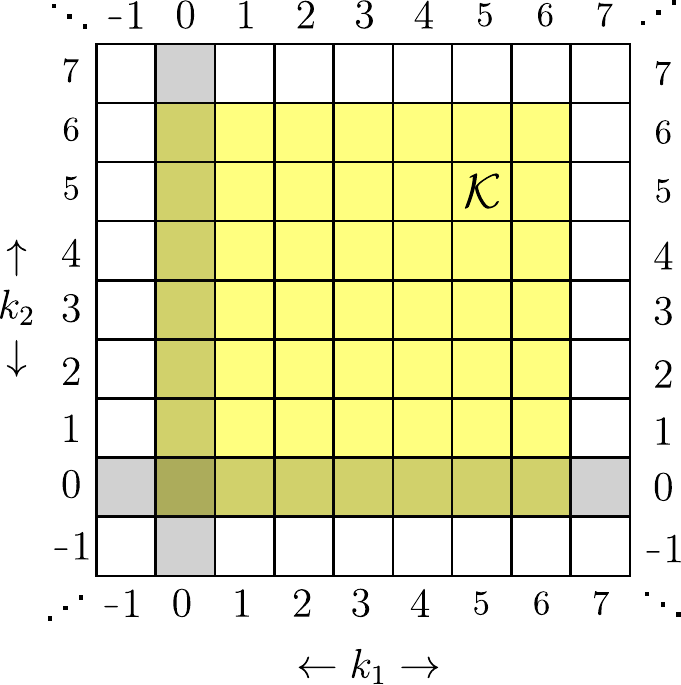}
        \vspace{-0.2\baselineskip}
        \caption{Dirichlet BC.}
        \label{fig:fnocg-mode-selection-dir2d}
    \end{subfigure}
    \hfill
    \begin{subfigure}[b]{0.32\textwidth}
        \centering
        \includegraphics[scale=0.69]{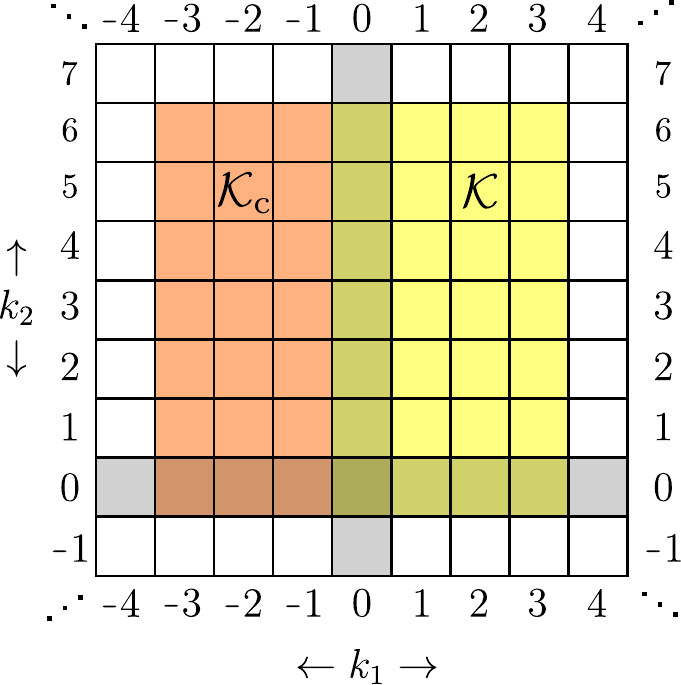}
        \vspace{-0.2\baselineskip}
        \caption{Mixed BC (per. x, Dir. y).}
        \label{fig:fnocg-mode-selection-mixed2d}
    \end{subfigure}
    \caption{Visualization of the mode selection strategy from \cref{tab:fnocg-mode-selection-2d} based on~$\modeSet$ (yellow) and~$\modeSetConj$ (orange) for an exemplary hyperparameter of~$\modeParam = 3$ and 2D problems ($\nDims=2$) having different boundary conditions.
    }
    \label{fig:fnocg-mode-selection}
\end{figure}

These variants of local parametrizations~$\localParam$ and their derivatives with respect to~$\weights\modeIdx$ for~$0 \leq \modeI \leq n$ can be evaluated easily in~$\orderOf{n}$, which makes them convenient to use in the preconditioner training.
The selection of the mode sets~$\modeSet, \modeSetConj \subsetneq \ffZ^{\nDims}$ is performed depending on the problem dimension~$\nDims$, the deployed unitary transform~$\trafo$, and the boundary conditions of the parametric PDE to consider.
As explained in \cref{sssec:fno} for FNOs, a hyperparameter~$\modeParam \in \ffN$ is also used here to determine the modes in the sets~$\modeSet$ and~$\modeSetConj$ in a straightforward way.
Examples of suitable mode selection strategies for different boundary conditions are given in \cref{tab:fnocg-mode-selection-2d} and visualized in \cref{fig:fnocg-mode-selection}.

\begin{table}[ht!]
	\caption{Mode selection strategies in \UNOCG{} based on a hyperparameter~$\modeParam \in \ffN$, the Fourier transform~$\trafoFFT$, and the sine transform~$\trafoDST$ that are suitable for different BCs and dimensions~$\nDims$. For simplicity, mixed BC are only considered for~$\nDims=2$.}
	\centering
	\begin{tabular}{ccccc}
		\toprule 
		Boundary condition & Learned modes $\modeSet$ & Conjugate modes $\modeSetConj$ & $\nModes=\norm{\modeSet}$
        & $\trafo$ \\ 
		\midrule
		Periodic & $\discrIntv{{-M}}{M}^{\nDims - 1} \times \discrIntv{0}{M} \setminus \{ \ul{0} \}$ & $\discrIntv{{-M}}{M}^{\nDims - 1} \times \discrIntv{{-M}}{{-1}}$ & $(2 \modeParam + 1)^{\nDims - 1}(\modeParam + 1) - 1$
        & $\trafoFFT$ \\
		Dirichlet & $\discrIntv{0}{2M}^{\nDims}$ & $\emptySet$ & $(2 \modeParam + 1)^{\nDims}$
        & $\trafoDST$ \\
		Mixed (e.g., per. $\mathrm{x}$, Dir. $\mathrm{y}$) & $\discrIntv{0}{M} \times \discrIntv{0}{2M}$ & $\discrIntv{{-M}}{{-1}} \times \discrIntv{0}{2M}$ & $(2 \modeParam + 1)(\modeParam + 1)$
        & $\trafoFFT^{\rm x} \circ \trafoDST^{\rm y}$ \\
		\bottomrule
	\end{tabular}
	\label{tab:fnocg-mode-selection-2d}
\end{table}

Similar to the fundamental solution~$\FundSolMatFourier$ of \FANS{} used in~\eqref{eq:FANS-prec-operator}, a matrix-valued function~$\fnocgFundSolMat: \ffZ^{\nDims} \to \ffR^{\nComp \times \nComp}$ can be defined for the \UNO{} preconditioner.
With that, the application of the \UNO{} preconditioner can be formulated similar to \FANS{}---see~\eqref{eq:FANS-prec-operator}---as
\begin{align} \label{eq:fnocg-op-fundsol}
    \precResSymbol_i(\fx) = \trafoInv \left( \sum_{j=1}^\nComp \fnocgFundSolSymbol_{ij} \odot \trafo \left( \resSymbol_j(\fx) \right) \right) \,, && 1 \leq i, j \leq \nComp \,, &&
    \fnocgFundSolMat(\mode) = \fnocgBypass + \fnocgConvKernelFourier(\mode) \,, &&
    \mode \in \ffZ^{\nDims} \,.
\end{align}
When the Fourier transform is used in the \UNO{} preconditioner, that is, if~$\trafo = \trafoFFT$, then \FANS{} could be \emph{exactly} reproduced by setting~$\fnocgFundSolMat = \FundSolMatFourier$ instead of the parametrization approach using a finite number of modes as in~\eqref{eq:fnocg-parametrization} and~\eqref{eq:fnocg-op-fundsol}.
However, for other unitary transforms~$\trafo$, there is no such interpretation.
Still, we can also learn a \UNO{} preconditioner, e.g., for the sine transform~$\trafo = \trafoDST$, and use it in a hybrid solver.
A comparison of the \UNO{} preconditioner with the Fourier layer used in FNOs (\cref{sssec:fno}) and the \FANS{} preconditioner (\cref{sssec:fans}) is available in \cref{tab:fnocg-comparison}.

\begin{table}[ht!]
	\caption{\normalsize Comparison of the structure of a Fourier layer used in FNOs, the \UNO{} preconditioner, and \FANS{}.}
	\centering \begin{tabular}{cccc}
		\toprule 
		 & Fourier layer in FNOs & \UNO{} preconditioner & \FANS{} \\ 
		\midrule 
		Operator & $\fnoAct \left( \fnoBypass \bullet + \trafoFFTInv \left( \fnoConvKernelFourier \odot \trafoFFT \left( \bullet \right) \right) \right)$ & $\trafoInv \left( \fnocgBypass \, \trafo \left( \bullet \right) + \fnocgConvKernelFourier \odot \trafo \left( \bullet \right) \right)$ & $\trafoFFTInv \left( \FundSolMatFourier \odot \trafoFFT \left( \bullet \right) \right)$ \\
		\midrule 
		\makecell{Learnable \\ weights} & $\fnoWeightTensor \in \ffC^{\nModes \times \nChannelsLayer \times \nChannelsLayerPrev}$, $\fnoBypass \in \ffR^{\nChannelsLayer \times \nChannelsLayerPrev}$ &
		$\weights \in \ffR^{\left(\nModes + 1\right) \nPrecC}$
		& --- \\
		\midrule
		\makecell{Fundamental \\ solution} & --- & $\fnocgFundSolMat: \ffZ^\nDims \to \ffR^{\nComp \times \nComp}$ & $\FundSolMatFourier: \ffZ^\nDims \to \ffR^{\nComp \times \nComp}$ \\
		\midrule 
		\makecell{Structure} & \cref{fig:fno-architecture} & \cref{fig:fnocg-architecture} & \cref{fig:FANS-architecture} \\
		\bottomrule
	\end{tabular}
	\label{tab:fnocg-comparison}
\end{table}

Using the \UNO{} preconditioner for the \CG{} method from \cref{alg:pcg} gives rise to a hybrid solver that we denote as \UNOCG{}.
It can be interpreted as a machine-learned generalization of the classical solver \FANS{} summarized in \cref{sssec:fans}, since instead of the analytically derived \FANS{} fundamental solution~$\FundSolMatFourier$ in Fourier space, a machine-learned fundamental solution~$\fnocgFundSolMat$ is used, see also \cref{tab:fnocg-comparison}.
While \FANS{} is only applicable to homogenization problems with periodic BC, \UNOCG{} is not limited to periodic BC but is, in theory, applicable to general PDEs with various boundary conditions that are discretized on a regular grid.
{\color{revision}Note that using this parametrization, the $\UNO{}$ preconditioner as a continuous operator is, in general, invariant of the resolution since only modes in a transformed space are learned. Thus, a \UNO{} preconditioner learned on data with a fixed resolution could also be applied to problems with different resolutions.}

\subsubsection{Analysis on an algebraic level}
\label{sssec:fnocg-analysis}

While the \UNO{} preconditioner is an operator that maps between function spaces just like general neural operators, in a hybrid solver, it is always evaluated on a discrete grid with a fixed resolution.
To facilitate further analysis, we introduce an algebraic notation to describe its action in the discrete setting.
For that, we assume that the discretization of the learned preconditioner~$\Prec_\weights$ is given by a matrix~$\PrecD_\weights$ that is determined by the learnable weights~$\weights$, and maps the residual vector~$\resD = \vecOp{ \field{ \res } } \in \ffR^{\nComp n}$ to the preconditioned residual vector~$\precResD = \vecOp{ \field{ \precRes } } \in \ffR^{\nComp n}$ via~$\precResD = \PrecD_\weights \, \resD$.
On a discrete grid with~$n$ nodes and~$\nComp$ nodal DOF, the unitary transform~$\trafo$ of the 
\UNO{} layer can be represented by a unitary matrix~$\trafoD \in \ffC^{\nComp n \times \nComp n}$ with~$\trafoDInv = \trafoDHer$ such that
\begin{align}
	\hat{\resD} = \trafoD \, \resD \,, &&
	\resD = \trafoDHer \, \hat{\resD} && \iff &&
	\resD = \vecOp{ \field{ \res } } \,, && \hat{\resD} = \vecOp{ \field{ \hat{\res} } } \,, &&
	\hat{\res} = \trafo \left( \res \right) \,.
\end{align}
Based on that, the mapping in \eqref{eq:fnocg-operator} can be expressed on an algebraic level.
In fact, the linear operator that maps~$\resD \in \ffR^{\nComp n}$ to~$\precResD \in \ffR^{\nComp n}$ can be formulated explicitly as a matrix-vector product with the \UNO{} preconditioner matrix~$\PrecD_\weights \in \ffR^{\nComp n \times \nComp n}$ via
\begin{align} \label{eq:FNOCG-prec-fields-disc}
	\precResD = \PrecD_\weights \, \resD \,, && \PrecD_\weights = \trafoDHer \PrecSparse \, \trafoD \,,
\end{align}
with~$\PrecSparse$ containing the values of the matrix-valued function~$\fnocgFundSolMat: \ffZ^{\nDims} \to \ffR^{\nComp \times \nComp}$ from \eqref{eq:fnocg-op-fundsol} on a discrete grid with~$n$ nodes.
The sparse matrix~$\PrecSparse$ is a band matrix that is defined as
\begin{align} \label{eq:prec-sparse-def}
	\PrecSparse =
	\begin{bmatrix}
		\diag{ \vecOp{ \field{ \fnocgFundSolSymbol_{11} } } } & \cdots & \diag{ \vecOp{ \field{ \fnocgFundSolSymbol_{1 \nComp} } } } \\
		\vdots & \ddots & \vdots \\
		\diag{ \vecOp{ \field{ \fnocgFundSolSymbol_{\nComp 1} } } } & \cdots & \diag{ \vecOp{ \field{ \fnocgFundSolSymbol_{\nComp \nComp} } } }
	\end{bmatrix} \in \ffR^{\nComp n \times \nComp n} \,.
\end{align}
Recall that the vectorization operator $\vecOp{\bullet}$ defined in \cref{ssec:notation} unwraps the values of $\field{\fnocgFundSolSymbol_{ij}}$ on a discrete grid into a vector according to a given mode ordering $\mode\modeIdx$ for $1 \leq \modeI \leq n$, such as
\begin{align} \label{eq:fnocg-local-matrix}
    \vecOp{ \field{ \fnocgFundSolSymbol_{lm} }} = \begin{bmatrix} \fnocgFundSolSymbol\modeIdxFirst_{lm} & \cdots & \fnocgFundSolSymbol^{\langle n \rangle}_{lm} \end{bmatrix} \in \ffR^{n} \,, &&
    \fnocgFundSolSymbol\modeIdx_{lm} = \fnocgFundSolSymbol_{lm} \left( \mode\modeIdx \right) \,, &&
    1 \leq l, m \leq c \,.
\end{align}
A possible interpretation of the \UNO{} preconditioner on an algebraic level is a Block-Jacobi-like preconditioner that consists of the sparse matrix $\PrecSparse$ that is, however, applied in a transformed space that is constructed by the given unitary transform $\trafo$ in the continuous setting---or the unitary matrix $\trafoD$ in the discrete setting.
Due to its symmetry, the unique nonzero entries of the band matrix~$\PrecSparse$ can be collected in a vector~$\PrecEntries \in \ffR^{\nPrecC n}$, which can also be separated into parts~$\PrecEntries_i \in \ffR^n$, where each of them contains the diagonal entries of a unique block of the symmetric band matrix~$\PrecSparse$ that is indexed with~$1 \leq i \leq \nPrecC$, as in
\begin{align}
    \PrecEntries = \T{\begin{bmatrix} \T{\PrecEntries}_1 & \cdots & \T{\PrecEntries}_\nPrecC \end{bmatrix}} \in \ffR^{\nPrecC n} \,, &&
    \PrecEntries_i \in \ffR^n \,, && 1 \leq i \leq \nPrecC = \frac{\nComp(\nComp + 1)}{2} \,.
\end{align}
To directly obtain the values~$\PrecEntries \in \ffR^{\nPrecC n}$ based on the learnable weights~$\weights \in \ffR^{\nWeights}$, we introduce a global parametrization as
\begin{align} \label{eq:global-param}
    \globalParam: \ffR^\nWeights \to \ffR^{\nPrecC n} \,, && \weights \mapsto \globalParam(\weights) = \PrecEntries \,,
\end{align}
that is implicitly defined based on the local parametrization~$\localParam$ together with~\eqref{eq:fnocg-parametrization},~\eqref{eq:fnocg-op-fundsol}, and~\eqref{eq:prec-sparse-def} and can be evaluted in~$\orderOf{n}$.

\subsubsection{Convergence analysis of \texorpdfstring{\UNOCG{}}{UNO-CG}} 
\label{sssec:fnocg-convergence}

While $\PrecSparse$ is a sparse matrix with only $\nComp^2 n$ nonzero entries, the preconditioner matrix~$\PrecD_\weights$ itself is a dense matrix. However, it does not need to be assembled since the preconditioner is always applied {\color{revision}by a sparse matrix-vector product with the sparse band matrix $\PrecSparse$ in Fourier space}.
Still, the algebraic form in~\eqref{eq:FNOCG-prec-fields-disc} is helpful for further analysis and especially to prove convergence.
For that, the spectrum of the preconditioner is of interest.

\begin{lemma}\label{spectrum-lemma}
The spectrum (set of eigenvalues) of the \UNO{} preconditioner matrix~$\PrecD_\weights$ is given by
\begin{align}
	\eig{\PrecD_\weights} = \eig{ \mathrm{blockdiag}\left(\fnocgFundSolMat\modeIdx\right) } = \bigcup\limits_{i=1}^{n} \eig{\fnocgFundSolMat\modeIdx} \,, && 1 \leq \modeI \leq n \,,
\end{align}
where $\eig{\bullet}$ denotes the set of eigenvalues of $\bullet$, and~$\fnocgFundSolMat\modeIdx$ is defined in \eqref{eq:fnocg-local-matrix} based on \eqref{eq:fnocg-op-fundsol} and \eqref{eq:fnocg-parametrization}.
\end{lemma}
See \cref{sec:appendix-proofs} for a proof.
Following \cref{spectrum-lemma}, the preconditioner matrix~$\PrecD_{\weights}$ is positive definite for suitable local parametrizations~$\localParam$ as introduced in \cref{sssec:fnocg-parametrization}.
Together with the symmetry of $\PrecD_{\weights} \in \ffR^{\nComp n \times \nComp n}$, the \CG{} solver from \cref{alg:pcg} with the \UNO{} preconditioner converges in exact arithmetic after at most $\nComp n$ iterations to the exact solution $\solDParam \in \ffR^{\nComp n}$ following the established convergence theory in the literature\cite{Hackbusch1994}.

Using \UNO{} as a preconditioner for the \CG{} method leads to our proposed hybrid solver \UNOCG{} with guaranteed convergence that is described in \cref{alg:unocg}.
The number of FLOPs (floating-point operations) that are required for applying the \UNO{} preconditioner in the \UNOCG{} hybrid solver in a discrete setting on a grid with $n$ nodes, is approximately given by
\begin{align} \label{eq:fno-flops}
	\underbrace{ \nComp n \log n }_{\text{apply } \trafo} \; + \; \underbrace{ \nComp^2 n }_{\text{apply $\PrecSparse$}} \; + \; \underbrace{ \nComp n \log n }_{\text{apply } \trafoInv} \; \in \; \orderOf{ n \log n } \,.
\end{align}

\begin{algorithm}
    \caption{\UNOCG{} hybrid solver (\cref{alg:pcg} with machine-learned \UNO{} preconditioner)}
\label{alg:unocg}
    \hspace*{\algorithmicindent} \textbf{Data:} $\StiffMat\in \ffR^{\nComp n \times \nComp n}$, $\rhsDParam \in \ffR^{\nComp n}$, $\solD\cgIdxFirst_{\param} \in \ffR^{\nComp n}$, $\PrecSparse \in \ffR^{\nComp n \times \nComp n}$, $\cgTol > 0$ \\
    \hspace*{\algorithmicindent} \textbf{Result:} $\solDParam \in \ffR^{\nComp n}$
    \begin{algorithmic}[1]
\State{$\resD \gets \rhsDParam - \StiffMat\, \solD\cgIdxFirst_{\param}$} \algorithmiccomment{Initial residual}
\State{$\cgSearchDir \gets \ul{0} \in \ffR^{\nComp n}$; $\cgDelta_1 \gets 1 \in \ffR$} \algorithmiccomment{Initialization}
	\While{ $\normg{\resD}_\infty > \cgTol $ }
        \State{Reshape $\resD$ to $\field{\res}$; $\hat{\resD} \gets \vecOp{\trafo \left( \field{\res} \right) }$} \algorithmiccomment{\fboxc{uniSlightblue}{$\orderOf{n \log n}$}}
		\State{$\hat{\precResD} \gets \PrecSparse \, \hat{\resD}$} \hspace*{14em}
        \rlap{\smash{$\left.\begin{array}{@{}c@{}}\\{}\\{}\\{}\end{array}\color{uniSlightblue}\right\}
          \color{uniSlightblue}\begin{tabular}{l}$\text{\UNO{} preconditioner}$\end{tabular}$}}\algorithmiccomment{\fboxc{uniSlightblue}{$\orderOf{n}$} since $\PrecSparse$ is a band matrix}
		\State{Reshape $\hat{\precResD}$ to $\field{\hat{\precRes}}$; $\precResD \gets \vecOp{\trafoInv \left( \field{\hat{\precRes}} \right) }$} \algorithmiccomment{\fboxc{uniSlightblue}{$\orderOf{n \log n}$}}
		\State{$\cgDelta_0 \gets \cgDelta_1; \; \cgDelta_1 \gets \resD \cdot \precResD$}
        
		\State{$\cgSearchDir \gets \precResD + \frac{\cgDelta_1}{\cgDelta_0} \, \cgSearchDir$} \algorithmiccomment{Update search direction}
		\State{$\cgMatvecResult \gets \StiffMat\, \cgSearchDir$} \algorithmiccomment{Matrix-vector product}
		\State{$\cgStepWidth \gets \cgDelta_1 \, / \left(\cgSearchDir \cdot \cgMatvecResult \right)$} \algorithmiccomment{Compute optimal step width}
		\State{$\resD \gets \resD - \cgStepWidth \cgMatvecResult$} \algorithmiccomment{Update residual}
		\State{$\solDParam \gets \solDParam + \cgStepWidth \cgSearchDir$} \algorithmiccomment{Update solution approximation}
    \EndWhile
    \end{algorithmic}
\end{algorithm}

\begin{remark*}
There remains a critical aspect that must be taken into account.
If~$\PrecD_{\weights}$ exhibits eigenvalues close to zero---which is not prevented by construction---the convergence of the hybrid solver in floating-point arithmetic may fail.
Hence, it is advised to perform the following safety check with a given tolerance $\epsilon > 0$ (e.g., in the same order of magnitude as the \emph{machine epsilon}):
\begin{align} \label{eq:fnocg-sanity-check}
    0< \epsilon < \lambda_1 \leq \lambda_2 \leq \dots \leq \lambda_{\nComp n} \,, &&
	\bigcup\limits_{i=1}^{n} \eig{\fnocgFundSolMat\modeIdx} = \{ \lambda_1, \dots, \lambda_{\nComp n}\} \,.
\end{align}
This criterion can be checked {\color{revision}in~$\orderOf{n}$ by computing the eigenvalues of the symmetric matrices~$\fnocgFundSolMat\modeIdx \in \ffR^{c \times c}$ for~$1 \leq i \leq n$, or evaluating Sylvester's criterion on ~$\fnocgFundSolMat\modeIdx - \epsilon \ull{I} \in \ffR^{c \times c}$, which is computationally inexpensive and straightforward to parallelize}.
\end{remark*}

\begin{remark*}
When using periodic boundary conditions or Neumann boundary conditions, the stiffness matrix can exhibit a null space, i.e., it is only positive semi-definite.
As long as the right-hand side vector is in the range or image of the stiffness matrix, the (preconditioned) \CG{} solver still converges to the solution.
In this case, the \UNO{} preconditioner $\PrecD_{\weights}$ is allowed to share the null space of the stiffness matrix (i.e, $\PrecD_{\weights}$ is allowed to have a certain number of zero eigenvalues), and convergence of \UNOCG{} can still be proven similarly to \cite{Kaasschieter1988}.
\end{remark*}

\subsection{Naive first-order training procedure for \UNOCG}
\label{ssec:fnocg-naive-training}

What remains to be explained is how the \UNO{} preconditioner~$\Prec_{\weights}$ can be trained, i.e., how its learnable weights~$\weights$ can be determined on the basis of training data.
For this purpose, a naive method similar to the standard techniques in scientific machine learning is presented first.
In a discrete setting, the training data from \cref{ssec:mlprec} that is given by the continuous fields~$\PrecOutSample$ and~$\PrecInSample$ corresponds to vectors that contain their DOF according to a given discretization as in
\begin{align}
	\PrecOutDSample = \vecOp{\field{ \PrecOutSample }} \in \ffR^{\nComp n} \,, &&
	\PrecInDSample = \vecOp{\field{ \PrecInSample }} \in \ffR^{\nComp n} \,, &&
	1 \leq \sampleI \leq \nSamples \,.
\end{align}
Then, the discrete operator $\PrecD_{\weights}$ of the machine-learned preconditioner~$\Prec_{\weights}$ on a fixed regular grid with~$n$ nodes is trained by learning the relation between these discrete representations.
For that, we consider the given training data in the form of sets of vectors
\begin{align} \label{eq:prec-train-data-discr}
	\PrecOutDSet = \left\{ \PrecOutDSample: 1 \leq \sampleI \leq \nSamples \right\} \subsetneq \ffR^{\nComp  n} \,, &&
	\PrecInDSet = \left\{ \PrecInDSample: 1 \leq \sampleI \leq \nSamples \right\} \subsetneq \ffR^{\nComp  n} \,,
\end{align}
 that satisfy the linear system of the discretized \ParPDE{} from \eqref{eq:fem-linear-system} for given parameters, i.e.,
\begin{align}
    \StiffMatParamSample \, \PrecOutDSample = \PrecInDSample \,, && \param\sampleIdx \in \paramSpace \,, && 1 \leq \sampleI \leq \nSamples \,,
\end{align}
where $\StiffMatParamSample$ denotes the stiffness matrix of the discretized \ParPDE{} for different chosen parameters~$\param\sampleIdx \in \paramSpace$.
As a specific realization of the loss function in \eqref{eq:prec-loss}, we use the MSE loss that is often used for regression problems in scientific machine learning on the DOF vectors, which leads to
\begin{align} \label{eq:fnocg-loss}
	\!\min_\weights \quad \Loss \,, &&
	\text{where} && 
	\Loss = \frac{1}{\nSamples} \sum_{\sampleI = 1}^\nSamples \lossSample \,, &&
	\text{with} &&
	\lossSample = \ell \left( \PrecD_\weights \, \PrecInDSample ,\, \PrecOutDSample \right) = \normg{ \PrecD_\weights \, \PrecInDSample - \PrecOutDSample }^2 \,.
\end{align}
In practice, the machine-learned preconditioner works best in a \CG{} solver if the vectors $\PrecInDSample$, $\PrecOutDSample$ are close to the vectors that actually occur in the solver, see \cref{alg:pcg}.
Thus, a sensible approach is often to use
\begin{align}
\PrecInDSample = \rhsDParamSample &&
\text{and} &&
\PrecOutDSample = \StiffMat^{-1}_{\param\sampleIdx} \, \rhsDParamSample = \solDParamSample \,, &&
1 \leq \sampleI \leq \nSamples \,,
\end{align} as part of the training data sets $\PrecInDSet$ and $\PrecOutDSet$.
Here,~$\rhsDParam\sampleIdx$ is the actual right-hand side of the linear system to solve, that is, $\StiffMatParamSample \, \solD\sampleIdx = \rhsDParam\sampleIdx$.
In order to gain a \UNO{} preconditioner that works on an entire problem class of stiffness matrices $\StiffMatParamSample$, many different parameters~$\param\sampleIdx \in \paramSpace$ for the \ParPDEs{} should be considered in the training data.
For the naive approach of preconditioner training, we apply the standard training procedure for machine learning models based on gradient descent and its variants (SGD \cite{Saad1999}, Adam \cite{Adam}, etc.).
This requires derivatives of the evaluated loss with respect to the learnable weights as in
\begin{align}
	\optimGrad = \frac{\partial \Loss}{\partial \weights} \in \ffR^{\nWeights} \,.
\end{align}
Here, backpropagation through the preconditioner, including the transformation~$\trafo${\color{revision},} is used to determine these derivatives in each epoch.
However, this is memory exhaustive and computationally expensive, having a time complexity of~$\orderOf{\nSamples n \log n}$.
In addition, these first-order optimization schemes suffer from poor convergence, especially when the global optimum has to be reached with high accuracy.
The just introduced naive approach is formalized in \cref{alg:fnocg-training-naive}.

\begin{algorithm}
    \caption{Naive first-order \UNO{} preconditioner training}
	\label{alg:fnocg-training-naive}
    \hspace*{\algorithmicindent} \textbf{Data:} sets $\PrecInDSet$, $\PrecOutDSet$ with $\PrecOutDSample, \PrecInDSample \in \ffR^{\nComp n}$ for $\param\sampleIdx \in \paramSpace$, $1 \leq j \leq \nSamples$; parametrization $\globalParam$; $\nEpochs \in \ffN$ \\
    \hspace*{\algorithmicindent} \textbf{Result:} $\PrecSparse \in \ffR^{\nComp n \times \nComp n}$ (sparse matrix with $\nComp^2 n$ nonzero entries)
    \begin{algorithmic}[1]
\State{Load data $\PrecOutDSample, \PrecInDSample \in \ffR^{\nComp n}$, $1 \leq j \leq \nSamples$}
\State{Initialize weights $\weights \in \ffR^{\nWeights}$}
	\For{ $1 \leq e \leq \nEpochs $ }
		\State{Compute loss $\Loss$ on data $\PrecOutDSample, \PrecInDSample \in \ffR^{\nComp n}$, $1 \leq j \leq \nSamples$ } \algorithmiccomment{\fboxc{uniSlightblue}{$\orderOf{\nSamples n \log n}$}}
		\State{Compute gradient $\optimGrad = \frac{\partial \Loss}{\partial \weights}$ using backpropagation} \algorithmiccomment{\fboxc{uniSlightblue}{$\orderOf{\nSamples n \log n}$}}
		\State{Update weights $\weights \gets \weights - \learningRate \optimGrad$} \algorithmiccomment{with a learning rate $\learningRate > 0$}
    \EndFor
		\State{Assemble $\PrecSparse \in \ffR^{\nComp n \times \nComp n}$ based on $\PrecEntries = \globalParam \left( \weights \right) \in \ffR^{\nPrecC n}$}
    \end{algorithmic}
\end{algorithm}

\subsection{Efficient second-order training procedure for \texorpdfstring{\UNOCG}{UNO-CG}}
\label{ssec:fnocg-fasttraining}

As outlined in \cref{ssec:fnocg-naive-training}, the standard training procedure, consisting of backpropagation and first-order optimization methods such as gradient descent, is not efficient and suffers from poor convergence.
By exploiting the special advantageous structure of the \UNO{} preconditioner on an algebraic level as analyzed in \cref{sssec:fnocg-analysis}, the training procedure can be drastically improved.
In fact, it is possible to derive analytical expressions for the gradients and the Hessian with respect to the learnable weights that can be evaluated in~$\orderOf{n}$ after preprocessing the training data.
This gives rise to a Newton-Raphson method for the preconditioner training that converges quadratically to the global optimum as long as the initial guess is good enough.
In contrast to the training procedure typically used for neural operators, the method derived in the following is {\color{revision}sufficiently fast in practice} that it does not require the use of high-end GPUs to achieve short training times---even for large-scale high-dimensional problems.

We achieve this by interpreting the optimization problem from \eqref{eq:fnocg-loss} as a nonlinear least squares problem, {\color{revision}where the nonlinearity with respect to the weights~$\weights$ arises from the local parametrization $\localParam$ introduced in \cref{sssec:fnocg-parametrization}.}
To this end, we reformulate the sample-wise loss term by recalling the definition of $\PrecD_\weights$ from \eqref{eq:FNOCG-prec-fields-disc} as
\begin{align*}
    \lossSample = \normg{\PrecD_\weights \, \PrecInDSample - \PrecOutDSample}^2
    &= \Her{\left( \PrecD_\weights \, \PrecInDSample - \PrecOutDSample \right)} \left( \PrecD_\weights \, \PrecInDSample - \PrecOutDSample \right) \\
    &= \Her{\left( \PrecInDSample \right)} \trafoDHer \PrecSparse \, \trafoD \, \trafoDHer \PrecSparse \, \trafoD \, \PrecInDSample - \Her{\left( \PrecInDSample \right)} \trafoDHer \PrecSparse \, \trafoD \, \PrecOutDSample - \Her{\left( \PrecOutDSample \right)} \trafoDHer \PrecSparse \, \trafoD \, \PrecInDSample + \normg{\PrecOutDSample}^2 \\
    &= \Her{\left( \trafoD \, \PrecInDSample \right)} \PrecSparse \, \PrecSparse \left( \trafoD \, \PrecInDSample \right) - \Her{\left( \trafoD \, \PrecInDSample \right)} \PrecSparse \left( \trafoD \, \PrecOutDSample \right) - \Her{\left( \trafoD \, \PrecOutDSample \right)} \PrecSparse \left( \trafoD \, \PrecInDSample \right) + \normg{\PrecOutDSample}^2
    .
\end{align*}
Herein, the properties of the real symmetric matrix~$\PrecSparse$, and the unitary matrix~$\trafoD$ corresponding to the unitary transform~$\trafo$ are exploited.
In this result, the transformed samples of the training data
\begin{align}
    \PrecInDFwdSample = \trafoD \, \PrecInDSample = \vecOp{ \trafo \left( \field{ \PrecInSample } \right) } \,, &&
    \PrecOutDFwdSample = \trafoD \, \PrecOutDSample = \vecOp{ \trafo \left( \field{ \PrecOutSample } \right) } \,, && 1 \leq \sampleI \leq \nSamples \,.
\end{align}
can be precomputed in a one-time effort after data generation and before training.
This is possible with a computational effort of~$\orderOf{\nSamples n \log n}$ under the assumption that fast algorithms for the unitary transform~$\trafo$ exist, similar to the FFT.
These quantities are independent of the model's weights and hyperparameters and can hence be reused in every epoch of the optimization scheme.
Furthermore, a scalar quantity~$\featureDelta$ in the loss function~$\Loss$ can be precomputed via~$\featureDelta = \frac{1}{\nSamples} \sum_{\sampleI = 1}^\nSamples \featureDeltaSample$, where~$\featureDeltaSample = \normg{\PrecOutDSample}^2$.
Based on that, the loss function can be simplified to~$\Loss = \sum_{\sampleI = 1}^\nSamples \lossSample$ with the sample-wise loss term
\begin{align} \label{eq:loss-precomputed}
	\lossSample = \Her{\left( \PrecInDFwdSample \right)} \PrecSparse \, \PrecSparse \, \PrecInDFwdSample - \Her{\left( \PrecInDFwdSample \right)} \PrecSparse \, \PrecOutDFwdSample - \Her{\left( \PrecOutDFwdSample \right)} \PrecSparse \, \PrecInDFwdSample + \featureDeltaSample \,.
\end{align}
Since~$\PrecSparse$ is a sparse band matrix with in total~$\nComp^2 n$ nonzero entries, the computational effort to evaluate the expression for the loss function in \eqref{eq:loss-precomputed} is in~$\orderOf{n}$.
Further, derivatives of the loss function with respect to the symmetric matrix~$\PrecSparse$ are also easily available.
In the following, special cases are presented that account for a special structure of~$\PrecSparse$.
These cases build up on each other and are sufficient for all problems tackled in \cref{sec:problems}, but can also be extended easily.

\paragraph*{Special case of one nodal DOF (\texorpdfstring{$\nComp=1$}{\nComp=1})}
\label{sssec:special1dof}

The first case of having only one nodal DOF is, for instance, given for thermal homogenization problems.
If the \ParPDEs{} to solve is only scalar-valued, this reduces to determining the nonzero entries~$\PrecEntries \in \ffR^{n}$ of~$\PrecSparse = \diag{\PrecEntries}$ as introduced in \cref{sssec:fnocg-analysis}.
Following this, the loss~$\lossSample$ for the~$\sampleI$-th sample in~\eqref{eq:loss-precomputed} can be expressed as
\begin{align}
	\lossSample =
	\left( \left(\PrecInDFwdSample\right)^* \odot \PrecInDFwdSample \right) \cdot (\PrecEntries \odot \PrecEntries) - 2 \, \realPart{ \left(\PrecInDFwdSample\right)^* \odot \PrecOutDFwdSample } \cdot \PrecEntries + \featureDeltaSample \,.
\end{align}
This result allows for to direct computation of the loss~$\Loss$ as well as its first-order and second-order derivatives with respect to the nonzero entries~$\PrecEntries$ with a computational effort in~$\orderOf{n}$ via
\begin{gather}
	\Loss = \frac{1}{\nSamples} \sum_{\sampleI = 1}^\nSamples \lossSample = \featureA_1 \cdot \left( \PrecEntries \odot \PrecEntries \right) - \featureB_1 \cdot \PrecEntries + \featureDelta \in \ffR \,, \\
	\frac{\partial \Loss}{\partial \PrecEntries} = \frac{1}{\nSamples} \sum_{\sampleI = 1}^\nSamples \frac{\partial \lossSample}{\partial \PrecEntries} = 2 \featureA_1 \odot \PrecEntries - \featureB_1 \in \ffR^n \,, \qquad\qquad
	\frac{\partial^2 \Loss}{\partial \PrecEntries \partial \PrecEntries} = \frac{1}{\nSamples} \sum_{\sampleI = 1}^\nSamples \frac{\partial^2 \lossSample}{\partial \PrecEntries\partial \PrecEntries} = 2 \diag{\featureA_1} \in \ffR^{n \times n} \,,
\end{gather}
by using the precomputed features~$\featureA_1, \featureB_1 \in \ffR^n$ that are defined as
\begin{align} \label{eq:features-c-1}
    \featureA_1 = \frac{1}{\nSamples} \sum_{\sampleI = 1}^\nSamples \left( \PrecInDFwdSample \right)^* \odot \PrecInDFwdSample \,, &&
    \featureB_1 = \frac{2}{\nSamples} \sum_{\sampleI = 1}^\nSamples \realPart{ \left( \PrecInDFwdSample \right)^* \odot \PrecOutDFwdSample } \,.
\end{align}

\paragraph*{Special case of two nodal DOF (\texorpdfstring{$\nComp=2$}{\nComp=2})}
\label{sssec:special2dof}

For mechanical homogenization problems in~$d=2$ dimensions, $\nComp=2$ nodal DOF are present, which requires a slightly more complicated procedure due to $\nPrecC = \nComp (\nComp+1)/2=3$.
Therein, we introduce a split of the DOF vectors~$\PrecInDSample, \PrecOutDSample \in \ffR^{2n}$ into vectors containing the corresponding nodal DOF~$\PrecInDSample_{\mathrm{x}}, \PrecInDSample_{\mathrm{y}}, \PrecOutDSample_{\mathrm{x}}, \PrecOutDSample_{\mathrm{y}} \in \ffR^{n}$ analogously to
\begin{align}
\label{eq:dof-split-2}
\begin{split}
    \PrecInDSample = 
    \T{\begin{bmatrix} \left(\PrecInDSample_{\mathrm{x}}\right)_1 & \left(\PrecInDSample_{\mathrm{y}}\right)_1 & \cdots & \left(\PrecInDSample_{\mathrm{x}}\right)_n & \left(\PrecInDSample_{\mathrm{y}}\right)_n \end{bmatrix}} \in \ffR^{2n} \,.
\end{split}
\end{align}
In this case, the problem reduces to learning a symmetric band matrix that is parametrized via
\begin{align}
	\PrecSparse =
	\begin{bmatrix}
		\diag{\PrecEntries_1} & \diag{\PrecEntries_3} \\
		\diag{\PrecEntries_3} & \diag{\PrecEntries_2}
	\end{bmatrix} \,,
	&& \PrecEntries = \T{\begin{bmatrix} \T{\PrecEntries}_1 \quad \T{\PrecEntries}_2 \quad \T{\PrecEntries}_3 \end{bmatrix}} \in \ffR^{3n} \,.
\end{align}
Based on that, the loss~$\Loss$ over all samples as well as its first-order and second-order derivatives with respect to the nonzero entries~$\PrecEntries \in \ffR^{3 n}$ of~$\PrecSparse$ can again be computed in~$\orderOf{n}$ using
\begin{gather}
	\Loss =
	\featureA_1 \cdot \left( \PrecEntries_1 \odot \PrecEntries_1 + \PrecEntries_3 \odot \PrecEntries_3 \right)
	+ \featureA_2 \cdot \left( \PrecEntries_2 \odot \PrecEntries_2 + \PrecEntries_3 \odot \PrecEntries_3 \right)
	+ \featureA_3 \cdot \left( \PrecEntries_1 \odot \PrecEntries_3 + \PrecEntries_2 \odot \PrecEntries_3 \right)
	- \featureB_1 \cdot \PrecEntries_1
	- \featureB_2 \cdot \PrecEntries_2
	- \featureB_3 \cdot \PrecEntries_3
	+ \featureDelta \,, \\
	\frac{\partial \Loss}{\partial \PrecEntries} = \T{\begin{bmatrix}
		2 \featureA_1 \odot \PrecEntries_1
		+ \featureA_3 \odot \PrecEntries_3
		- \featureB_1 \\
		2 \featureA_2 \odot \PrecEntries_2
		+ \featureA_3 \odot \PrecEntries_3
		- \featureB_2 \\
		2 \featureA_1 \odot \PrecEntries_3
		+ 2 \featureA_2 \odot \PrecEntries_3
		+ \featureA_3 \odot \left( \PrecEntries_1 + \PrecEntries_2 \right)
		- \featureB_3
	\end{bmatrix}} \,, \quad
	\frac{\partial^2 \Loss}{\partial \PrecEntries\partial \PrecEntries} = \begin{bmatrix}
		2 \diag{\featureA_1} & \ull{0} & \diag{\featureA_3} \\
		\ull{0} & 2 \diag{\featureA_2} & \diag{\featureA_3} \\
		\diag{\featureA_3} & \diag{\featureA_3} & 2 \diag{\featureA_1 + \featureA_2}
	\end{bmatrix} \,,
\end{gather}
with the precomputed features $\featureA_1, \featureA_2, \featureA_3, \featureB_1, \featureB_2, \featureB_3 \in \ffR^n$ that are defined as
\begin{align*} \label{eq:features-c-2}
    \featureA_1 &= \frac{1}{\nSamples} \sum_{\sampleI = 1}^\nSamples \PrecInDAdjXSample \odot \PrecInDFwdXSample \,, &&
    \featureA_2 = \frac{1}{\nSamples} \sum_{\sampleI = 1}^\nSamples \PrecInDAdjYSample \odot \PrecInDFwdYSample \,, &&
    \featureA_3 = \frac{1}{\nSamples} \sum_{\sampleI = 1}^\nSamples \PrecInDAdjXSample \odot \PrecInDFwdYSample + \PrecInDAdjYSample \odot \PrecInDFwdXSample \,, \\
    \featureB_1 &= \frac{2}{\nSamples} \sum_{\sampleI = 1}^\nSamples \realPart{ \PrecInDAdjXSample \odot \PrecOutDFwdXSample } \,, &&
    \featureB_2 = \frac{2}{\nSamples} \sum_{\sampleI = 1}^\nSamples \realPart{ \PrecInDAdjYSample \odot \PrecOutDFwdYSample } \,, &&
    \featureB_3 = \frac{2}{\nSamples} \sum_{\sampleI = 1}^\nSamples \realPart{ \PrecInDAdjXSample \odot \PrecOutDFwdYSample + \PrecInDAdjYSample \odot \PrecOutDFwdXSample } \,.
\end{align*}

\paragraph*{Special case of three nodal DOF (\texorpdfstring{$\nComp=3$}{\nComp=3})}

For mechanical homogenization problems in~$d=3$ dimensions,~$\nComp=3$ nodal DOF are involved.
A similar procedure for computing the loss and its derivatives in~$\orderOf{n}$ is available in \cref{sec:training-c-3}.

\paragraph*{Newton-Raphson method}
In all these cases, the gradient~$\optimGrad$ and the Hessian~$\optimHess = \T{\optimHess}$ of the loss~$\Loss$ with respect to the learnable weights~$\weights$ are available based on the parametrization~$\globalParam$ via the chain rule in a straightforward way as in
\begin{align}
    \optimGrad = \frac{\partial \Loss}{\partial \weights} = \frac{\partial \Loss}{\partial \PrecEntries} \frac{\partial \PrecEntries}{\partial \weights} = \frac{\partial \Loss}{\partial \PrecEntries} \frac{\partial \globalParam(\weights)}{\partial \weights} \,, &&
    \optimHess = \frac{\partial^2 \Loss}{\partial \weights \partial \weights} = \frac{\partial^2 \Loss}{\partial \PrecEntries \partial \weights} \frac{\partial \PrecEntries}{\partial \weights} + \frac{\partial \Loss}{\partial \PrecEntries} \frac{\partial^2 \globalParam(\weights)}{\partial \weights \partial \weights} \,.
\end{align}
The explicit gradient and Hessian evaluation in the proposed \UNOCG{} training enables a straightforward implementation of a Newton-Raphson scheme for the training, yielding second-order convergence and \emph{significant} savings over the naive approach described in \cref{ssec:fnocg-naive-training}.
Before training, the features~$\featureA\featureIdx$ and~$\featureB\featureIdx$ for~$1 \leq \featureI \leq \nPrecC$ have to be computed in a preprocessing step for the available training data at a time complexity of~$\orderOf{\nSamples n \log n}$.
The actual training procedure for the~\UNO{} preconditioner is performed with a time complexity of~$\orderOf{n}$ as it is presented in \cref{alg:fnocg-newton-training-fast}.
Note that the unitary transform~$\trafo$, which represents the computationally expensive part of the \UNO{} preconditioner, only has to be applied once to each training sample, but \emph{never} during the iterative training procedure.
It is also worth noting that this condenses the training data consisting of~$\nSamples \nComp n$ values into a reduced representation of~$\nPrecC n \ll \nSamples \nComp n$ values without losing information relevant for the \UNO{} preconditioner training, facilitating data handling and reducing memory requirements massively (by a factor of~$\frac{2}{\nComp+1}\nSamples$).

\begin{algorithm}
    \caption{Second-order \UNO{} preconditioner training with linear time complexity}
\label{alg:fnocg-newton-training-fast}
    \hspace*{\algorithmicindent} \textbf{Data:} features $\featureA\featureIdx, \featureB\featureIdx \in \ffR^n$, $1 \leq m \leq \nPrecC $; parametrization $\globalParam$; $\nEpochs \in \ffN$ \\
    \hspace*{\algorithmicindent} \textbf{Result:} $\PrecSparse \in \ffR^{\nComp n \times \nComp n}$ (sparse matrix with $\nComp^2 n$ nonzero entries)
    \begin{algorithmic}[1]
\State{Load features $\featureA\featureIdx, \featureB\featureIdx \in \ffR^n$, $1 \leq m \leq \nPrecC $}
\State{Initialize weights $\weights \in \ffR^{\nWeights}$}
	\For{ $1 \leq e \leq \nEpochs $ }
		\State{Compute gradient $\optimGrad = \frac{\partial \Loss}{\partial \weights}$ using features $\featureA\featureIdx, \featureB\featureIdx \in \ffR^n$, $1 \leq m \leq \nPrecC$} \algorithmiccomment{\fboxc{uniSlightblue}{$\orderOf{n}$}}
		\State{Compute Hessian $\optimHess = \frac{\partial^2 \Loss}{\partial \weights \partial \weights}$ using features $\featureA\featureIdx, \featureB\featureIdx \in \ffR^n$, $1 \leq m \leq \nPrecC $} \algorithmiccomment{\fboxc{uniSlightblue}{$\orderOf{n}$}}
	\State{Compute Newton step $\optimStep = \optimHess^{-1} \optimGrad$} \algorithmiccomment{\fboxc{uniSlightblue}{$\orderOf{n}$} since $\optimHess$ is a band matrix}
    \State{Update weights $\weights \gets \weights + \optimStep$}
    \EndFor
		\State{Assemble $\PrecSparse \in \ffR^{\nComp n \times \nComp n}$ based on $\PrecEntries = \globalParam \left( \weights \right) \in \ffR^{\nPrecC n}$}
    \end{algorithmic}
\end{algorithm}

\paragraph*{\color{revision}Amortization of training costs}

{\color{revision}
The generation of training data~$\PrecInDSet$, $\PrecOutDSet$ and the computation of features~$\featureA\featureIdx$, $\featureB\featureIdx$ for the \UNO{} preconditioner incur an initial cost. However, this cost can be quickly amortized in many-query scenarios comprising, e.g., multi-scale simulations and materials design.
Once trained, the preconditioner can be reused for multiple simulations with different parameters without retraining. Additionally, even a preconditioner trained on a small number of samples can be used to generate further high-quality training data due to the guaranteed convergence of \UNOCG{} to continuously improve the learned \UNO{} preconditioner. Moreover, as discussed in \cref{sssec:architecture}, the \UNO{} preconditioner operates in a transformed space, making it generally resolution-invariant. Therefore, a preconditioner trained on data at a given resolution could also be applied to problems with different resolutions. This highlights the trade-off between an initial training effort and long-term computational benefits in many-query scenarios.
}

\subsection{GPU-boosted matrix-free implementation}
\label{ssec:fnocg-implementation}

We provide an open-source implementation of the \UNOCG{} hybrid solver as a GPU-accelerated solver based on PyTorch \cite{PyTorch,PyTorch2}.
It enables the batched solution of parametric PDEs using the \CG{} method with different preconditioners and is available in \cite{UNOCGgithub2025}.
In particular, matrix-free implementations of the preconditioned \CG{} method enable excellent scaling when solving parametric PDEs on modern GPU architectures in a many-query scenario.
The learned \UNO{} preconditioner can also be used in CPU-based high-performance computing (HPC) environments.
For this purpose, we provide an efficient CPU-based implementation based on PETSc \cite{Petsc1997} and FFTW \cite{FFTW05} in our software repository \cite{UNOCGgithub2025} that can, e.g., be used to solve parametric PDEs using FEniCS \cite{LoggWells2010} or Firedrake \cite{FiredrakeUserManual}.
In particular, we use these frameworks to validate our GPU-accelerated implementation.

\section{Problem formulations}
\label{sec:problems}

\begin{figure}[ht!]
    \centering
    \begin{subfigure}[b]{0.495\textwidth}
            \includegraphics[scale=0.49]{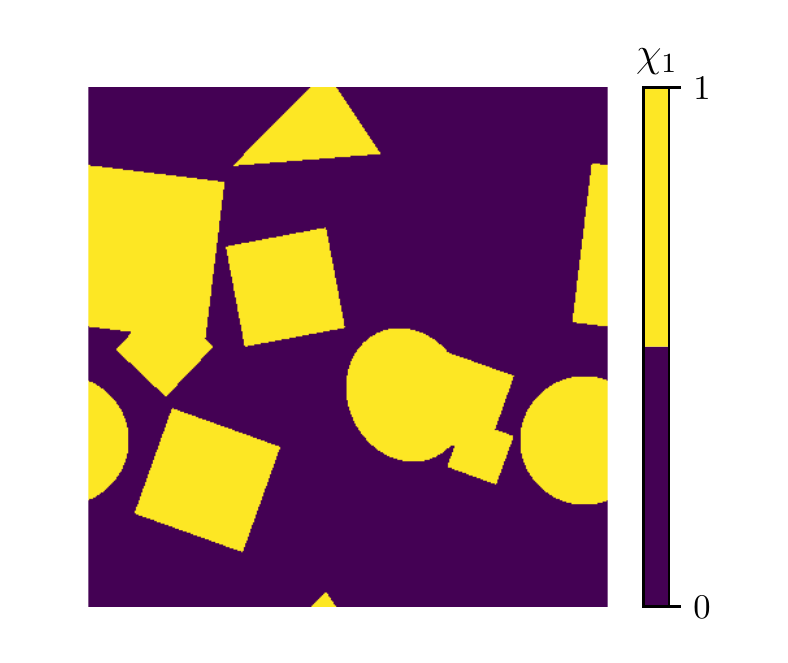}
            \includegraphics[scale=0.49]{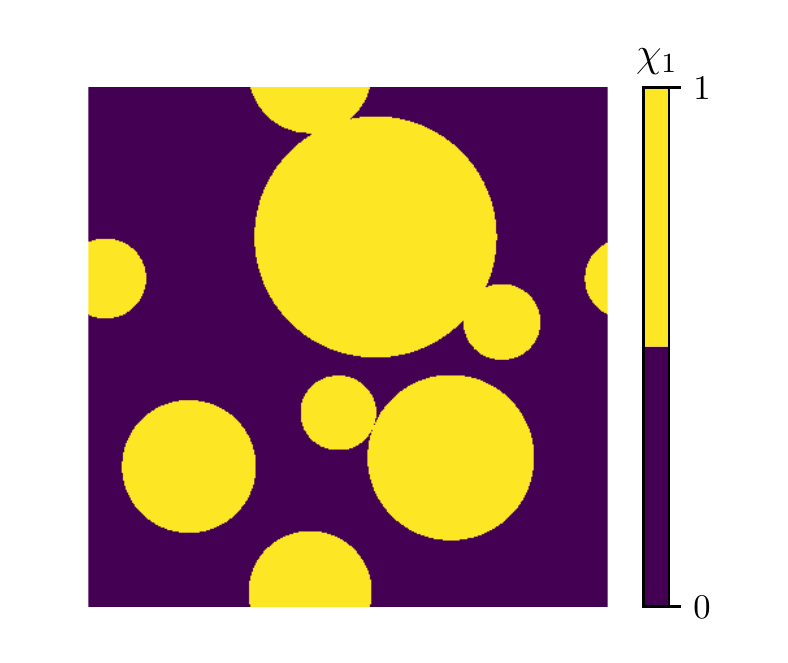}
        \vspace{-0.5\baselineskip}
        \caption{Examples of 2D microstructures ($400^2$ pixels).}
        \label{fig:microstructure-2d}
    \end{subfigure}
    \hfill
    \begin{subfigure}[b]{0.495\textwidth}
        \begin{center}
            \includegraphics[scale=0.14]{fig/ms3d\_4.png}
            \includegraphics[scale=0.14]{fig/ms3d\_1.png}
        \end{center}
        \vspace{-0.5\baselineskip}
        \caption{Examples of 3D microstructures ($192^3$ voxels).}
        \label{fig:microstructure-3d}
    \end{subfigure}
    \vspace{-0.5\baselineskip}
    \caption{\normalsize Examples of microstructures that we consider for homogenization problems.}
    \label{fig:microstructures}
\end{figure}

While the framework of \UNOCG{} is designed for general \ParPDEs{} as introduced in \cref{ssec:parPDE}, we apply it here to various homogenization problems, a common application scenario for which the \UNOCG{} hybrid solver can be benchmarked against existing, specialized solvers.
In computational homogenization, the overall goal is to determine the effective material behavior of a heterogeneous material based on a given microstructure using numerical simulations.
For that, the microstructure is assumed to be a periodic continuation of a representative volume element (RVE) with the domain~$\domain \subset \ffR^\nDims$.
The microscopic position in the RVE is denoted by~$\fx \in \domain$.
We consider homogenization problems on both a two-dimensional and a three-dimensional RVE, i.e.,~$\nDims \in \{2, 3\}$.
The two-dimensional problems are more convenient for visualization and analysis, whereas an application to three-dimensional problems shows the scalability towards real-world problems of our hybrid~\UNOCG{} solver.
The domain of the RVE~$\domain = [-\rveLen/2, \rveLen/2]^\nDims$ ($\rveLen>0$) has boundaries
\begin{align}
    \edgesI = \left\{ \fx \in \partial\domain \; \Big\vert \; \fx \cdot \unitVecI = \pm \frac{\rveLen}{2} \right\} \,,
    && 1 \leq \dimI \leq \nDims \,,
    && \edges = \bigcup\limits_{\dimI=1}^{\nDims} \edgesI \,.
\end{align}
Each RVE is assumed to consist of two phases.
Hence, the domain $\domain$ is decomposed into subdomains (i.e., material phases)~$\domain_0$, $\domain_1$ with corresponding indicator functions $\chi_0, \chi_1: \domain \to \{0, 1\}$.
The assumption of two phases is only made for simplicity and does not represent a restriction of the presented methods.
Often, the microstructures for homogenization problems are given through images, i.e., each pixel or voxel of an image defines the assignment of the corresponding area or volume of the microstructure to one of the two phases.
Hence, it is convenient to use the underlying regular grid of these images as a discretization of the domain~$\domain$.
On the one hand, we consider data sets of two-dimensional microstructures with a resolution of~$400 \times 400$ pixels that are available via~\cite{Darus} and have been used in \cite{Lissner2019, Lissner2023}.
On the other hand, we will use a data set of three-dimensional microstructures with a resolution of~$192 \times 192 \times 192$ voxels (that is $\approx 7.1 \cdot 10^6$ voxels in total) that are available in~\cite{Prifling2021} and lead to large-scale FEM problems.
{\color{revision}While these microstructures serve as proof-of-concept examples, this choice is primarily driven by the availability of a large number of samples. The memory-efficient training of \UNOCG{} also allows an application to larger microstructures.}
Examples from these microstructure data sets are shown in \cref{fig:microstructures}.
For all homogenization problems presented in the following, the discretization is performed using linear FEM on a regular grid as introduced in \cref{ssec:discretization}.

\subsection{Homogenization problems of linear heat conduction}
\label{ssec:thermal-problems}

As the first type of homogenization problems, we consider thermal problems to compute the effective thermal conductivity of microstructured materials.
Let~$\temperature(\fx)$ be the temperature for~$\fx \in \domain$. The temperature gradient~$\tempGrad$ and the heat flux $\flux$ are related by Fourier's law as in
\begin{align}
    \tempGrad(\fx) &= \temperature(\fx) \otimes \nabla = \pd{\temperature(\fx)}{x_i} \, \unitVecI, &
    \flux(\fx) &= - \heatCondTensor (\fx) \tempGrad(\fx) \,,
\end{align}
with a conductivity tensor~$\heatCondTensor \in \spdSpace{{\color{revision}\ffR^{\nDims \times \nDims}}}$.
The steady-state heat conduction problem then reads
\begin{align}
    \flux(\fx) \cdot \nabla &= 0 \, .
\end{align}
The heat conductivity tensor~$\heatCondTensor = \left( \heatCond_0 \chi_0(\fx) + \heatCond_1 \chi_1(\fx) \right) \fI$ is phasewise constant and assumed isotropic with conductivities $\heatCond_0$, $\heatCond_1$ for each phase.
The assumption of isotropy is only made for simplicity and does not represent a restriction on the presented algorithms.
For homogenization problems, it is convenient to decompose the primary variable into parts that are either related to macroscopic or microscopic quantities.
For thermal homogenization problems, we decompose the temperature field into
\begin{align} \label{eq:thermal-decomposition}
	\temperature(\fx) = \tempMacro + \tempGradMacro \cdot \fx + \tempFluct(\fx) \,, && \fx \in \domain \,,
\end{align}
with the macroscopic temperature~$\tempMacro$, the temperature field~$\tempGradMacro \cdot \fx$ that is induced by a prescribed macroscopic temperature gradient~$\tempGradMacro$, and the temperature fluctuation field~$\tempFluct(\fx)$.
The latter satisfies the zero-mean property, i.e.,~$\int_{\domain} \tempFluct(\fx) \dInt{\fx} = 0$.
As boundary conditions for the fluctuation field~$\tempFluct$, we consider the possibilities
\begin{itemize}
    \item periodic boundary conditions, i.e.,
    $\displaystyle    \onEdgeP{ \tempFluct } = \onEdgeM{ \tempFluct } \,$,
    \item and Dirichlet boundary conditions, i.e.,
    $\displaystyle
        \onBoundary{ \tempFluct } = 0 \,$.
\end{itemize}
The differential operator of the linear heat conduction homogenization problem presented in the general framework of \ParPDEs{} from \cref{ssec:parPDE} and the right-hand side~$\rhsParam$ depend on parameters~$\param = \left\{ \param_1, \param_2 \right\}$ that are given by
\begin{align} \label{eq:operator-thermal}
	\diffOpParam(\bullet) = \left( -\left(\paramSymbol_1(\fx) \fI\right) \left( \bullet \otimes \nabla \right) \right) \cdot \nabla \,, &&
    \param_1 = \heatCond_0 \chi_0(\fx) + \heatCond_1 \chi_1(\fx) ,\, && \param_2 = \tempGradMacro \,.
\end{align}

\subsection{Homogenization problems of linear elasticity}
\label{ssec:mechanical-problems}

In addition to thermal homogenization problems, we consider mechanical homogenization problems to determine the effective stiffness of microstructured materials.
Again, two-dimensional ($d=2$) and three-dimensional ($d=3$) microstructures are considered.
Given the displacement field~$\displacement(\fx) \in \ffR^\nDims$ and the symmetric and positive definite fourth-order stiffness~$\stiffTensor \in \spdSpace{\ffR^{\nDims \times \nDims \times \nDims \times \nDims}}$, the strain~$\strain(\fx) \in \symSpace{ \ffR^{\nDims \times \nDims} }$ and the stress~$\stress(\fx) \in \symSpace{ \ffR^{\nDims \times \nDims} }$ are
\begin{align}
    \strain(x) &= \symOp{ \displacement(\fx) \otimes \nabla } = \frac{1}{2} \left( \displacement(\fx) \otimes \nabla + \nabla \otimes \displacement(\fx) \right), &
    \stress(\fx) &= \stiffTensor(\fx) \strain(\fx) \,.
\end{align}
The quasi-static balance of linear momentum then reads
\begin{align} \label{eq:PDE-mechanical}
    \stress(\fx) \cdot \nabla & = \zero \, .
\end{align}

Since~$\strain(\fx)$ and~$\stress(\fx)$ are symmetric second-order tensors, these can be expressed using vectors  $\strainMandel, \stressMandel \in \ffR^{\nDims (\nDims + 1)/2}$ based on the Mandel notation---see \cref{sec:appendix-notation}.
For two-dimensional mechanical problems ($\nDims = 2$), the plane strain assumption is made.
We assume isotropic materials, for which the fourth-order stiffness tensor~$\stiffTensor(\fx)$ is defined as
\begin{align}
    \stiffTensor(\fx) = \stiffTensor_0 \chi_0(\fx) + \stiffTensor_1 \chi_1(\fx) \,, &&
    \stiffTensor_i = 3 \lameLambda_i \, \PIso_1 + 2 \lameMu_i \, \Isym \in \ffR^{\nDims \times \nDims \times \nDims \times \nDims} \,, && i \in \{0, 1\} \,,
\end{align}
with the Lamé coefficients~$\lameLambda_i, \lameMu_i>0$ and the isotropic projectors as given in \cref{sec:appendix-notation}. Again, the assumption of isotropic material behavior for each phase represents no general restriction of the presented algorithms.
As for thermal homogenization problems, we employ a decomposition for the field of the displacement field~$\displacement(\fx)$ into
\begin{equation} \label{eq:mechanical-decomposition}
    \displacement(\fx) = \strainMacro \cdot \fx + \dispFluct(\fx) \,,
\end{equation}
where the field~$\strainMacro \cdot \fx$ is induced by a prescribed macroscopic strain~$\strainMacro$ and~$\dispFluct(\fx)$ is the microscopic displacement fluctuation field satisfying~$\int_{\domain} \dispFluct(\fx) \dInt{\fx} = \sty{0}$.
The following boundary conditions are considered:
\begin{itemize}
    \item periodic boundary conditions, i.e.,
    $\displaystyle        \onEdgeP{ \dispFluct } = \onEdgeM{ \dispFluct } \,$,
    \item Dirichlet boundary conditions or uniform kinematic boundary conditions (UKBC), i.e.,
    $   \onBoundary{ \dispFluct } = \sty{0} \,$,
    \item and mixed boundary conditions as a combination of both, e.g., for some $1 \leq \dimJ \leq \nDims$,
    \begin{align}
        \onEdgeIP{ \dispFluct } = \onEdgeIM{ \dispFluct } \,, &&
        \text{for} && i \in \{1, \dots, \nDims\} \setminus \{ \dimJ \} \,, &&
        \text{and} &&
        \onEdgeJP{ \dispFluct} = \onEdgeJM{ \dispFluct} = \sty{0} \,.
    \end{align}
\end{itemize}
The different boundary conditions of homogenization problems are also discussed in, e.g., \cite{Schroder2014, Mercer2015}.
Mixed boundary conditions could be particularly useful for cold rolling applications \cite{Song2024}.
The linear elasticity homogenization problem with the Lamé coefficients~$\lameLambda(\fx)$,~$\lameMu(\fx)$ can be described in the framework of parametric PDEs via
\begin{align} \label{eq:operator-mechanical}
     \diffOpParam(\bullet) = \div{\left(3 \paramSymbol_1(\fx) \PIso_1 + 2 \paramSymbol_2(\fx) \Isym \right) \symOp{ \bullet \otimes \nabla } } \,, &&
     \param = \left\{ \param_1, \param_2, \param_3 \right\} \,,
\end{align}
\begin{align}
     \param_1 = \lameLambda_0 \chi_0(\fx) + \lameLambda_1 \chi_1(\fx) ,\, &&
     \param_2 = \lameMu_0 \chi_0(\fx) + \lameMu_1 \chi_1(\fx) ,\, &&
     \param_3 = \strainMandelMacro \,.
\end{align}

\section{Numerical results}
\label{sec:results}

In this section, the \UNOCG{} hybrid solver is applied to several large-scale homogenization problems that are summarized in \cref{tab:results-overview}.
These problems are given by the parametric PDEs defined in \cref{sec:problems}.
Each sample involves chosen parameters $\param$ that are discretized on the grid given by the two-dimensional or three-dimensional microstructures (see \cref{fig:microstructures}).
Both the stiffness matrix and the right-hand side vector depend on $\param$.

\begin{table}[ht!]
    \vspace{-2mm}
    \caption{Overview of the homogenization problems studied in the following and the amount of data that is involved for a single sample.}
    \centering \begin{tabular}{cccccc}
        \toprule 
        Problem & Domain & Parameters $\param$ & 
        Discretized parameters &
        Solution~$\solParam$ & DOF for~$\solParam$ \\ 
        \midrule 
        Heat conduction (\cref{ssec:thermal-problems}) & $\domain \subsetneq \ffR^2$ & $\Omega \to \ffR$, $\tempGradMacro \in \ffR^2$ & $\ffR^{400 \times 400}$, $\ffR^2$ & $\Omega \to \ffR$ & $\approx 1.6 \cdot 10^5$ \\
        Heat conduction (\cref{ssec:thermal-problems}) & $\domain \subsetneq \ffR^3$ & $\Omega \to \ffR$, $\tempGradMacro \in \ffR^3$ & $\ffR^{192 \times 192 \times 192}, \ffR^3$ & $\Omega \to \ffR$ & $\approx 7.1 \cdot 10^6$ \\ 
        \midrule 
        Linear elasticity (\cref{ssec:mechanical-problems}) & $\domain \subsetneq \ffR^2$ & $\Omega \to \ffR^2$, $\strainMandelMacro \in \ffR^3$ & $\ffR^{2 \times 400 \times 400}$, $\ffR^3$ 
        & $\Omega \to \ffR^2$ & $\approx 3.2 \cdot 10^5$ \\
        Linear elasticity (\cref{ssec:mechanical-problems}) & $\domain \subsetneq \ffR^3$ & $\Omega \to \ffR^2$, $\strainMandelMacro \in \ffR^6$ & $\ffR^{{\color{revision}3} \times 192 \times 192 \times 192}$, $\ffR^6$ 
        & $\Omega \to \ffR^3$ & $\approx 2.1 \cdot 10^7$ \\
        \bottomrule
    \end{tabular}
    \label{tab:results-overview}
    \vspace{-2mm}
\end{table}

For each problem and for different boundary conditions on~$\solParam$, training and test data sets are generated for different parameters $\param$.
While the \UNO{} preconditioner is learned based on training data, the resulting \UNOCG{} hybrid solver is evaluated only on test data in the following.
More information about the data sets is available in \cref{ssec:hyperparameters}, {\color{revision}and \cref{ssec:runtimes} reports the corresponding wall-clock runtimes for all solvers used in the following numerical experiments.}

\subsection{Linear heat conduction in 2D heterogeneous media}
\label{ssec:results-thermal-2d}

\begin{figure}[ht!]
	\centering
	\includegraphics[width=\textwidth]{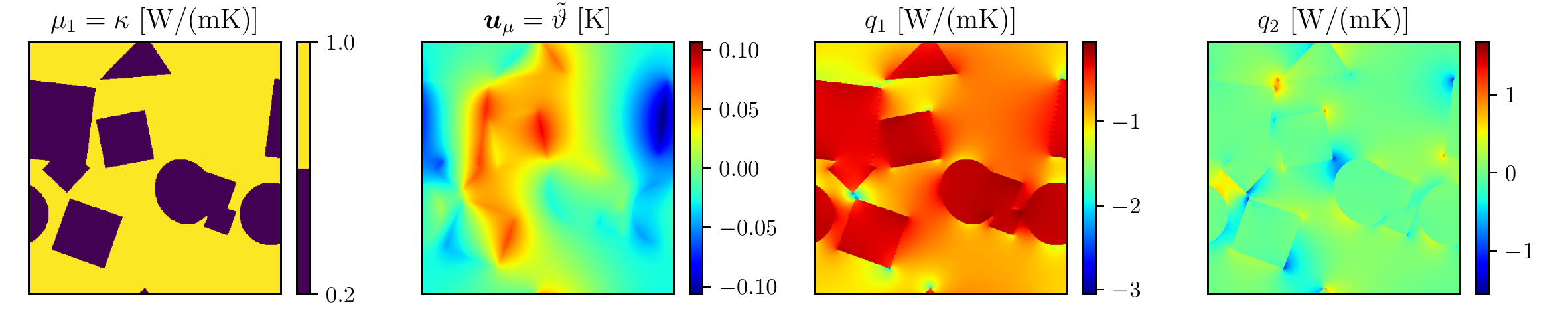}
    \vspace{-1.5\baselineskip}
	\caption{Parameter~$\param_1(\fx)$, solution~$\sol_{\param}(\fx)$ (temperature fluctuation~$\tempFluct(\fx)$), and the components of the heat flux field~$\flux(\fx)$ of the thermal 2D problem with periodic BC for the macroscopic temperature gradient $\tempGradMacro=\T{\begin{bmatrix} 1 & 0 \end{bmatrix}}$.}
	\label{fig:results_thermal_2d_per}
\end{figure}

First, we consider thermal homogenization problems as introduced in \cref{ssec:thermal-problems} on 2D microstructures with a resolution~$400 \times 400$, resulting in~$1.6 \cdot 10^5$ DOF.
Following the dataset \cite{Darus}, the material parameters of the two phases are chosen to be~$\heatCond_0 = 1 \,\fluxUnit$ and~$\heatCond_1 = 0.2 \,\fluxUnit$.
This corresponds to a phase contrast of~$R=\frac{\heatCond_0}{\heatCond_1} = 5$.
The parameter~$\param_1(\fx)$, solution~$\sol_{\param}(\fx)$ (i.e., temperature fluctuation~$\tempFluct(\fx)$), and the components of the heat flux field~$\flux(\fx)$ of the thermal 2D problem with periodic BC are shown in \cref{fig:results_thermal_2d_per} for one exemplary microstructure and a prescribed macroscopic temperature gradient of~$\tempGradMacro=\T{\begin{bmatrix} 1 & 0 \end{bmatrix}}$.

As described in \cref{ssec:fnocg-fasttraining}, the \UNO{} preconditioner is trained on the available data using a second-order optimization scheme, leading to the \UNOCG{} hybrid solver.
For comparison, we train the preconditioner on the same data using the naive training procedure outlined in \cref{ssec:fnocg-naive-training} and denote the resulting hybrid solver as \UNOCGnaive{}.
For both approaches, a separate tuning of the hyperparameter~$\modeParam$ is performed.
Since the improved preconditioner training procedure from \cref{ssec:fnocg-fasttraining} exhibits better convergence properties, it allows the use of higher~$\modeParam$ without risking overfitting and instabilities during training.

\begin{figure}[ht!]
    \centering
    \includegraphics[width=\textwidth]{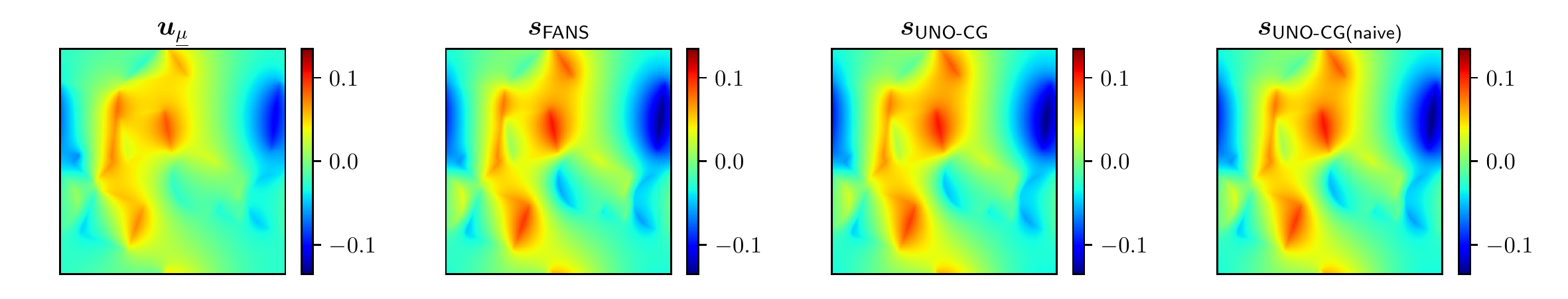}
    \vspace{-1.5\baselineskip}
    \caption{Comparison of the final solution $\solParam$ and the preconditioned residual $\precRes$ in the first {\color{revision}solver iteration} for the different preconditioned \CG{} solvers \FANS{}, \UNOCG{}, and \UNOCGnaive{}.}
    \label{fig:prec_action_thermal_2d_per}
\end{figure}

The \UNOCG{} solver is then tested on unseen parameters $\param\sampleIdx$ from a test data set.
More details about the considered training data and test data are available in \cref{ssec:hyperparameters}.
As a quality assessment of the different preconditioned solvers \FANS{}, \UNOCGnaive{}, and \UNOCG{}, the preconditioned residual in the first iteration $\precRes = \Prec \, \res^{(0)} = \Prec \, \rhsParam$ (since $\sol\cgIdxFirst_{\param} = \sty{0}$) is {\color{revision}qualitatively} compared with the final solution $\solParam$ in \cref{fig:prec_action_thermal_2d_per}.

In the following, the convergence behaviour of the \UNOCG{} hybrid solver is compared with unpreconditioned \CG{}, Jacobi-preconditioned \CG{} ($\JacCG$), and \FANS{}.
The condition numbers and convergence ratios for the exemplary sample shown in \cref{fig:results_thermal_2d_per} are computed and summarized in \cref{tab:thermal-precs}.
Based on the error bound in~\eqref{eq:cg-convergence},~$\cgItersBound$ is available as the estimated maximum number of iterations required until the stopping criterion~$\cgTol = 10^{-6}$ is reached, and~$\cgIters$ denotes the actual number of iterations.
Note that~$\cgItersBound$ is based on the error in the energy norm and is not necessarily a strict bound on the relative residual used to measure~$\cgIters$.
In addition, the effectiveness of the different preconditioners is compared in \cref{tab:thermal-precs}.

\begin{table}[ht!]
    \caption{Properties of solvers for the 2D thermal problem with the exemplary sample from \cref{fig:results_thermal_2d_per}. $\cgRate$ is the convergence rate based on $\cond$ used to estimate the maximum number of iterations $\cgItersBound$ based on~\eqref{eq:cg-convergence}. The effectiveness of the different preconditioners is measured in the reduction in $\cgIters$ over unpreconditioned \CG{}.}
    \centering \begin{tabular}{cccccccc}
        \toprule 
        Solver & $\eigMax \left( \PrecD\,\StiffMat \right)$ & $\eigMin \left( \PrecD\,\StiffMat \right)$ & $\cond \left( \PrecD\,\StiffMat \right)$ & $\cgRate \left( \PrecD\,\StiffMat \right)$ & $\cgItersBound$ & $\cgIters$ & Prec. effectiveness \\ 
        \midrule 
        \CG{} ($\PrecD=\IdentD$) & $3.9998$ & $0.0001$ & $30925.8650$ & $0.9887$ & $1276$ & $\sty{1039}$ & --- \\
        \JacCG & $1.5000$ & $0.0001$ & $23160.4012$ & $0.9870$ & $1104$ & $\sty{895}$ & $\approx 1.2$ \\
        \UNOCGnaive & $6.6966$ & $0.1659$ & $40.3699$ & $0.7280$ & $46$ & $\sty{31}$ & $\approx 33.5$ \\
        \UNOCG & $1.8778$ & $0.2504$ & $7.4983$ & $0.4650$ & $20$ & $\sty{20}$ & $\approx 52.0$ \\
        \FANS{} & $1.6667$ & $0.3333$ & $5.0000$ & $0.3820$ & $15$ & $\sty{14}$ & $\approx 74.2$ \\
        \bottomrule
    \end{tabular}
    \label{tab:thermal-precs}
\end{table}

\begin{figure}[ht!]
    \centering
    \includegraphics[width=\textwidth]{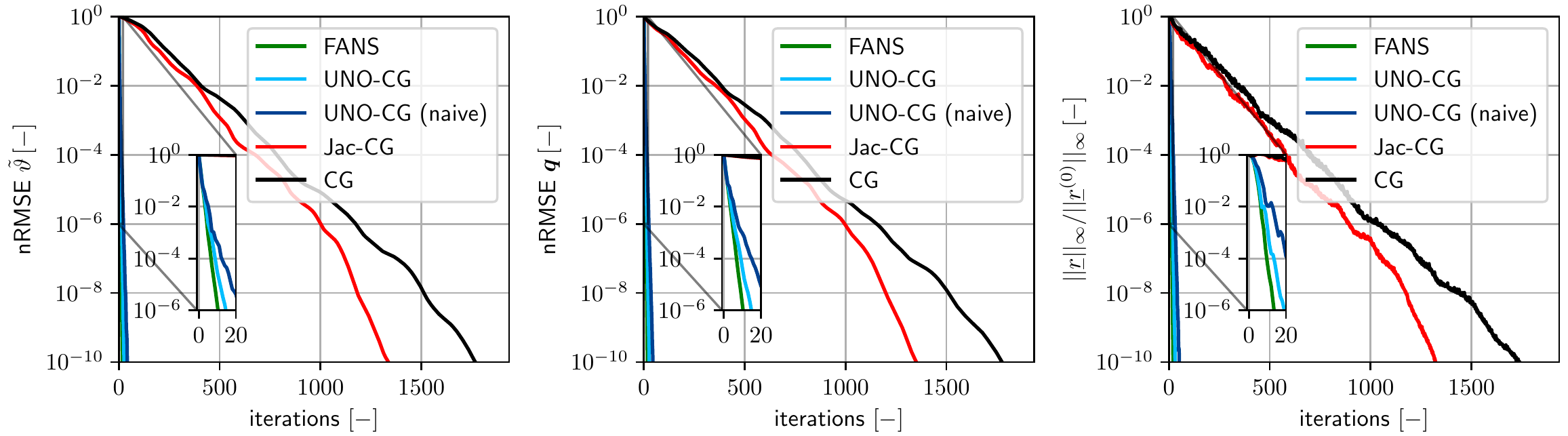}
    \vspace{-1.5\baselineskip}
    \caption{Convergence of \FANS, \UNOCG, \UNOCGnaive, \JacCG, and unpreconditioned \CG{} in the nRMSE on~$\tempFluct$ (left), nRMSE on~$\flux$ (middle), and the relative residual norm~$\normg{\resD}_{\infty} / \normg{\resD^{(0)}}_{\infty}$ (right) for the 2D thermal problem.}
    \label{fig:convergence_thermal_2d_per}
\end{figure}

For this problem, the convergence is measured in the maximum norm of the residual vector~$\normg{\resD}_{\infty}$.
Further, it is measured in the normalized root mean square error (nRMSE) between the iterative approximation~$\bullet_{\rm pred}$ and the final solution~$\bullet_{\rm ref}$ on the temperature fluctuation~$\tempFluct$ and heat flux~$\flux$, which is defined as
\begin{align}
    \text{nRMSE}(\bullet_{\rm pred}, \bullet_{\rm ref}) = \frac{\langle \normg{\bullet_{\rm pred} - \bullet_{\rm ref}}_2 \rangle}{\langle \normg{\bullet_{\rm ref}}_2 \rangle} \,, && \text{where} &&
    \langle \bullet \rangle = \frac{1}{\norm{\domain}} \int_{\domain} \bullet \dInt{\fx} \,.
\end{align}
These error measures are shown in \cref{fig:convergence_thermal_2d_per} for the same exemplary sample.
It can be seen that \UNOCG{} and \UNOCGnaive{} clearly outperform unpreconditioned \CG{} and $\JacCG$, while \UNOCG{} itself converges significantly faster than \UNOCGnaive{} (see also preconditioner effectiveness in \cref{tab:thermal-precs}).
However, it is not possible for \UNOCG{} to outperform the \FANS{} preconditioner that is analytically derived for this specific problem formulation.

\begin{figure}[ht!]
    \centering
    \begin{subfigure}[b]{0.49\textwidth}
        \includegraphics[scale=0.69]{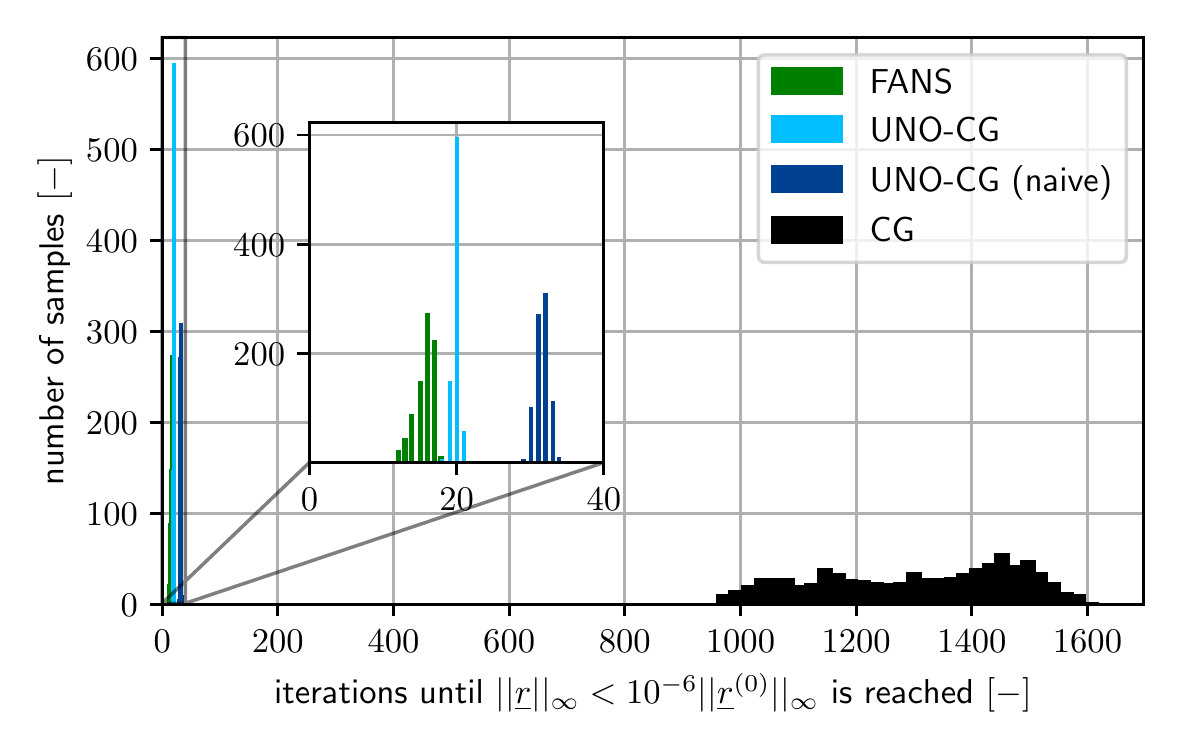}
        \vspace{-1.0\baselineskip}
        \caption{Robustness study on all 800 test samples.}
        \label{fig:convergence_histogram_thermal_2d_per}
    \end{subfigure}
    \hfill
    \begin{subfigure}[b]{0.49\textwidth}
        \includegraphics[scale=0.69]{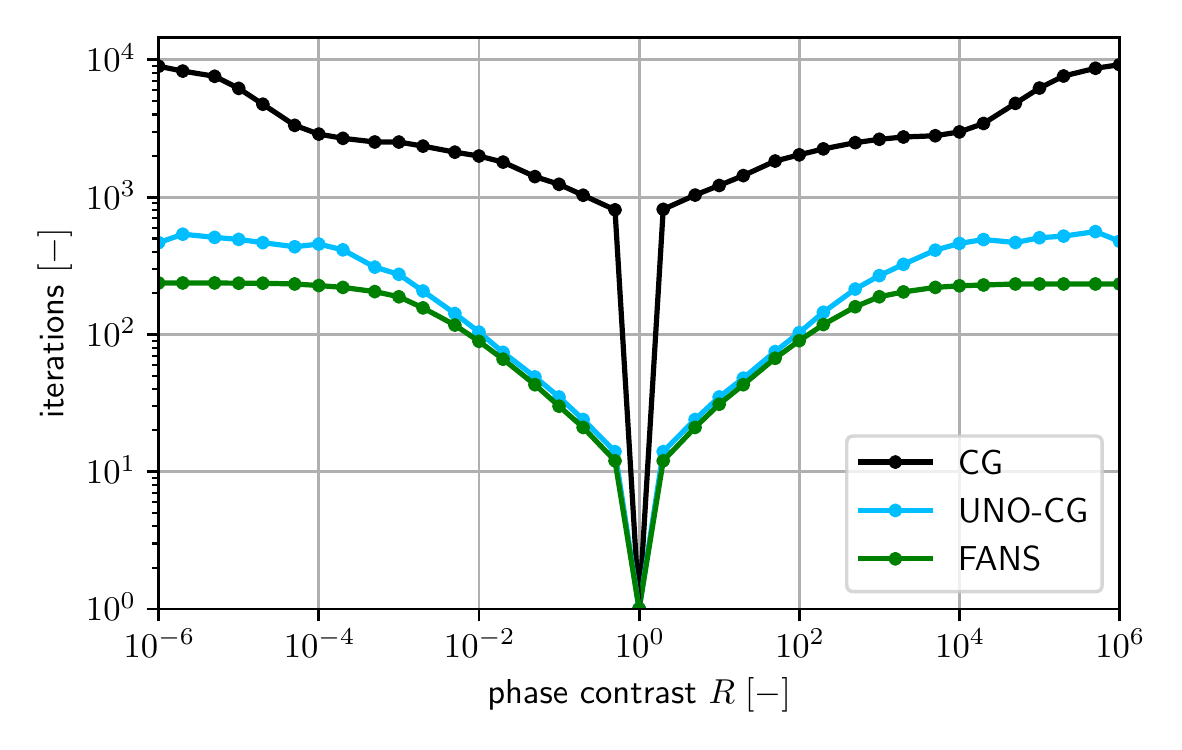}
        \vspace{-1.0\baselineskip}
        \caption{Generalization to varying phase contrasts $R=\kappa_0/\kappa_1$.}
        \label{fig:scaling_thermal_2d_per}
    \end{subfigure}
\vspace{-0.2\baselineskip}
    \caption{Evaluation of robustness and generalization capabilities of \UNOCG{} compared with \FANS{} and unpreconditioned \CG{} for the thermal 2D problem.}
    \label{fig:sweep_workflow_thermal_2d_per}
\end{figure}

For a robustness study on the entire test dataset containing 800 samples (see also \cref{ssec:hyperparameters}), the histogram in \cref{fig:convergence_histogram_thermal_2d_per} shows the number of iterations that is needed per sample to reach a relative tolerance of~$10^{-6}$ for the residual~$\normg{\resD}_{\infty}$.
Iteration counts follow a sharp distribution, confirming the finding for the single sample across the entire test set.
The generalization capabilities of the \UNO{} preconditioner are tested with respect to increasingly high phase contrasts~$R = \frac{\heatCond_0}{\heatCond_1}\gg 1$ or~$R \ll 1$ in \cref{fig:scaling_thermal_2d_per}, while it is only trained using data gained for~$R = 5$.
Here, \UNOCG{} shows a similar scaling with~$R$ as \FANS{}.

To summarize, these first results show that the hybrid solver \UNOCG{} is able to resemble the analytically derived preconditioner \FANS{} closely and almost matches its performance while only being based on data.
In the following, the capabilities of \UNOCG{} are also examined for scenarios where \FANS{} can not be applied.

\subsection{Linear heat conduction in 3D heterogeneous media}
\label{ssec:results-thermal-3d}

\begin{figure}[ht!]
    \centering
    \begin{subfigure}[b]{0.495\textwidth}
        \begin{center}
            \includegraphics[scale=0.15]{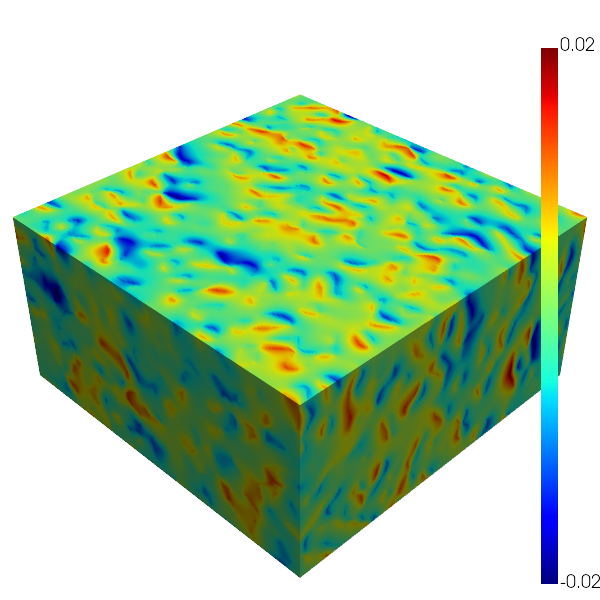}
            \includegraphics[scale=0.15]{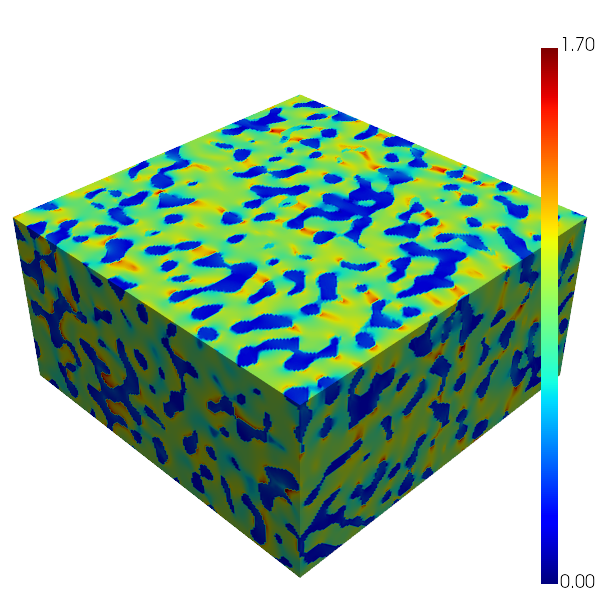}
        \end{center}
        \caption{\color{revision}Fields~$\tempFluct$ (left) and~$\normg{\flux}$ (right) for periodic BC.}
        \label{fig:results_thermal_3d_per}
    \end{subfigure}
    \hfill
    \begin{subfigure}[b]{0.495\textwidth}
        \begin{center}
            \includegraphics[scale=0.15]{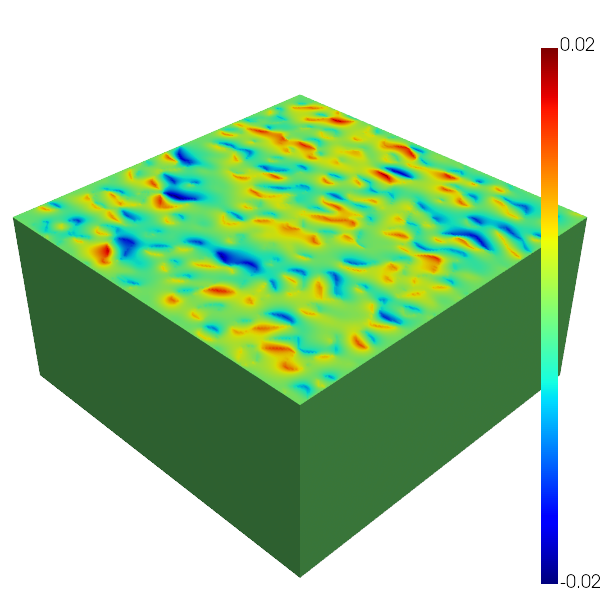}
            \includegraphics[scale=0.15]{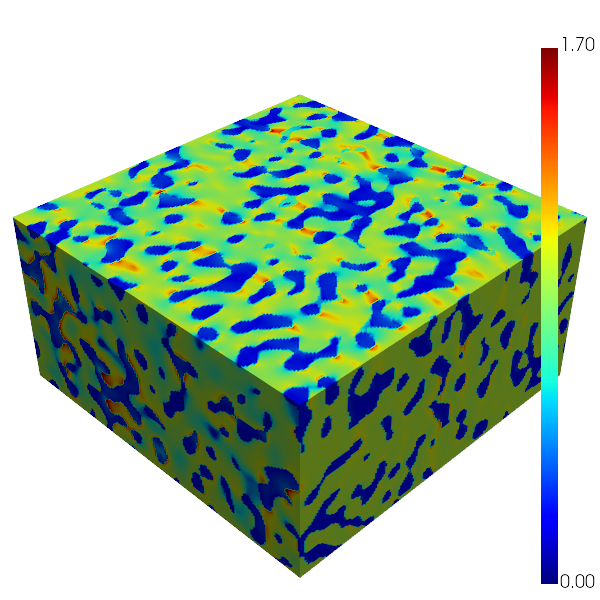}
        \end{center}
        \caption{\color{revision}Fields~$\tempFluct$ (left) and~$\normg{\flux}$ (right) for Dirichlet BC.}
        \label{fig:results_thermal_3d_dir}
    \end{subfigure}
    \caption{Results of the thermal homogenization problem in 3D for different boundary conditions, where each the macroscopic temperature gradient~$\tempGradMacro=\T{\begin{bmatrix} 1 & 0 & 0 \end{bmatrix}}$ is prescribed. For better visualization, the domain is clipped to the lower half.}
    \label{fig:results_thermal_3d}
\end{figure}

The thermal homogenization problem from \cref{ssec:mechanical-problems} is now solved on the 3D microstructures with a resolution $192^3$ and $\approx 7.1 \cdot 10^6$ DOF.
Again, the phase contrast is $R=\frac{\heatCond_0}{\heatCond_1}=5$.
To explore the capabilities of \UNOCG{} for other boundary conditions where \FANS{} is not applicable, we compare periodic BC and Dirichlet BC on the temperature fluctuations~$\tempFluct$, respectively.
For these boundary conditions, the temperature fluctuation field $\tempFluct(\fx)$, and the norm of the heat flux field $\normg{\flux(\fx)}$ are shown in \cref{fig:results_thermal_3d} for a prescribed macroscopic temperature gradient of~$\tempGradMacro=\T{\begin{bmatrix} 1 & 0 & 0 \end{bmatrix}}$.

\begin{figure}[ht!]
    \centering
    \begin{subfigure}[b]{0.495\textwidth}
        \begin{center}
            \includegraphics[scale=0.59]{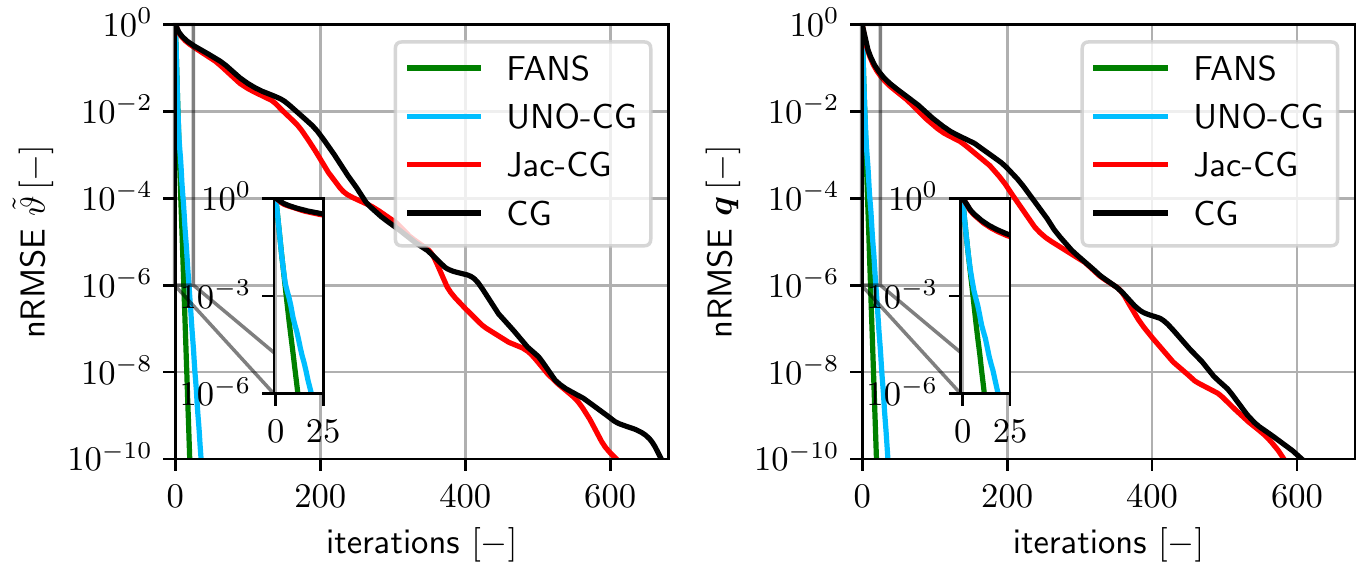}
			\vspace{-1\baselineskip}
        \end{center}
        \caption{Convergence of $\tempFluct$ (left) and $\flux$ (right) for periodic BC.}
        \label{fig:convergence_thermal_3d_per}
    \end{subfigure}
    \hfill
    \begin{subfigure}[b]{0.495\textwidth}
        \begin{center}
            \includegraphics[scale=0.59]{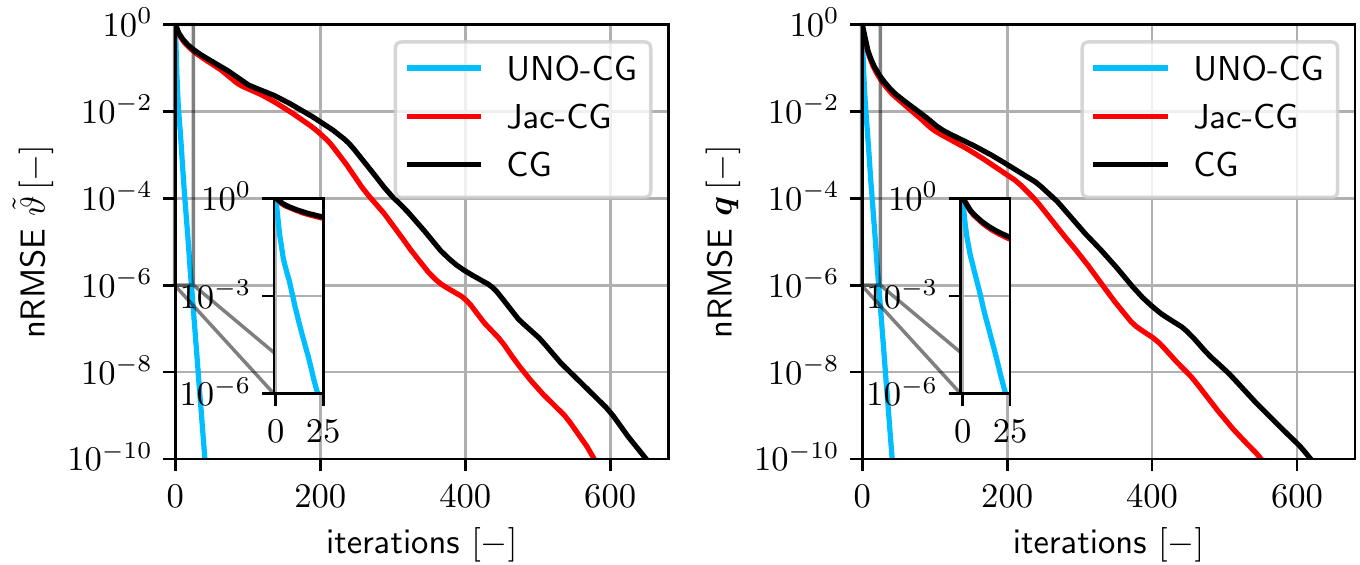}
			\vspace{-1\baselineskip}
        \end{center}
        \caption{Convergence of $\tempFluct$ (left) and $\flux$ (right) for Dirichlet BC.}
        \label{fig:convergence_thermal_3d_dir}
    \end{subfigure}
    \vspace{-0.5\baselineskip}
    \caption{Convergence of \FANS, \UNOCG, Jacobi-preconditioned \CG{} (\JacCG), and unpreconditioned \CG{} for the thermal homogenization problem in 3D with different boundary conditions as shown in \cref{fig:results_thermal_3d}.}
    \label{fig:convergence_thermal_3d}
\end{figure}

Again, the \UNO{} preconditioner is exclusively trained on discrete solutions of the parametric PDE.
It was not possible to consider \UNOCGnaive{} for the 3D problems due to the extensive memory requirements and computational costs for backpropagation, which are induced by the massive size of the problem (millions of double-precision values as inputs per sample).
In the case of Dirichlet BC, the sine transform is used as a unitary transform $\trafo = \trafoDST$ in \UNOCG{}.
The convergence behavior of the different preconditioned \CG{} solvers in various metrics is shown in \cref{fig:convergence_thermal_3d}.
\UNOCG{} can closely resemble the performance of \FANS{} but is not able to outperform it (\cref{fig:convergence_thermal_3d_per}). However, in contrast to \FANS{}, \UNOCG{} can also be applied for Dirichlet BC (\cref{fig:convergence_thermal_3d_dir}), where it achieves a reduction of iterations over unpreconditioned \CG{} that is similar to periodic BC.

\subsection{Linear elasticity in 2D heterogeneous media}
\label{ssec:results-mechanical-2d}

\begin{figure}[ht!]
    \centering
    \begin{subfigure}[b]{0.32\textwidth}
        \begin{center}
            \includegraphics[width=1.0\textwidth]{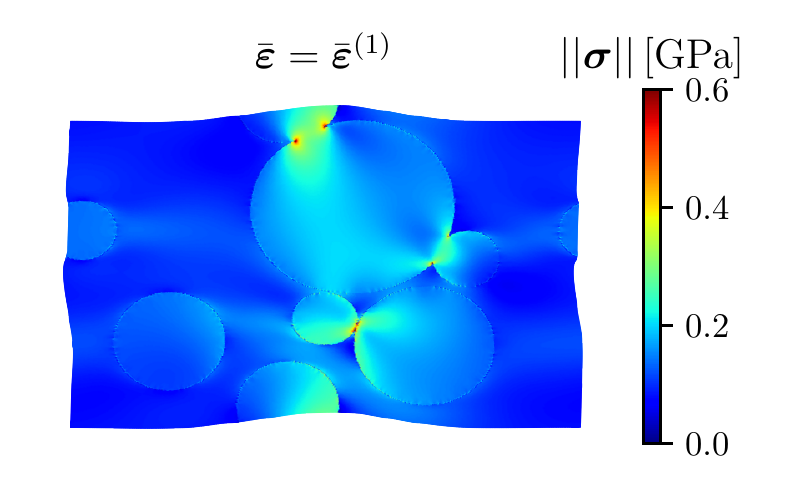}
        \end{center}
        \vspace{-1.7\baselineskip}
        \caption{Periodic BC.}
        \label{fig:results_mechanical_2d_per}
    \end{subfigure}
    \hfill
    \begin{subfigure}[b]{0.32\textwidth}
        \begin{center}
            \includegraphics[width=1.0\textwidth]{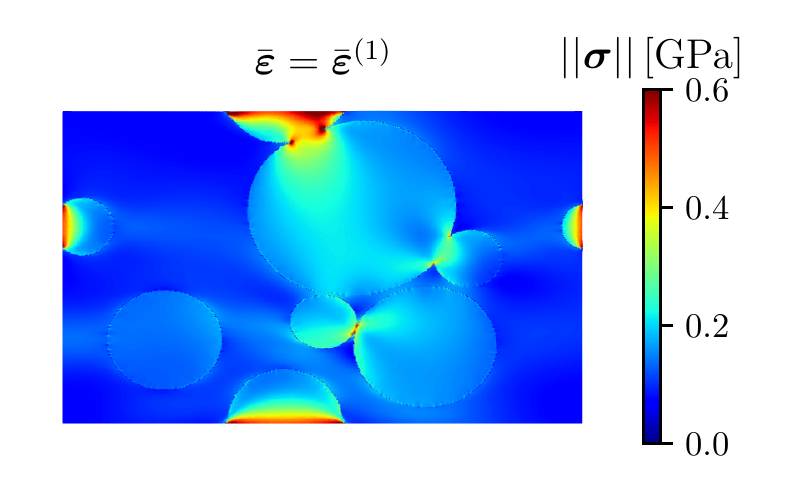}
        \end{center}
        \vspace{-1.7\baselineskip}
        \caption{Dirichlet BC.}
        \label{fig:results_mechanical_2d_dir}
    \end{subfigure}
    \begin{subfigure}[b]{0.32\textwidth}
        \begin{center}
            \includegraphics[width=1.0\textwidth]{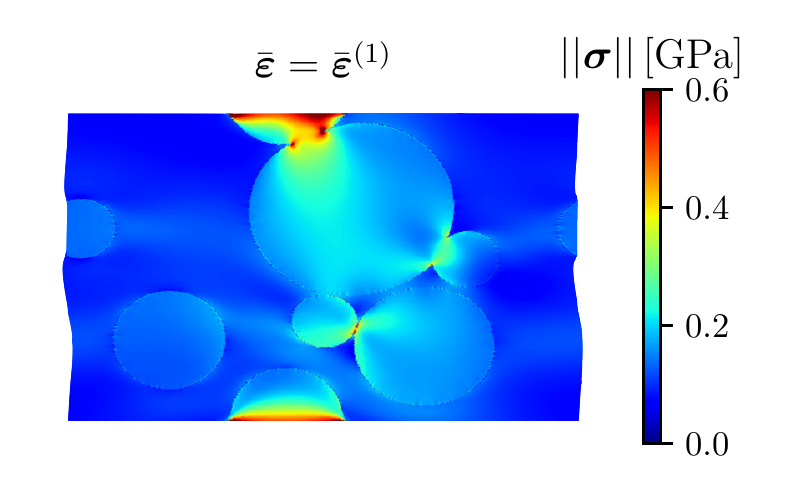}
        \end{center}
        \vspace{-1.7\baselineskip}
        \caption{Mixed BC.}
        \label{fig:results_mechanical_2d_mixed}
    \end{subfigure}
    \centering
    \vspace{-0.2\baselineskip}
    \caption{Stress norm~$\normg{\stress}$ on the deformed RVEs for the mechanical 2D problem with different BC for the macroscopic strain $\strainMacro = \strainMacro^{(1)} = \diag{\begin{bmatrix} 0.05 & {-0.05} \end{bmatrix}}$. The deformations are scaled for better visualization.}
    \label{fig:results_mechanical_2d}
\end{figure}

In the following, the mechanical homogenization as introduced in \cref{ssec:mechanical-problems} is considered.
Here, $\nComp=2$ nodal DOF are required, which results in $\approx 3.2 \cdot 10^5$ DOF for the considered 2D microstructures.
The material parameters of the phases are chosen as~$E_0 = 1 \, \stressUnit$,~$\nu_0 = 0 \, \noUnit$,~$E_1 = 10 \, \stressUnit$, and~$\nu_1 = 0.3 \, \noUnit$.
As boundary conditions on the displacement fluctuation field~$\dispFluct$, periodic BC, Dirichlet BC, and mixed BC (i.e., Dirichlet BC on $\Gamma^\pm_2$, periodic BC on $\Gamma^\pm_1$) are considered.
The Frobenius norm of the stress field~$\normg{\stress}$ is visualized on the deformed RVEs in \cref{fig:results_mechanical_2d} for the different boundary conditions.

\begin{figure}[ht!]
    \centering
    \begin{subfigure}[b]{0.32\textwidth}
        \begin{center}
            \includegraphics[scale=0.69]{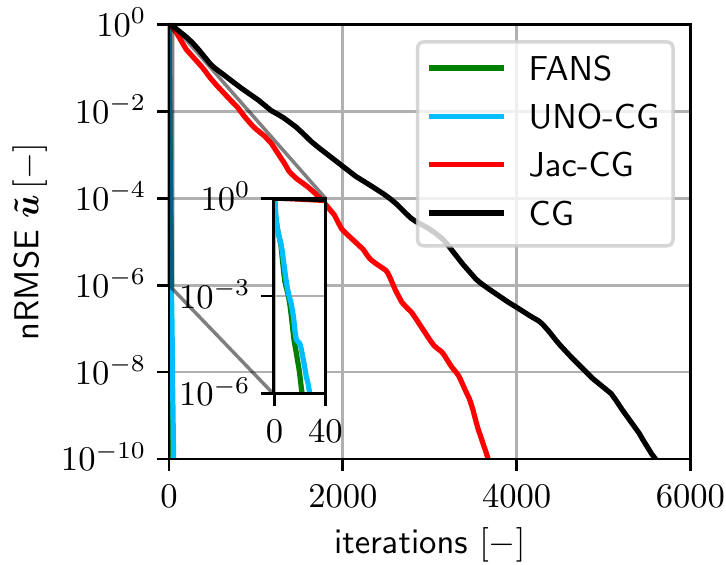}
        \end{center}
        \vspace{-1.0\baselineskip}
        \caption{Periodic BC.}
        \label{fig:convergence_mechanical_2d_per}
    \end{subfigure}
    \hfill
    \begin{subfigure}[b]{0.32\textwidth}
        \begin{center}
            \includegraphics[scale=0.69]{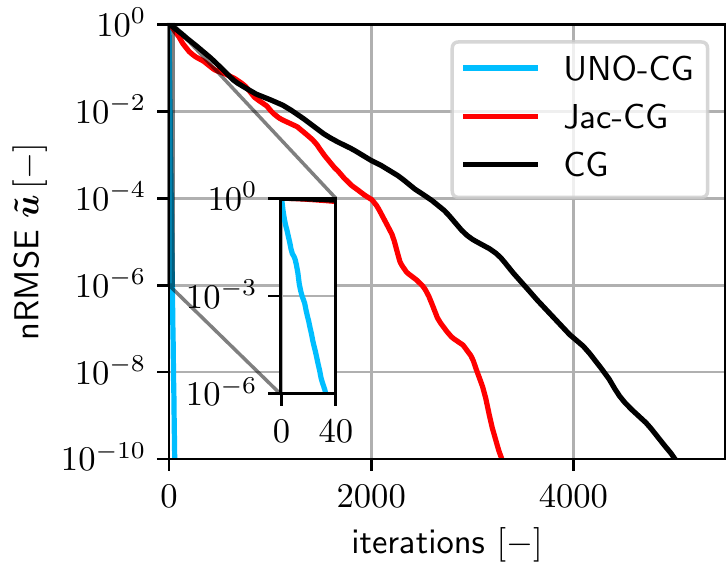}
        \end{center}
        \vspace{-1.0\baselineskip}
        \caption{Dirichlet BC.}
        \label{fig:convergence_mechanical_2d_2d_dir}
    \end{subfigure}
    \begin{subfigure}[b]{0.32\textwidth}
        \begin{center}
            \includegraphics[scale=0.69]{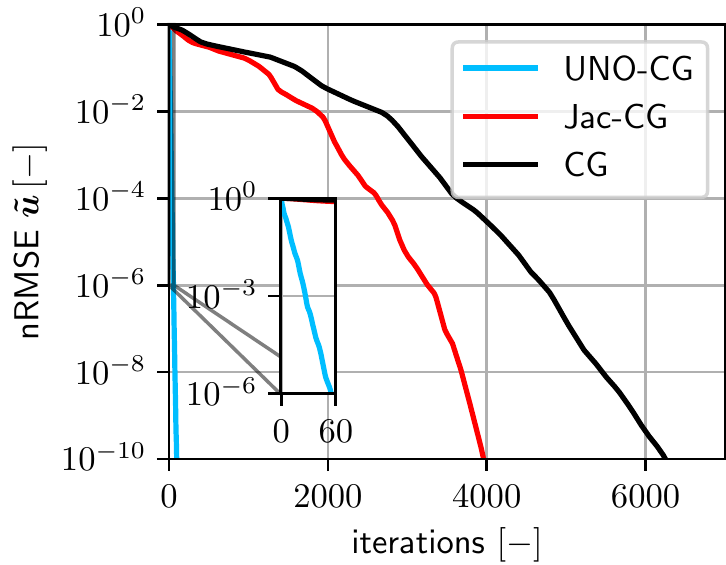}
        \end{center}
        \vspace{-1.0\baselineskip}
        \caption{Mixed BC.}
        \label{fig:convergence_mechanical_2d_mixed}
    \end{subfigure}
    \centering
    \vspace{-0.2\baselineskip}
    \caption{Convergence of \FANS{}, \UNOCG{}, and unpreconditioned \CG{} in the solution~$\dispFluct$ for the mechanical 2D problem with periodic BC (left), Dirichlet BC (middle), and mixed BC (right). \FANS{} is only applicable for periodic BC. Convergence in $\stress$ and $\resD$ is qualitatively similar.}
    \label{fig:convergence_mechanical_2d}
\end{figure}

For periodic BC, the convergence behavior of \UNOCG{} is compared with \FANS{}, \JacCG{}, and unpreconditioned \CG{} in \cref{fig:convergence_mechanical_2d}.
For Dirichlet and mixed BC, it is not possible to apply \FANS{}.
However, \UNOCG{} is able to achieve a convergence similar to the case of periodic BC, as it can also be observed in \cref{fig:convergence_mechanical_2d}.
The preconditioner effectiveness of \UNOCG{}, i.e., the reduction in $\cgIters$, for this problem is $\approx89$ for periodic BC, $\approx77$ for Dirichlet BC, and $\approx69$ for mixed BC (measured as in \cref{tab:thermal-precs}).

\begin{figure}[ht!]
    \centering
    \begin{subfigure}[b]{0.24\textwidth}
        \begin{center}
            \includegraphics[scale=0.55]{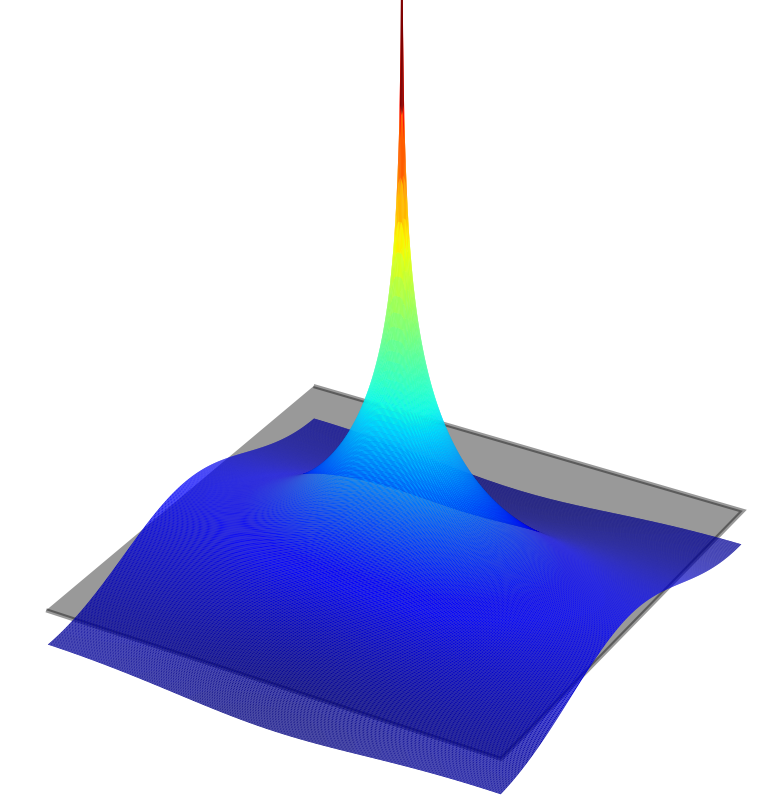}
        \end{center}
        \vspace{-1.0\baselineskip}
        \caption{\FANS{} for periodic BC.}
        \label{fig:fundamental_solution_mechanical_2d_fans}
    \end{subfigure}
    \hfill
    \begin{subfigure}[b]{0.24\textwidth}
        \begin{center}
            \includegraphics[scale=0.55]{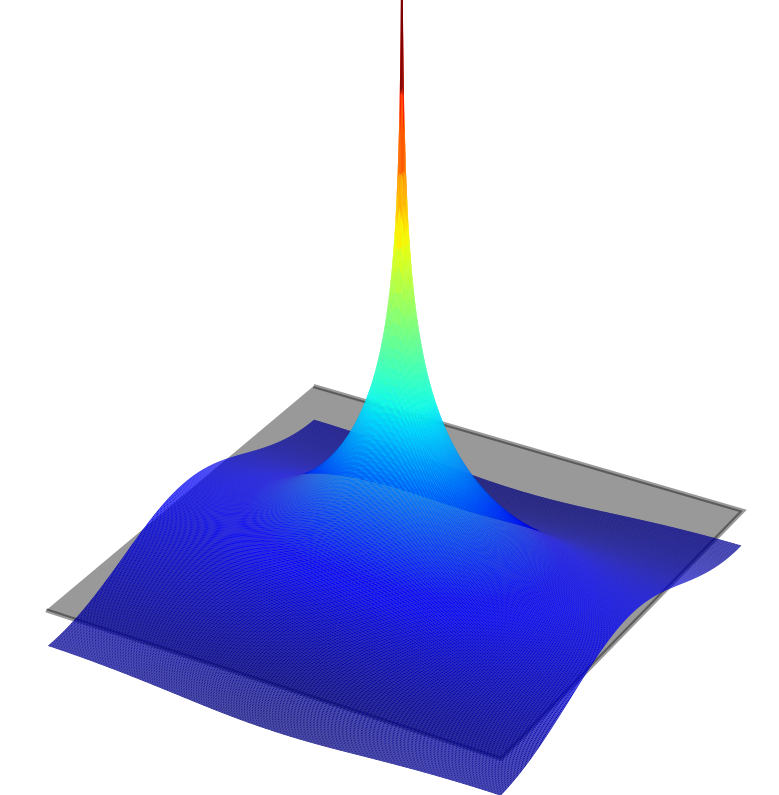}
        \end{center}
        \vspace{-1.0\baselineskip}
        \caption{\UNOCG{} for periodic BC.}
        \label{fig:fundamental_solution_mechanical_2d_per}
    \end{subfigure}
    \hfill
    \begin{subfigure}[b]{0.24\textwidth}
        \begin{center}
            \includegraphics[scale=0.55]{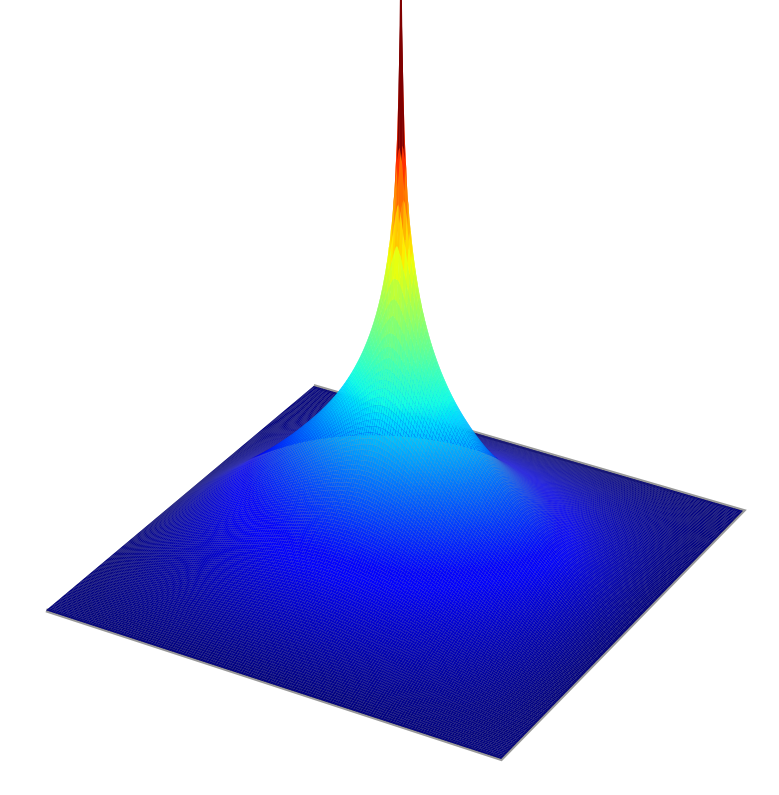}
        \end{center}
        \vspace{-1.0\baselineskip}
        \caption{\UNOCG{} for Dirichlet BC.}
        \label{fig:fundamental_solution_mechanical_2d_dir}
    \end{subfigure}
    \begin{subfigure}[b]{0.24\textwidth}
        \begin{center}
            \includegraphics[scale=0.55]{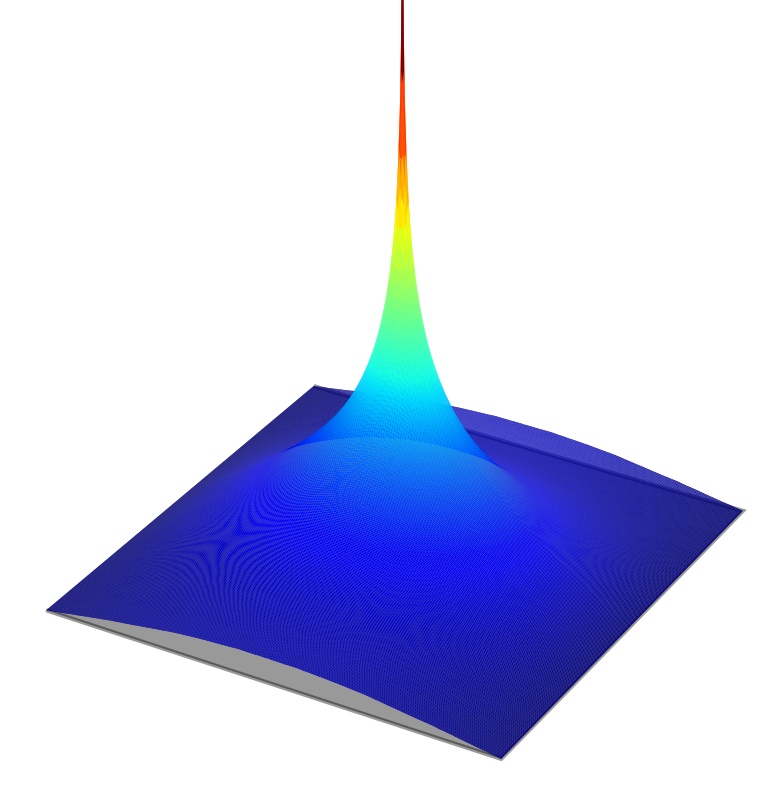}
        \end{center}
        \vspace{-1.0\baselineskip}
        \caption{\UNOCG{} for mixed BC.}
        \label{fig:fundamental_solution_mechanical_2d_mixed}
    \end{subfigure}
    \vspace{-0.2\baselineskip}
    \caption{Fundamental solution~$\FundSol_{22}$ (representing corrections of the $\mathrm{y}$ displacement due to a unit residual force in the $\mathrm{y}$ direction) of \FANS{} together with its machine-learned analogies extracted from \UNOCG{} for the mechanical 2D problem with different BC.
    The \UNO{} preconditioner is trained for each of the different BC separately.}
    \label{fig:fundamental_solution_mechanical_2d}
\end{figure}

Similar to \FANS{}, a matrix-valued fundamental solution~$\FundSolMat$ (i.e., a numerical Green's function) can be extracted from the \UNO{} preconditioner by applying~$\Prec_{\weights}$ on a Dirac delta~$\diracDelta_{\fx_0}$.
If the activation is centered in the origin, i.e.~$\fx_0 = \sty{0}$, the extracted fundamental solution from the machine-learned \UNO{} preconditioner matches the fundamental solution of \FANS{} for periodic BC---as it can be seen in \cref{fig:fundamental_solution_mechanical_2d_fans,fig:fundamental_solution_mechanical_2d_per} for the component~$\FundSol_{22}$.
Also for Dirichlet BC and mixed BC, fundamental solutions can be obtained from the learned \UNO{} preconditioner that satisfy the corresponding BC and are visualized in \cref{fig:fundamental_solution_mechanical_2d_dir,fig:fundamental_solution_mechanical_2d_mixed}.

\subsection{Linear elasticity in 3D heterogeneous media}
\label{ssec:results-mechanical-3d}

\begin{figure}[ht!]
    \centering
    \begin{subfigure}[b]{0.32\textwidth}
        \begin{center}
            \includegraphics[scale=0.2]{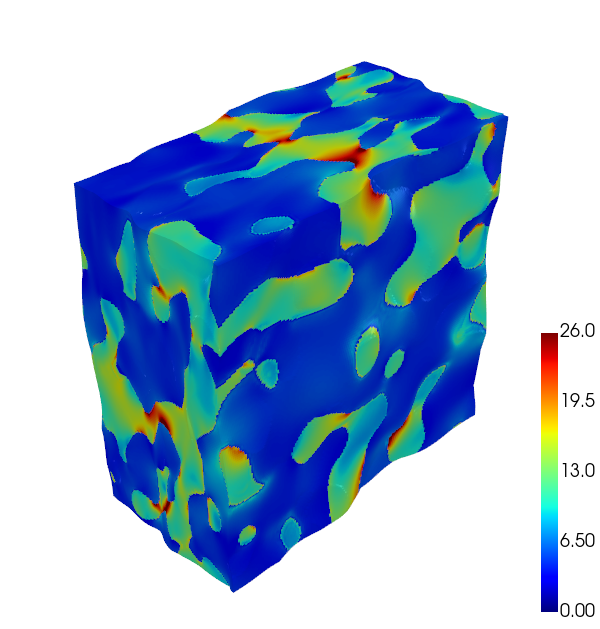}
        \end{center}
        \caption{Periodic BC.}
        \label{fig:results_mechanical_3d_per}
    \end{subfigure}
    \hfill
    \begin{subfigure}[b]{0.32\textwidth}
        \begin{center}
            \includegraphics[scale=0.2]{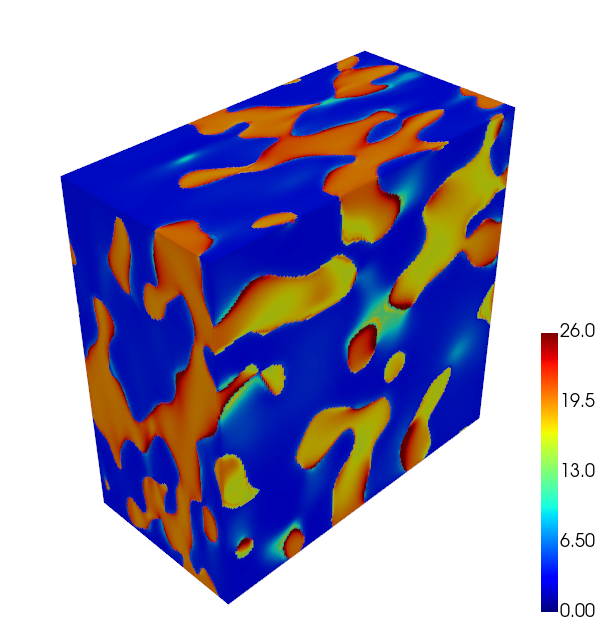}
        \end{center}
        \caption{Dirichlet BC.}
        \label{fig:results_mechanical_3d_dir}
    \end{subfigure}
    \begin{subfigure}[b]{0.32\textwidth}
        \begin{center}
            \includegraphics[scale=0.2]{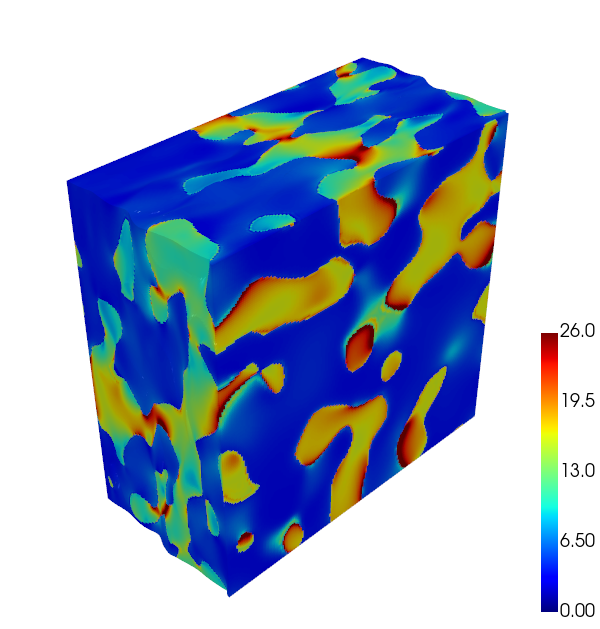}
        \end{center}
        \caption{Mixed BC}
        \label{fig:results_mechanical_3d_mixed}
    \end{subfigure}
    \centering
    \caption{Norm of the stress field~$\normg{\stress}$ $\mathrm{[GPa]}$ on the deformed RVEs for the mechanical 3D problem with periodic BC (left), Dirichlet BC (middle), and mixed BC (right). Deformations are scaled.}
    \label{fig:results_mechanical_3d}
\end{figure}

Finally, \UNOCG{} is applied to the mechanical homogenization problem as introduced in \cref{ssec:mechanical-problems} on 3D microstructures with a resolution of~$192^3$, which leads to a problem with~$\approx 2.1 \cdot 10^7$ DOF.
The material parameters of the phases are in the following chosen to be~$E_0 = 75 \, \stressUnit$,~$\nu_0 = 0.3 \, \noUnit$,~$E_1 = 400 \, \stressUnit$, and~$\nu_1 = 0.2 \, \noUnit$.
Similar to the 2D mechanical problems previously considered, periodic BC, Dirichlet BC, and mixed BC (here, Dirichlet BC on~$\Gamma^\pm_2$, periodic BC on~$\Gamma^\pm_1 \cup \Gamma^\pm_3$) are imposed on the displacement fluctuation fields~$\dispFluct$, respectively.
The norm of the stress field~$\normg{\stress}$ is visualized in \cref{fig:results_mechanical_3d} on the RVE that is deformed by~$\displacement$.
In this example, the macroscopic strain is chosen to be~$\strainMacro = \diag{\begin{bmatrix} 0.025 & -0.05 & 0.025 \end{bmatrix}}$.

\begin{figure*}[ht!]
    \centering
    \begin{subfigure}[b]{0.32\textwidth}
        \begin{center}
            \includegraphics[scale=0.69]{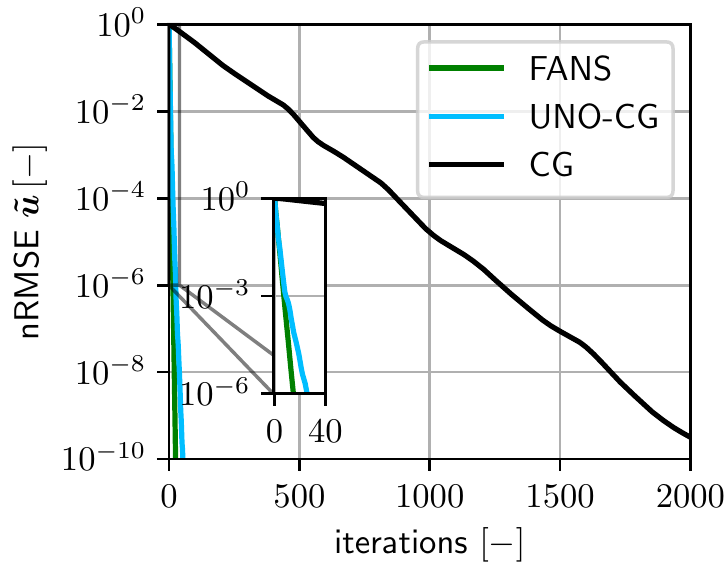}
        \end{center}
        \vspace{-1.0\baselineskip}
        \caption{{\color{revision}Periodic BC.}}
        \label{fig:convergence_mechanical_3d_per}
    \end{subfigure}
    \hfill
    \begin{subfigure}[b]{0.32\textwidth}
        \begin{center}
            \includegraphics[scale=0.69]{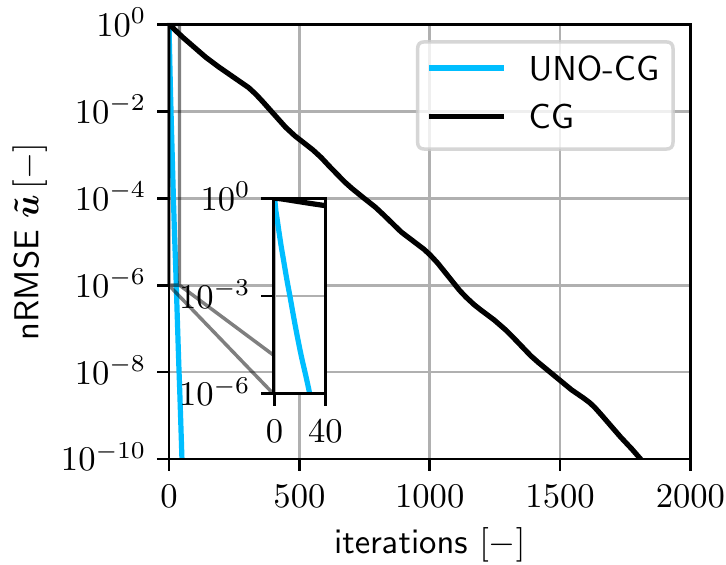}
        \end{center}
        \vspace{-1.0\baselineskip}
        \caption{{\color{revision}Dirichlet BC.}}
        \label{fig:convergence_mechanical_3d_dir}
    \end{subfigure}
    \begin{subfigure}[b]{0.32\textwidth}
        \begin{center}
            \includegraphics[scale=0.69]{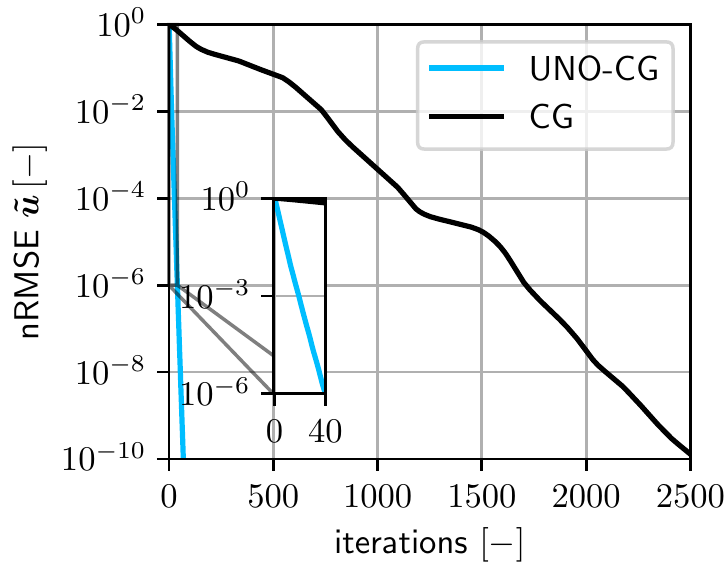}
        \end{center}
        \vspace{-1.0\baselineskip}
        \caption{{\color{revision}Mixed BC.}}
        \label{fig:convergence_mechanical_3d_mixed}
    \end{subfigure}
    \centering
    \caption{Convergence of \FANS{}, \UNOCG{}, and unpreconditioned \CG{} in the solution~$\dispFluct$ for {\color{revision}different boundary conditions}. \FANS{} is only applicable for periodic BC.}
    \label{fig:convergence_mechanical_3d}
\end{figure*}

For the considered boundary conditions, the convergence behavior in the displacement fluctuation fields~$\dispFluct$ is shown in \cref{fig:convergence_mechanical_3d}.
Similar reductions of the number of iterations for \FANS{} and \UNOCG{} over unpreconditioned \CG{} can be observed as for the two-dimensional mechanical problem from \cref{ssec:results-mechanical-2d}.

\section{Résumé}
\label{sec:resume}

\subsection{Methodological summary}
\label{ssec:resume-methods}

With \UNOCG{}, we propose a novel hybrid solver that accelerates conjugate gradient (\CG{}) solvers using our machine-learned \UNO{} preconditioner.
The solver is formulated in the general framework of parametric PDEs. In the present study, it is successfully applied to thermal and mechanical homogenization problems, but it can readily be applied to many other elliptic problems, e.g., to solve steady-state diffusion or permeability problems {\color{revision}based on Fick's or Darcy's law}.
While many previous approaches for machine-learned preconditioners\cite{Xiang2022, Zhang2024} are limited to use with \GMRES{} or similar solvers requiring additional memory and compute time, \UNOCG{} can unleash the well-known benefits of the \CG{} method commonly used for such problems.
\emph{UNOs} feature a general unitary transform instead of the Fourier transform used in Fourier Neural Operators (FNOs \cite{Li_fno_2021}). By limiting attention to a single linear Fourier layer with several constraints, and by replacing the Fourier transform by an arbitrary, more general unitary transform,  a \emph{Unitary Neural Operator} (\UNO{}) is obtained that, by construction, denotes a linear, symmetric, and positive definite operator. The designated layout of UNO renders it admissible as a preconditioner for the \CG{} scheme, and convergence of the gained hybrid \UNOCG{} solver can be guaranteed.
Interestingly, through a suitable choice of the unitary transform, it is easily possible to train \UNOCG{} for various, practically relevant boundary conditions beyond the common periodic BC, where classical solvers based on the FFT, such as $\FANS{}$, are not applicable straightforwardly.

Due to its architecture, \UNOCG{} can also be interpreted as a machine-learned extension of the existing Fourier-Accelerated Nodal Solvers (\FANS{} \cite{Leuschner2018,Keshav2022}) for homogenization problems. 
This gives rise to a physical interpretation of this machine-learned preconditioner, e.g., in terms of fundamental solutions, see, e.g., \cref{fig:fundamental_solution_mechanical_2d}. Hence, a direct comparison with existing and well-studied solvers based on explicit numerical Green's functions is possible, whereas many approaches in the machine learning literature are based on pure black-box models.
In addition to the naive training methods for machine-learned preconditioners from the literature, usually based on computationally involved backpropagation (herein used for \UNOCGnaive{}), we present a novel algorithm for highly memory-efficient and fast training of the \UNO{} preconditioner that accounts for its special structure.
After a single preprocessing step of the training data with a computational cost in $\orderOf{\nSamples n \log n}$, every optimization step in the training can be performed in $\orderOf{n}$.

\subsection{Result summary}
\label{ssec:resume-results}

In this work, \UNOCG{} is applied to linear thermal and mechanical homogenization problems, which are solved on two-dimensional and three-dimensional microstructures. This leads to linear systems of equations with up to 21 million DOF, for which direct solvers or standard iterative solvers take a prohibitively long time.
Even the use of simple algebraic preconditioners for \CG{}, such as the Jacobi preconditioner, cannot sufficiently reduce the time to solution.
Further, the use of most algebraic preconditioners in a many-query scenario is hampered by the fact that they explicitly depend on the stiffness matrix and have to be reassembled for each new parameter.
By using the \UNOCG{} hybrid solver, it is possible to significantly reduce the number of iterations (measured by the preconditioner effectiveness) compared to unpreconditioned \CG{} for all examined types of boundary conditions (periodic, Dirichlet, mixed). The number of iterations to convergence for periodic BC comes close to that of the analytically derived preconditioner \FANS{}. Various analyses show that \UNOCG{} for periodic BC emulates \FANS{} with high accuracy, although the machine-learned preconditioner is only trained on the basis of training data without physical insights such as, e.g., information on the local constitutive model.
A robustness study over the entire test dataset underlines the reliability of \UNOCG{} independent of the phase contrast~$R$, which only mildly affects the iteration count.
The proposed novel training procedure for \UNOCG{} not only significantly reduces the cost of preconditioner training over \UNOCGnaive{}, but the learned preconditioner also shows a much faster convergence behavior in the hybrid solver compared to naively trained \UNOCG{}.
Through this seminal improvement of the training, \UNOCG{} training can be performed for large-scale 3D problems without the need for excessive compute power and memory. For this reason, the use of \UNOCGnaive{} for the 3D problems in the current study is infeasible. However, the numerical results indicate that the new training method is not only faster and less memory-demanding but also less prone to overfitting. It also leads to a much improved preconditioner effectiveness of \UNOCG{} vs. \UNOCGnaive{}, see \cref{tab:thermal-precs}.

\subsection{Future perspective}
\label{ssec:resume-outlook}

The current study is confined to linear, parametric PDEs.
In future work, it would be interesting to investigate \UNOCG{} also for nonlinear parametric PDEs and especially homogenization problems with nonlinear material behavior.
For \FANS{}, it has already been shown in \cite{Leuschner2018,Keshav2022} that a preconditioner designed for a linearized problem can also perform well for nonlinear problems involving plasticity or hyperelasticity at finite strains.
We think that \UNOCG{} might even outperform the existing solvers without needing to invest in the selection of hyperparameters (e.g., the reference material must be chosen suitably in \FANS{}).
Regardless of this, \UNOCG{} could also provide acceleration for other boundary conditions where \FANS{} is not applicable, also in the nonlinear regime.
Furthermore, the \UNOCG{} hybrid solver could be applied to multiphysical coupled homogenization problems such as electro-magneto-mechanical problems in a straightforward way.
Since our PyTorch-based implementation of \UNOCG{} is fully differentiable, it is well suited for solving inverse problems, which could also be explored in the future.
While we have restricted ourselves to regular grids in this work, it might be worthwhile to investigate whether an extension of \UNOCG{} could also operate on unstructured meshes similar to approaches in \cite{Li2022,Bagchi1999}.
Furthermore, the generalization of the preconditioner to depend on the parameters, i.e., on the microstructure, the material properties of the individual phases therein, and on the loading condition, is the subject of current investigations, but well beyond the scope of the current study.

\bmsection*{Author contributions}

\textbf{Julius Herb:} writing -- original draft; investigation; methodology; formal analysis; software; visualization; data curation; conceptualization (equal); validation (equal).
\textbf{Felix Fritzen:} writing -- review \& editing; supervision; conceptualization (equal); validation (equal); funding acquisition; project administration; resources.

\bmsection*{Acknowledgments}
This research was partially funded by the Ministry of Science, Research, and the Arts (MWK) Baden-Württemberg, Germany, within the Artificial Intelligence Software Academy (AISA).
Contributions by Felix Fritzen are funded by Deutsche Forschungsgemeinschaft (DFG, German Research Foundation) within the Heisenberg program DFG-FR2702/10 project 517847245. Funded by Deutsche Forschungsgemeinschaft (DFG, German Research Foundation) under Germany's Excellence Strategy - EXC 2075 – 390740016. We acknowledge the support by the Stuttgart Center for Simulation Science (SimTech).
We thank Prof. Mathias Niepert, Institute for Artificial Intelligence, University of Stuttgart, for inspiring discussions and feedback on the manuscript draft.

\bmsection*{Financial disclosure}

None reported.

\bmsection*{Conflict of interest}

The authors declare no potential conflict of interest.

\bmsection*{Data Availability Statement}

A software package for \UNOCG{}, including a GPU-accelerated implementation of the hybrid solver in PyTorch \cite{PyTorch,PyTorch2}, an implementation for PETSc \cite{Petsc1997}, and the training procedures proposed in this article, is available on GitHub \cite{UNOCGgithub2025}.
All results and figures in this article can be reproduced using this software package\cite{UNOCGgithub2025}.
All considered data sets are available on DaRUS \cite{UNOCGdarus2025}.

\bibliography{literature}

@article{Prifling2021,
	title        = {{Large-Scale Statistical Learning for Mass Transport Prediction in Porous Materials Using 90, 000 Artificially Generated Microstructures}},
	author       = {Prifling,  Benedikt and R\"{o}ding,  Magnus and Townsend,  Philip and Neumann,  Matthias and Schmidt,  Volker},
	year         = 2021,
	month        = dec,
	journal      = {Frontiers in Materials},
	publisher    = {Frontiers Media SA},
	volume       = 8,
	doi          = {10.3389/fmats.2021.786502},
	issn         = {2296-8016},
	url          = {http://dx.doi.org/10.3389/fmats.2021.786502}
}

@article{Saad1986,
	title        = {{GMRES: A Generalized Minimal Residual Algorithm for Solving Nonsymmetric Linear Systems}},
	author       = {Saad, Youcef and Schultz, Martin H.},
	year         = 1986,
	journal      = {SIAM Journal on Scientific and Statistical Computing},
	volume       = 7,
	number       = 3,
	pages        = {856--869},
	doi          = {10.1137/0907058},
	url          = {https://doi.org/10.1137/0907058}
}

@article{LoggWells2010,
	title        = {{{DOLFIN:} Automated Finite Element Computing}},
	author       = {Logg, Anders and Wells, Garth N.},
	year         = 2010,
	journal      = {{ACM} Transactions on Mathematical Software},
	volume       = 37,
	doi          = {10.1145/1731022.1731030}
}

@article{Adam,
	title        = {{Adam: A Method for Stochastic Optimization}},
	author       = {Kingma,  Diederik P. and Ba,  Jimmy},
	year         = 2015,
	journal      = {ICLR},
	publisher    = {arXiv},
	doi          = {10.48550/ARXIV.1412.6980},
	url          = {https://arxiv.org/abs/1412.6980},
	copyright    = {arXiv.org perpetual,  non-exclusive license}
}

@article{Darus,
	title        = {{2d microstructure data}},
	author       = {Lissner, Julian},
	year         = 2020,
	journal      = {DaRUS repository},
	publisher    = {DaRUS},
	doi          = {10.18419/darus-1151},
	url          = {https://doi.org/10.18419/darus-1151},
	version      = {V2}
}

@book{Saad1999,
	title        = {{On-Line Learning in Neural Networks}},
	author       = {Saad, David},
	year         = 1999,
	month        = jan,
	publisher    = {Cambridge University Press},
	doi          = {10.1017/cbo9780511569920},
	isbn         = 9780511569920,
	url          = {http://dx.doi.org/10.1017/CBO9780511569920}
}

@book{Saad2003,
	title        = {{Iterative Methods for Sparse Linear Systems}},
	author       = {Saad, Yousef},
	year         = 2003,
	month        = jan,
	publisher    = {Society for Industrial and Applied Mathematics},
	doi          = {10.1137/1.9780898718003},
	isbn         = 9780898718003,
	url          = {http://dx.doi.org/10.1137/1.9780898718003}
}

@article{Duester2008,
	title        = {{The finite cell method for three-dimensional problems of solid mechanics}},
	author       = {D\"{u}ster,  A. and Parvizian,  J. and Yang,  Z. and Rank,  E.},
	year         = 2008,
	month        = aug,
	journal      = {Computer Methods in Applied Mechanics and Engineering},
	publisher    = {Elsevier BV},
	volume       = 197,
	number       = {45–48},
	pages        = {3768–3782},
	doi          = {10.1016/j.cma.2008.02.036},
	issn         = {0045-7825},
	url          = {http://dx.doi.org/10.1016/j.cma.2008.02.036}
}

@article{Hestenes1952,
	title        = {{Methods of conjugate gradients for solving linear systems}},
	author       = {Hestenes,  M.R. and Stiefel,  E.},
	year         = 1952,
	month        = dec,
	journal      = {Journal of Research of the National Bureau of Standards},
	publisher    = {National Institute of Standards and Technology (NIST)},
	volume       = 49,
	number       = 6,
	pages        = 409,
	doi          = {10.6028/jres.049.044},
	issn         = {0091-0635},
	url          = {http://dx.doi.org/10.6028/jres.049.044}
}

@article{Melchers2024,
	title        = {{Neural Green's Operators for Parametric Partial Differential Equations}},
	author       = {Melchers,  Hugo and Prins,  Joost and Abdelmalik,  Michael},
	year         = 2024,
	journal      = {arXiv},
	publisher    = {arXiv},
	doi          = {10.48550/ARXIV.2406.01857},
	url          = {https://arxiv.org/abs/2406.01857},
	copyright    = {Creative Commons Attribution 4.0 International}
}

@misc{Li2025,
      title={Neural Preconditioning Operator for Efficient PDE Solves}, 
      author={Zhihao Li and Di Xiao and Zhilu Lai and Wei Wang},
      year={2025},
      eprint={2502.01337},
      archivePrefix={arXiv},
      primaryClass={cs.CE},
      url={https://arxiv.org/abs/2502.01337}, 
      howpublished={https://arxiv.org/abs/2502.01337}
}

@phdthesis{Xiang2022,
	title        = {{Solution of large linear systems with a massive number of right-hand sides and machine learning}},
	author       = {Xiang, Yan-Fei},
	year         = 2022,
	month        = dec,
	address      = {Bordeaux},
	number       = {2022BORD0383},
	url          = {https://theses.hal.science/tel-03967557},
	school       = {{Universit{\'e} de Bordeaux}},
	type         = {Theses},
	hal_id       = {tel-03967557},
	hal_version  = {v1}
}

@incollection{PyTorch,
	title        = {{PyTorch: An Imperative Style, High-Performance Deep Learning Library}},
	author       = {Paszke, Adam and Gross, Sam and Massa, Francisco and Lerer, Adam and Bradbury, James and Chanan, Gregory and Killeen, Trevor and Lin, Zeming and Gimelshein, Natalia and Antiga, Luca and Desmaison, Alban and Kopf, Andreas and Yang, Edward and DeVito, Zachary and Raison, Martin and Tejani, Alykhan and Chilamkurthy, Sasank and Steiner, Benoit and Fang, Lu and Bai, Junjie and Chintala, Soumith},
	year         = 2019,
	booktitle    = {Advances in Neural Information Processing Systems 32},
	publisher    = {Curran Associates, Inc.},
	pages        = {8024--8035},
	url          = {http://papers.neurips.cc/paper/9015-pytorch-an-imperative-style-high-performance-deep-learning-library.pdf}
}

@article{PyTorch2,
	title        = {{PyTorch 2: Faster Machine Learning Through Dynamic Python Bytecode Transformation and Graph Compilation}},
	author       = {Jason Ansel and Edward Yang and Horace He and Natalia Gimelshein and Animesh Jain and Michael Voznesensky and Bin Bao and David Berard and Geeta Chauhan and Anjali Chourdia and Will Constable and Alban Desmaison and Zachary DeVito and Elias Ellison and Will Feng and Jiong Gong and Michael Gschwind and Brian Hirsh and Sherlock Huang and Laurent Kirsch and Michael Lazos and Yanbo Liang and Jason Liang and Yinghai Lu and CK Luk and Bert Maher and Yunjie Pan and Christian Puhrsch and Matthias Reso and Mark Saroufim and Helen Suk and Michael Suo and Phil Tillet and Eikan Wang and Xiaodong Wang and William Wen and Shunting Zhang and Xu Zhao and Keren Zhou and Richard Zou and Ajit Mathews and Gregory Chanan and Peng Wu and Soumith Chintala},
	date         = {2024-02-06},
	langid       = {english}
}

@book{Hackbusch1994,
	title        = {{Iterative Solution of Large Sparse Systems of Equations}},
	author       = {Hackbusch,  Wolfgang},
	year         = 1994,
	journal      = {Applied Mathematical Sciences},
	publisher    = {Springer New York},
	doi          = {https://doi.org/10.1007/978-1-4612-4288-8},
	isbn         = 9781461242888,
	issn         = {0066-5452},
	url          = {http://dx.doi.org/10.1007/978-1-4612-4288-8}
}

@inproceedings{Li2020MultipoleGN,
author = {Li, Zongyi and Kovachki, Nikola and Azizzadenesheli, Kamyar and Liu, Burigede and Bhattacharya, Kaushik and Stuart, Andrew and Anandkumar, Anima},
title = {Multipole graph neural operator for parametric partial differential equations},
year = {2020},
isbn = {9781713829546},
publisher = {Curran Associates Inc.},
address = {Red Hook, NY, USA},
articleno = {567},
numpages = {12},
location = {Vancouver, BC, Canada},
series = {NIPS '20},
organization = {NIPS}
}

@article{Carslaw1925,
	title        = {{A historical note on Gibbs' phenomenon in Fourier's series and integrals}},
	author       = {H. S. Carslaw},
	year         = 1925,
	journal      = {Bulletin of the American Mathematical Society},
	publisher    = {American Mathematical Society},
	volume       = 31,
	number       = 8,
	pages        = {420 -- 424}
}

@book{Olver2014,
	title        = {{Introduction to Partial Differential Equations}},
	author       = {Olver,  Peter J.},
	year         = 2014,
	journal      = {Undergraduate Texts in Mathematics},
	publisher    = {Springer International Publishing},
	doi          = {10.1007/978-3-319-02099-0},
	isbn         = 9783319020990,
	issn         = {2197-5604},
	url          = {http://dx.doi.org/10.1007/978-3-319-02099-0}
}

@article{Raissi2019,
	title        = {{Physics-informed neural networks: {A} deep learning framework for solving forward and inverse problems involving nonlinear partial differential equations}},
	shorttitle   = {Physics-informed neural networks},
	author       = {Raissi, M. and Perdikaris, P. and Karniadakis, G. E.},
	year         = 2019,
	month        = feb,
	journal      = {Journal of Computational Physics},
	volume       = 378,
	pages        = {686--707},
	doi          = {10.1016/j.jcp.2018.10.045},
	issn         = {0021-9991},
	url          = {https://www.sciencedirect.com/science/article/pii/S0021999118307125},
	urldate      = {2023-08-30}
}

@inproceedings{
Li_fno_2021,
title={Fourier Neural Operator for Parametric Partial Differential Equations},
author={Zongyi Li and Nikola Borislavov Kovachki and Kamyar Azizzadenesheli and Burigede liu and Kaushik Bhattacharya and Andrew Stuart and Anima Anandkumar},
booktitle={International Conference on Learning Representations},
year={2021},
url={https://openreview.net/forum?id=c8P9NQVtmnO},
organization={ICLR}
}

@article{li_physics-informed-no_2023,
author = {Li, Zongyi and Zheng, Hongkai and Kovachki, Nikola and Jin, David and Chen, Haoxuan and Liu, Burigede and Azizzadenesheli, Kamyar and Anandkumar, Anima},
title = {Physics-Informed Neural Operator for Learning Partial Differential Equations},
year = {2024},
issue_date = {September 2024},
publisher = {Association for Computing Machinery},
address = {New York, NY, USA},
volume = {1},
number = {3},
url = {https://doi.org/10.1145/3648506},
doi = {10.1145/3648506},
journal = {ACM / IMS J. Data Sci.},
month = may,
articleno = {9},
numpages = {27}
}

@inproceedings{
Li2020,
title={Neural Operator: Graph Kernel Network for Partial Differential Equations},
author={Anima Anandkumar and Kamyar Azizzadenesheli and Kaushik Bhattacharya and Nikola Kovachki and Zongyi Li and Burigede Liu and Andrew Stuart},
booktitle={ICLR 2020 Workshop on Integration of Deep Neural Models and Differential Equations},
year={2019},
url={https://openreview.net/forum?id=fg2ZFmXFO3},
organization={ICLR}
}

@article{Lu2021,
	title        = {{Learning nonlinear operators via DeepONet based on the universal approximation theorem of operators}},
	author       = {Lu,  Lu and Jin,  Pengzhan and Pang,  Guofei and Zhang,  Zhongqiang and Karniadakis,  George Em},
	year         = 2021,
	month        = mar,
	journal      = {Nature Machine Intelligence},
	publisher    = {Springer Science and Business Media LLC},
	volume       = 3,
	number       = 3,
	pages        = {218--229},
	doi          = {10.1038/s42256-021-00302-5},
	issn         = {2522-5839},
	url          = {http://dx.doi.org/10.1038/s42256-021-00302-5}
}

@article{Kovachki2023,
author = {Kovachki, Nikola and Li, Zongyi and Liu, Burigede and Azizzadenesheli, Kamyar and Bhattacharya, Kaushik and Stuart, Andrew and Anandkumar, Anima},
title = {Neural operator: learning maps between function spaces with applications to PDEs},
year = {2023},
issue_date = {January 2023},
publisher = {JMLR.org},
volume = {24},
number = {1},
issn = {1532-4435},
journal = {J. Mach. Learn. Res.},
month = jan,
articleno = {89},
numpages = {97},
}

@article{
Zhao_incremental_2022,
title={Incremental Spatial and Spectral Learning of Neural Operators for Solving Large-Scale {PDE}s},
author={Robert Joseph George and Jiawei Zhao and Jean Kossaifi and Zongyi Li and Anima Anandkumar},
journal={Transactions on Machine Learning Research},
issn={2835-8856},
year={2024},
url={https://openreview.net/forum?id=xI6cPQObp0},
note={}
}

@article{Lissner2019,
	title        = {{Data-{Driven} {Microstructure} {Property} {Relations}}},
	author       = {Lissner, Julian and Fritzen, Felix},
	year         = 2019,
	month        = jun,
	journal      = {Mathematical and Computational Applications},
	volume       = 24,
	number       = 2,
	pages        = 57,
	doi          = {10.3390/mca24020057},
	issn         = {2297-8747},
	url          = {https://www.mdpi.com/2297-8747/24/2/57},
	urldate      = {2023-08-30},
	copyright    = {http://creativecommons.org/licenses/by/3.0/},
	language     = {en}
}

@article{Lissner2023,
	title        = {{Double {U}-{Net}: {Improved} multiscale modeling via fully convolutional neural networks}},
	shorttitle   = {Double {U}-{Net}},
	author       = {Lissner, Julian and Fritzen, Felix},
	year         = 2023,
	journal      = {PAMM},
	pages        = {e202300205},
	doi          = {10.1002/pamm.202300205},
	issn         = {1617-7061},
	url          = {https://onlinelibrary.wiley.com/doi/abs/10.1002/pamm.202300205},
	urldate      = {2023-10-13},
	language     = {en}
}

@book{Stoer2002,
	title        = {{Introduction to Numerical Analysis}},
	author       = {Stoer,  J. and Bulirsch,  R.},
	year         = 2002,
	journal      = {Texts in Applied Mathematics},
	publisher    = {Springer New York},
	doi          = {10.1007/978-0-387-21738-3},
	isbn         = 9780387217383,
	issn         = {0939-2475},
	url          = {http://dx.doi.org/10.1007/978-0-387-21738-3}
}

@article{FFTW05,
	title        = {{The Design and Implementation of {FFTW3}}},
	author       = {Frigo, Matteo and Johnson, Steven~G.},
	year         = 2005,
	journal      = {Proceedings of the IEEE},
	volume       = 93,
	number       = 2,
	pages        = {216--231}
}

@article{Leuschner2018,
	title        = {{Fourier-{Accelerated} {Nodal} {Solvers} ({FANS}) for homogenization problems}},
	author       = {Leuschner, Matthias and Fritzen, Felix},
	year         = 2018,
	month        = sep,
	journal      = {Computational Mechanics},
	volume       = 62,
	number       = 3,
	pages        = {359--392},
	doi          = {10.1007/s00466-017-1501-5},
	issn         = {1432-0924},
	url          = {https://doi.org/10.1007/s00466-017-1501-5},
	urldate      = {2023-08-10},
	language     = {en}
}

@article{Lippmann1950,
	title        = {{Variational Principles for Scattering Processes. I}},
	author       = {Lippmann, B. A. and Schwinger, Julian},
	year         = 1950,
	journal      = {Phys. Rev.},
	publisher    = {American Physical Society},
	volume       = 79,
	pages        = {469--480},
	doi          = {10.1103/PhysRev.79.469},
	url          = {https://link.aps.org/doi/10.1103/PhysRev.79.469},
	issue        = 3,
	numpages     = {0}
}

@article{Yvonnet2012,
	title        = {{A fast method for solving microstructural problems defined by digital images: a space Lippmann-Schwinger scheme}},
	author       = {Yvonnet, Julien},
	year         = 2012,
	journal      = {International Journal for Numerical Methods in Engineering},
	volume       = 92,
	number       = 2,
	pages        = {178--205},
	doi          = {https://doi.org/10.1002/nme.4334},
	url          = {https://onlinelibrary.wiley.com/doi/abs/10.1002/nme.4334}
}

@article{vondrejc_fft-based_2014,
	title        = {{An {FFT}-based {Galerkin} method for homogenization of periodic media}},
	author       = {Vondrejc, Jaroslav and Zeman, Jan and Marek, Ivo},
	year         = 2014,
	month        = aug,
	journal      = {Computers \& Mathematics with Applications},
	volume       = 68,
	number       = 3,
	pages        = {156--173},
	doi          = {10.1016/j.camwa.2014.05.014},
	issn         = {0898-1221},
	url          = {https://www.sciencedirect.com/science/article/pii/S0898122114002077},
	urldate      = {2023-08-11},
	language     = {en}
}

@article{zeman_finite_2017,
	title        = {{A finite element perspective on nonlinear {FFT}-based micromechanical simulations}},
	author       = {Zeman, J. and de Geus, T. W. J. and Vondrejc, J. and Peerlings, R. H. J. and Geers, M. G. D.},
	year         = 2017,
	journal      = {International Journal for Numerical Methods in Engineering},
	volume       = 111,
	number       = 10,
	pages        = {903--926},
	doi          = {10.1002/nme.5481},
	issn         = {1097-0207},
	url          = {https://onlinelibrary.wiley.com/doi/abs/10.1002/nme.5481},
	urldate      = {2023-08-11},
	copyright    = {Copyright _ 2016 John Wiley \& Sons, Ltd.},
	language     = {en}
}

@article{Schneider2017,
	title        = {{{FFT}-based homogenization for microstructures discretized by linear hexahedral elements}},
	author       = {Schneider, Matti and Merkert, Dennis and Kabel, Matthias},
	year         = 2017,
	journal      = {International Journal for Numerical Methods in Engineering},
	volume       = 109,
	number       = 10,
	pages        = {1461--1489},
	doi          = {10.1002/nme.5336},
	issn         = {1097-0207},
	url          = {https://onlinelibrary.wiley.com/doi/abs/10.1002/nme.5336},
	urldate      = {2023-08-20},
	copyright    = {Copyright _ 2016 John Wiley \& Sons, Ltd.},
	language     = {en}
}

@article{Brunton2023,
	title        = {{Promising directions of machine learning for partial differential equations}},
	author       = {Brunton,  Steven L. and Kutz,  J. Nathan},
	year         = 2024,
	month        = jun,
	journal      = {Nature Computational Science},
	publisher    = {Springer Science and Business Media LLC},
	volume       = 4,
	number       = 7,
	pages        = {483–494},
	doi          = {10.1038/s43588-024-00643-2},
	issn         = {2662-8457},
	url          = {http://dx.doi.org/10.1038/s43588-024-00643-2}
}

@article{Zeman2010,
	title        = {{Accelerating a {FFT}-based solver for numerical homogenization of periodic media by conjugate gradients}},
	author       = {Zeman, Jan and Vondrejc, Jaroslav and Novak, Jan and Marek, Ivo},
	year         = 2010,
	month        = oct,
	journal      = {Journal of Computational Physics},
	volume       = 229,
	number       = 21,
	pages        = {8065--8071},
	doi          = {10.1016/j.jcp.2010.07.010},
	issn         = {00219991},
	url          = {https://linkinghub.elsevier.com/retrieve/pii/S0021999110003931},
	urldate      = {2023-08-31},
	language     = {en}
}

@book{Golub2013,
	title        = {{Matrix computations}},
	author       = {Golub, G.H. and Van Loan, C.F.},
	year         = 2013,
	publisher    = {Johns Hopkins University Press},
	series       = {Johns hopkins studies in the mathematical sciences},
	isbn         = {978-1-4214-0794-4},
	url          = {https://books.google.de/books?id=X5YfsuCWpxMC},
	collection   = {Johns hopkins studies in the mathematical sciences}
}

@article{Moulinec1998,
	title        = {{A numerical method for computing the overall response of nonlinear composites with complex microstructure}},
	author       = {Moulinec, H. and Suquet, P.},
	year         = 1998,
	month        = apr,
	journal      = {Computer Methods in Applied Mechanics and Engineering},
	volume       = 157,
	number       = 1,
	pages        = {69--94},
	doi          = {10.1016/S0045-7825(97)00218-1},
	issn         = {0045-7825},
	url          = {https://www.sciencedirect.com/science/article/pii/S0045782597002181},
	urldate      = {2023-09-04}
}

@article{Moulinec2014,
	title        = {{Comparison of three accelerated {FFT}-based schemes for computing the mechanical response of composite materials}},
	author       = {Moulinec, H. and Silva, F.},
	year         = 2014,
	journal      = {International Journal for Numerical Methods in Engineering},
	volume       = 97,
	number       = 13,
	pages        = {960--985},
	doi          = {10.1002/nme.4614},
	issn         = {1097-0207},
	url          = {https://onlinelibrary.wiley.com/doi/abs/10.1002/nme.4614},
	urldate      = {2023-09-07},
	copyright    = {Copyright _ 2014 John Wiley \& Sons, Ltd.},
	language     = {en}
}

@article{Keshav2022,
	title        = {{FFT-based homogenization at finite strains using composite boxels (ComBo)}},
	author       = {Keshav,  Sanath and Fritzen,  Felix and Kabel,  Matthias},
	year         = 2022,
	month        = oct,
	journal      = {Computational Mechanics},
	publisher    = {Springer Science and Business Media LLC},
	volume       = 71,
	number       = 1,
	pages        = {191–212},
	doi          = {10.1007/s00466-022-02232-4},
	issn         = {1432-0924},
	url          = {http://dx.doi.org/10.1007/s00466-022-02232-4}
}

@article{Schneider2019,
	title        = {{On polarization-based schemes for the {FFT}-based computational homogenization of inelastic materials}},
	author       = {Schneider, Matti and Wicht, Daniel and Bohlke, Thomas},
	year         = 2019,
	month        = oct,
	journal      = {Computational Mechanics},
	volume       = 64,
	number       = 4,
	pages        = {1073--1095},
	doi          = {10.1007/s00466-019-01694-3},
	issn         = {1432-0924},
	url          = {https://doi.org/10.1007/s00466-019-01694-3},
	urldate      = {2023-10-25},
	language     = {en}
}

@article{Schneider2021,
	title        = {{A review of nonlinear FFT-based computational homogenization methods}},
	author       = {Schneider,  Matti},
	year         = 2021,
	month        = mar,
	journal      = {Acta Mechanica},
	publisher    = {Springer Science and Business Media LLC},
	volume       = 232,
	number       = 6,
	pages        = {2051–2100},
	doi          = {10.1007/s00707-021-02962-1},
	issn         = {1619-6937},
	url          = {http://dx.doi.org/10.1007/s00707-021-02962-1}
}

@article{Lucarini2021,
	title        = {{{FFT} based approaches in micromechanics: fundamentals, methods and applications}},
	shorttitle   = {{FFT} based approaches in micromechanics},
	author       = {Lucarini, S. and Upadhyay, M. V. and Segurado, J.},
	year         = 2021,
	month        = dec,
	journal      = {Modelling and Simulation in Materials Science and Engineering},
	volume       = 30,
	number       = 2,
	pages        = {023002},
	doi          = {10.1088/1361-651X/ac34e1},
	issn         = {0965-0393},
	url          = {https://dx.doi.org/10.1088/1361-651X/ac34e1},
	urldate      = {2023-11-03},
	language     = {en}
}

@article{Ladecky2023,
	title        = {{An optimal preconditioned {FFT}-accelerated finite element solver for homogenization}},
	author       = {Ladecky, Martin and Leute, Richard J. and Falsafi, Ali and Pultarova, Ivana and Pastewka, Lars and Junge, Till and Zeman, Jan},
	year         = 2023,
	month        = jun,
	journal      = {Applied Mathematics and Computation},
	volume       = 446,
	pages        = 127835,
	doi          = {10.1016/j.amc.2023.127835},
	issn         = {0096-3003},
	url          = {https://www.sciencedirect.com/science/article/pii/S0096300323000048},
	urldate      = {2024-01-16}
}

@article{Moulinec1994,
	title        = {{A fast numerical method for computing the linear and nonlinear mechanical properties of composites}},
	author       = {Moulinec, H. and Suquet, Pierre},
	year         = 1994,
	month        = apr,
	journal      = {Comptes Rendus de l'Academie des sciences. Serie II. Mecanique, physique, chimie, astronomie},
	url          = {https://hal.science/hal-03019226},
	urldate      = {2024-01-29}
}

@inbook{Petsc1997,
	title        = {{Efficient Management of Parallelism in Object-Oriented Numerical Software Libraries}},
	author       = {Balay,  Satish and Gropp,  William D. and McInnes,  Lois Curfman and Smith,  Barry F.},
	year         = 1997,
	booktitle    = {Modern Software Tools for Scientific Computing},
	publisher    = {Birkh\"{a}user Boston},
	pages        = {163–202},
	doi          = {10.1007/978-1-4612-1986-6_8},
	isbn         = 9781461219866,
	url          = {http://dx.doi.org/10.1007/978-1-4612-1986-6_8}
}

@book{Barrett1994,
  title = {Templates for the Solution of Linear Systems: Building Blocks for Iterative Methods},
  ISBN = {9781611971538},
  url = {http://dx.doi.org/10.1137/1.9781611971538},
  DOI = {10.1137/1.9781611971538},
  publisher = {Society for Industrial and Applied Mathematics},
  author = {Barrett,  Richard and Berry,  Michael and Chan,  Tony F. and Demmel,  James and Donato,  June and Dongarra,  Jack and Eijkhout,  Victor and Pozo,  Roldan and Romine,  Charles and van der Vorst,  Henk},
  year = {1994},
  month = jan 
}

@misc{UNOCGgithub2025,
	title        = {{UNO-CG: Accelerating Conjugate Gradient Solvers with Unitary Neural Operators (software package)}},
	author       = {Herb, Julius and Fritzen, Felix},
	journal      = {GitHub repository},
	publisher    = {GitHub},
	url          = {https://github.com/DataAnalyticsEngineering/UNOCG},
	howpublished = {\url{https://github.com/DataAnalyticsEngineering/UNOCG}},
	version      = {1.0.0},
	date         = {2026-01-26},
    year         = 2026,
}

@misc{UNOCGdarus2025,
	title        = {{Supplemental data for "Accelerating Conjugate Gradient Solvers for Homogenization Problems with Unitary Neural Operators"}},
	author       = {Herb, Julius and Fritzen, Felix},
	year         = 2026,
	journal      = {DaRUS repository},
	publisher    = {DaRUS},
	howpublished = {\url{https://doi.org/10.18419/DARUS-5686}},
	version      = {1},
	date         = {2026-01-26}
}

@manual{FiredrakeUserManual,
	title        = {{Firedrake User Manual}},
	author       = {David A. Ham and Paul H. J. Kelly and Lawrence Mitchell and Colin J. Cotter and Robert C. Kirby and Koki Sagiyama and Nacime Bouziani and Sophia Vorderwuelbecke and Thomas J. Gregory and Jack Betteridge and Daniel R. Shapero and Reuben W. Nixon-Hill and Connor J. Ward and Patrick E. Farrell and Pablo D. Brubeck and India Marsden and Thomas H. Gibson and Miklos Homolya and Tianjiao Sun and Andrew T. T. McRae and Fabio Luporini and Alastair Gregory and Michael Lange and Simon W. Funke and Florian Rathgeber and Gheorghe-Teodor Bercea and Graham R. Markall},
	year         = 2023,
	month        = may,
	address      = {London},
	doi          = {10.25561/104839},
	organization = {Imperial College London and University of Oxford and Baylor University and University of Washington},
	edition      = {First edition}
}

@article{Mercer2015,
	title        = {{Novel formulations of microscopic boundary-value problems in continuous multiscale finite element methods}},
	author       = {Mercer, Brian S. and Mandadapu, Kranthi K. and Papadopoulos, Panayiotis},
	year         = 2015,
	month        = apr,
	journal      = {Computer Methods in Applied Mechanics and Engineering},
	volume       = 286,
	pages        = {268--292},
	doi          = {10.1016/j.cma.2014.12.021},
	issn         = {0045-7825}
}

@inbook{Schroder2014,
	title        = {{A numerical two-scale homogenization scheme: the FE2-method}},
	author       = {Schr\"{o}der,  J\"{o}rg},
	year         = 2014,
	booktitle    = {Plasticity and Beyond},
	publisher    = {Springer Vienna},
	pages        = {1–64},
	doi          = {10.1007/978-3-7091-1625-8\_1},
	isbn         = 9783709116258,
	issn         = {2309-3706},
	url          = {http://dx.doi.org/10.1007/978-3-7091-1625-8\%5F1}
}

@book{Yvonnet2019,
	title        = {{Computational Homogenization of Heterogeneous Materials with Finite Elements}},
	author       = {Yvonnet,  Julien},
	year         = 2019,
	journal      = {Solid Mechanics and Its Applications},
	publisher    = {Springer International Publishing},
	doi          = {10.1007/978-3-030-18383-7},
	isbn         = 9783030183837,
	issn         = {2214-7764},
	url          = {http://dx.doi.org/10.1007/978-3-030-18383-7}
}

@book{Chawla2013,
	title        = {{Metal Matrix Composites}},
	author       = {Chawla,  Nikhilesh and Chawla,  Krishan K.},
	year         = 2013,
	publisher    = {Springer New York},
	doi          = {10.1007/978-1-4614-9548-2},
	isbn         = 9781461495482,
	url          = {http://dx.doi.org/10.1007/978-1-4614-9548-2}
}

@inbook{Dirrenberger2019,
	title        = {{Computational Homogenization of Architectured Materials}},
	author       = {Dirrenberger,  Justin and Forest,  Samuel and Jeulin,  Dominique},
	year         = 2019,
	booktitle    = {Architectured Materials in Nature and Engineering},
	publisher    = {Springer International Publishing},
	pages        = {89–139},
	doi          = {10.1007/978-3-030-11942-3\_4},
	isbn         = 9783030119423,
	issn         = {2196-2812},
	url          = {http://dx.doi.org/10.1007/978-3-030-11942-3\%5F4}
}

@article{Risthaus2024,
	title        = {{Imposing different boundary conditions for thermal computational homogenization problems with FFT- and tensor-train-based Green's operator methods}},
	author       = {Risthaus, Lennart and Schneider, Matti},
	year         = 2024,
	journal      = {International Journal for Numerical Methods in Engineering},
	volume       = 125,
	number       = 7,
	pages        = {e7423},
	doi          = {10.1002/nme.7423},
	issn         = {1097-0207},
	rights       = {2024 The Authors. International Journal for Numerical Methods in Engineering published by John Wiley & Sons Ltd.},
	language     = {en}
}

@article{Song2024,
	title        = {{High-resolution simulating of grain substructure in cold rolling and its effects on primary recrystallization in annealing of ferritic stainless steel}},
	author       = {Kangjie Song and Haochen Ding and Chi Zhang and Liwen Zhang and Guanyu Deng and Huaibei Zheng},
	year         = 2024,
	journal      = {Journal of Materials Research and Technology},
	volume       = 30,
	pages        = {40--51},
	doi          = {https://doi.org/10.1016/j.jmrt.2024.03.065},
	issn         = {2238-7854},
	url          = {https://www.sciencedirect.com/science/article/pii/S2238785424005969}
}

@article{Aboudi2004,
	title        = {{The Generalized Method of Cells and High-Fidelity Generalized Method of Cells Micromechanical Models—A Review}},
	author       = {Aboudi, Jacob},
	year         = 2004,
	month        = jul,
	journal      = {Mechanics of Advanced Materials and Structures},
	publisher    = {Informa UK Limited},
	volume       = 11,
	number       = {4–5},
	pages        = {329–366},
	doi          = {10.1080/15376490490451543},
	issn         = {1537-6532},
	url          = {http://dx.doi.org/10.1080/15376490490451543}
}

@book{Bagchi1999,
	title        = {{The Nonuniform Discrete Fourier Transform and Its Applications in Signal Processing}},
	author       = {Bagchi,  Sonali and Mitra,  Sanjit K.},
	year         = 1999,
	publisher    = {Springer US},
	doi          = {10.1007/978-1-4615-4925-3},
	isbn         = 9781461549253,
	url          = {http://dx.doi.org/10.1007/978-1-4615-4925-3}
}

@article{Pooja2025,
	title        = {{Metal matrix composites: revolutionary materials for shaping the future}},
	author       = {Pooja,  Km. and Tarannum,  Nazia and Chaudhary,  Pallavi},
	year         = 2025,
	month        = feb,
	journal      = {Discover Materials},
	publisher    = {Springer Science and Business Media LLC},
	volume       = 5,
	number       = 1,
	doi          = {10.1007/s43939-025-00226-6},
	issn         = {2730-7727},
	url          = {http://dx.doi.org/10.1007/s43939-025-00226-6}
}

@article{Osanov2016,
	title        = {{Topology Optimization for Architected Materials Design}},
	author       = {Osanov,  Mikhail and Guest,  James K.},
	year         = 2016,
	month        = jul,
	journal      = {Annual Review of Materials Research},
	publisher    = {Annual Reviews},
	volume       = 46,
	number       = 1,
	pages        = {211–233},
	doi          = {10.1146/annurev-matsci-070115-031826},
	issn         = {1545-4118},
	url          = {http://dx.doi.org/10.1146/annurev-matsci-070115-031826}
}

@article{Zhang2024,
	title        = {{Blending neural operators and relaxation methods in PDE numerical solvers}},
	author       = {Zhang,  Enrui and Kahana,  Adar and Kopaničáková,  Alena and Turkel,  Eli and Ranade,  Rishikesh and Pathak,  Jay and Karniadakis,  George Em},
	year         = 2024,
	month        = oct,
	journal      = {Nature Machine Intelligence},
	publisher    = {Springer Science and Business Media LLC},
	volume       = 6,
	number       = 11,
	pages        = {1303–1313},
	doi          = {10.1038/s42256-024-00910-x},
	issn         = {2522-5839},
	url          = {http://dx.doi.org/10.1038/s42256-024-00910-x}
}

@article{Kopanicakova2025,
	title        = {{DeepONet Based Preconditioning Strategies for Solving Parametric Linear Systems of Equations}},
	author       = {Kopaničáková,  Alena and Karniadakis,  George Em},
	year         = 2025,
	month        = feb,
	journal      = {SIAM Journal on Scientific Computing},
	publisher    = {Society for Industrial & Applied Mathematics (SIAM)},
	volume       = 47,
	number       = 1,
	pages        = {C151–C181},
	doi          = {10.1137/24m162861x},
	issn         = {1095-7197},
	url          = {http://dx.doi.org/10.1137/24M162861X}
}

@article{Cuomo2022,
	title        = {{Scientific Machine Learning Through Physics--Informed Neural Networks: Where we are and What’s Next}},
	author       = {Cuomo,  Salvatore and Di Cola,  Vincenzo Schiano and Giampaolo,  Fabio and Rozza,  Gianluigi and Raissi,  Maziar and Piccialli,  Francesco},
	year         = 2022,
	month        = jul,
	journal      = {Journal of Scientific Computing},
	publisher    = {Springer Science and Business Media LLC},
	volume       = 92,
	number       = 3,
	doi          = {10.1007/s10915-022-01939-z},
	issn         = {1573-7691},
	url          = {http://dx.doi.org/10.1007/s10915-022-01939-z}
}

@article{Watson2025,
	title        = {{Machine Learning with Physics Knowledge for Prediction: A Survey}},
	author       = {Joe Watson and Chen Song and Oliver Weeger and Theo Gruner and An Thai Le and Kay Hansel and Ahmed Hendawy and Oleg Arenz and Will Trojak and Miles Cranmer and Carlo D'Eramo and Fabian Buelow and Tanmay Goyal and Jan Peters and Martin W Hoffmann},
	year         = 2025,
	journal      = {Transactions on Machine Learning Research},
	issn         = {2835-8856},
	url          = {https://openreview.net/forum?id=ZiJYahyXLU}
}

@article{Li2022,
	title        = {{Fourier Neural Operator with Learned Deformations for PDEs on General Geometries}},
	author       = {Li,  Zongyi and Huang,  Daniel Zhengyu and Liu,  Burigede and Anandkumar,  Anima},
	year         = 2023,
	journal      = {Journal of Machine Learning Research},
	publisher    = {arXiv},
	doi          = {10.48550/ARXIV.2207.05209},
	url          = {https://arxiv.org/abs/2207.05209},
	copyright    = {arXiv.org perpetual,  non-exclusive license}
}

@article{Kaasschieter1988,
	title        = {{Preconditioned conjugate gradients for solving singular systems}},
	author       = {Kaasschieter,  E.F.},
	year         = 1988,
	month        = nov,
	journal      = {Journal of Computational and Applied Mathematics},
	publisher    = {Elsevier BV},
	volume       = 24,
	number       = {1–2},
	pages        = {265–275},
	doi          = {10.1016/0377-0427(88)90358-5},
	issn         = {0377-0427},
	url          = {http://dx.doi.org/10.1016/0377-0427(88)90358-5}
}

\appendix

\bmsection{Additional notational details}\label{sec:appendix-notation}

The plane strain assumption is employed for all mechanical two-dimensional problems, i.e.,
\begin{equation} \label{eq:plain-strain}
    \strainComp_{33} = \strainComp_{13} = \strainComp_{23} = 0 \,.
\end{equation}
Symmetric tensors can be expressed in Mandel notation using the orthonormal basis that is formed by
\begin{align} \label{eq:mandel-basis}
    \mandelBasis{1} = \unitVec_{(1)} \otimes \unitVec_{(1)} \,, &&
    \mandelBasis{2} = \unitVec_{(2)} \otimes \unitVec_{(2)} \,, &&
    \mandelBasis{3} = \sqrt{2} \, \symOp{ \unitVec_{(1)} \otimes \unitVec_{(2)} } \,,
\end{align}
where~$\otimes$ denotes the tensor product.
For three-dimensional problems, we introduce the Mandel notation via the basis
\begin{align}
	&\mandelBasis{1} = \unitVec_{(1)} \otimes \unitVec_{(1)} \,,
	&&\mandelBasis{2} = \unitVec_{(2)} \otimes \unitVec_{(2)} \,,
	&&\mandelBasis{3} = \unitVec_{(3)} \otimes \unitVec_{(3)} \,, \\
	&\mandelBasis{4} = \sqrt{2} \, \symOp{ \unitVec_{(1)} \otimes \unitVec_{(2)} } \,,
	&&\mandelBasis{5} = \sqrt{2} \, \symOp{ \unitVec_{(1)} \otimes \unitVec_{(3)} } \,,
	&&\mandelBasis{6} = \sqrt{2} \, \symOp{ \unitVec_{(2)} \otimes \unitVec_{(3)} } \,.
\end{align}
In Mandel notation, symmetric tensors can be expressed as~$\strain \leftrightarrow \strainMandel \in \ffR^D$, $\stress \leftrightarrow \stressMandel \in \ffR^D$, $D = \nDims (\nDims + 1) / 2$.
We define $\PIso_1 = \frac{1}{3} \fI \otimes \fI$ and $\PIso_2 = \Isym - \PIso_1$ based on the identity~$\fI \in \ffR^{\nDims \times \nDims}$ and the identity on symmetric second-order tensors~$\Isym \in \ffR^{\nDims \times \nDims \times \nDims \times \nDims}$.

\bmsection{Efficient training of \texorpdfstring{\UNOCG}{UNO-CG} for three nodal DOF (\texorpdfstring{$\nComp=3$}{c=3})}

\bmsubsection{Local parametrization}
\label{sec:local-parametrization-c-3}
We extend the local parametrization from~\eqref{eq:local-param-2} to~$\localParam: \ffR^{6} \to \spdSpace{\ffR^{3 \times 3}}$, which is defined as
\begin{align} \label{eq:local-param-3}
	\localParam \left( \weights\modeIdx \right) =
    \begin{bmatrix}
 		\weight\modeIdx_1 \cdot \weight\modeIdx_1 & \weight\modeIdx_1 \cdot \weight\modeIdx_3 & \weight\modeIdx_1 \cdot \weight\modeIdx_6 \\
 		\weight\modeIdx_1 \cdot \weight\modeIdx_3 & \weight\modeIdx_2 \cdot \weight\modeIdx_2 + \weight\modeIdx_3 \cdot \weight\modeIdx_3 & \weight\modeIdx_2 \cdot \weight\modeIdx_5 + \weight\modeIdx_3 \cdot \weight\modeIdx_6 \\
 		\weight\modeIdx_1 \cdot \weight\modeIdx_6 & \weight\modeIdx_2 \cdot \weight\modeIdx_5 + \weight\modeIdx_3 \cdot \weight\modeIdx_6 & \weight\modeIdx_4 \cdot \weight\modeIdx_4 + \weight\modeIdx_5 \cdot \weight\modeIdx_5 + \weight\modeIdx_6 \cdot \weight\modeIdx_6
 	\end{bmatrix} \,, &&
	0 \leq \modeI \leq \nModes \,.
\end{align}

\bmsubsection{Computation of loss and derivatives}
\label{sec:training-c-3}
We introduce a split of the DOF vectors~$\PrecInDSample, \PrecOutDSample, \dots \in \ffR^{3n}$ into vectors containing the nodal DOF~$\PrecInDXSample, \PrecInDYSample, \PrecInDZSample, \dots \in \ffR^{n}$ as in
\begin{align}
	\PrecInDSample
	= \T{\begin{bmatrix} \left(\PrecInDXSample\right)_1 & \left(\PrecInDYSample\right)_1 & \left(\PrecInDZSample\right)_1 & \cdots & \left(\PrecInDXSample\right)_n & \left(\PrecInDYSample\right)_n & \left(\PrecInDZSample\right)_n \end{bmatrix}} \in \ffR^{3 n} \,.
\end{align}
The symmetric band matrix~$\PrecSparse$ can then be written as
\begin{align}
	\PrecSparse =
	\begin{bmatrix}
		\diag{\PrecEntries_1} & \diag{\PrecEntries_3} & \diag{\PrecEntries_6} \\
		\diag{\PrecEntries_3} & \diag{\PrecEntries_2} & \diag{\PrecEntries_5} \\
		\diag{\PrecEntries_6} & \diag{\PrecEntries_5} & \diag{\PrecEntries_4}
	\end{bmatrix} \,,
	&& \PrecEntries = \T{\begin{bmatrix} \T{\PrecEntries}_1 \quad \T{\PrecEntries}_2 \quad \T{\PrecEntries}_3 \quad  \T{\PrecEntries}_4 \quad \T{\PrecEntries}_5 \quad \T{\PrecEntries}_6 \end{bmatrix}} \in \ffR^{6n} \,.
\end{align}
Following this ansatz, it is possible to compute the loss $\Loss$ and its derivates with respect to $\PrecEntries$ via
\begin{align*}
	\Loss =
	&\featureA_{1} \cdot \left( \PrecEntries_1 \odot \PrecEntries_1 + \PrecEntries_3 \odot \PrecEntries_3 + \PrecEntries_6 \odot \PrecEntries_6 \right)
	+ \featureA_{2} \cdot \left( \PrecEntries_2 \odot \PrecEntries_2 + \PrecEntries_3 \odot \PrecEntries_3 + \PrecEntries_5 \odot \PrecEntries_5 \right)
	+ \featureA_{3} \cdot \left( \PrecEntries_1 \odot \PrecEntries_3 + \PrecEntries_2 \odot \PrecEntries_3 + \PrecEntries_5 \odot \PrecEntries_6 \right) \\
	&+ \featureA_{4} \cdot \left( \PrecEntries_4 \odot \PrecEntries_4 + \PrecEntries_5 \odot \PrecEntries_5 + \PrecEntries_6 \odot \PrecEntries_6 \right)
	+ \featureA_{5} \cdot \left( \PrecEntries_2 \odot \PrecEntries_5 + \PrecEntries_3 \odot \PrecEntries_6 + \PrecEntries_4 \odot \PrecEntries_5 \right)
	+ \featureA_{6} \cdot \left( \PrecEntries_1 \odot \PrecEntries_6 + \PrecEntries_3 \odot \PrecEntries_5 + \PrecEntries_4 \odot \PrecEntries_6 \right) \\
	&- \featureB_{1} \cdot \PrecEntries_1
	- \featureB_{2} \cdot \PrecEntries_2
	- \featureB_{3} \cdot \PrecEntries_3
	- \featureB_{4} \cdot \PrecEntries_4
	- \featureB_{5} \cdot \PrecEntries_5
	- \featureB_{6} \cdot \PrecEntries_6
	+ \featureDelta \in \ffR \,,
\end{align*}
\begin{align*}
	\frac{\partial \Loss}{\partial \PrecEntries} =& \T{\begin{bmatrix}
			2 \featureA_{1} \odot \PrecEntries_1
			+ \featureA_{3} \odot \PrecEntries_3
			+ \featureA_{6} \odot \PrecEntries_6
			- \featureB_{1} \\
			2 \featureA_{2} \odot \PrecEntries_2
			+ \featureA_{3} \odot \PrecEntries_3
			+ \featureA_{5} \odot \PrecEntries_5
			- \featureB_{2} \\
			2 \featureA_{1} \odot \PrecEntries_3
			+ 2 \featureA_{2} \odot \PrecEntries_3
			+ \featureA_{3} \odot \left( \PrecEntries_1 + \PrecEntries_2 \right)
			+ \featureA_{5} \odot \PrecEntries_6
			+ \featureA_{6} \odot \PrecEntries_5
			- \featureB_{3} \\
			2 \featureA_{4} \odot \PrecEntries_4
			+ \featureA_{5} \odot \PrecEntries_5
			+ \featureA_{6} \odot \PrecEntries_6
			- \featureB_{4} \\
			2 \featureA_{2} \odot \PrecEntries_5
			+ \featureA_{3} \odot \PrecEntries_6
			+ 2 \featureA_{4} \odot \PrecEntries_5
			+ \featureA_{5} \odot \left( \PrecEntries_2 + \PrecEntries_4 \right)
			+ \featureA_{6} \odot \PrecEntries_3
			- \featureB_{5} \\
			2 \featureA_{1} \odot \PrecEntries_6
			+ \featureA_{3} \odot \PrecEntries_5
			+ 2 \featureA_{4} \odot \PrecEntries_6
			+ \featureA_{5} \odot \PrecEntries_3
			+ \featureA_{6} \odot \left( \PrecEntries_1 + \PrecEntries_4 \right)
			- \featureB_{6}
	\end{bmatrix}} \in \ffR^{6n} \,,
\end{align*}
\begin{align*}
	\frac{\partial^2 \Loss}{\partial \PrecEntries\partial \PrecEntries} =& \begin{bmatrix}
		2 \diag{\featureA_1} & \ull{0} & \diag{\featureA_3} & \ull{0} & \ull{0} & \diag{\featureA_6} \\
		\ull{0} & 2 \diag{\featureA_2} & \diag{\featureA_3} & \ull{0} & \diag{\featureA_5} & \ull{0} \\
		\diag{\featureA_3} & \diag{\featureA_3} & 2 \diag{\featureA_1 + \featureA_2} & \ull{0} & \diag{\featureA_6} & \diag{\featureA_5} \\
		\ull{0} & \ull{0} & \ull{0} & 2 \diag{\featureA_4} & \diag{\featureA_5} & \diag{\featureA_6} \\
		\ull{0} & \diag{\featureA_5} & \diag{\featureA_6} & \diag{\featureA_5} & 2 \diag{\featureA_2 + \featureA_4} & \diag{\featureA_3} \\
		\diag{\featureA_6} & \ull{0} & 
		\diag{\featureA_5} & \diag{\featureA_6} & \diag{\featureA_3} & 2 \diag{\featureA_1 + \featureA_4}
	\end{bmatrix} \in \ffR^{6n \times 6n} \,,
\end{align*}
with the previous features from \cref{sssec:special2dof} and the additional features $\featureA_{\rm 4}, \featureA_{\rm 5}, \featureA_{\rm 6}, \featureB_{\rm 4}, \featureB_{\rm 5}, \featureB_{\rm 6} \in \ffR^n$ as
\begin{align*}
    \featureA_{\rm 4} &= \frac{1}{\nSamples} \sum_{\sampleI = 1}^\nSamples \PrecInDAdjZSample \odot \PrecInDFwdZSample \,, &&
    \featureA_{\rm 5} = \frac{1}{\nSamples} \sum_{\sampleI = 1}^\nSamples \PrecInDAdjYSample \odot \PrecInDFwdZSample + \PrecInDAdjZSample \odot \PrecInDFwdYSample \,, &&
    \featureA_{\rm 6} = \frac{1}{\nSamples} \sum_{\sampleI = 1}^\nSamples \PrecInDAdjXSample \odot \PrecInDFwdZSample + \PrecInDAdjZSample \odot \PrecInDFwdXSample \,, \\
    \featureB_{\rm 4} &= \frac{2}{\nSamples} \sum_{\sampleI = 1}^\nSamples \realPart{ \PrecInDAdjZSample \odot \PrecOutDFwdZSample } \,, &&
    \featureB_{\rm 5} = \frac{2}{\nSamples} \sum_{\sampleI = 1}^\nSamples \realPart{ \PrecInDAdjYSample \odot \PrecOutDFwdZSample + \PrecInDAdjZSample \odot \PrecOutDFwdYSample } \,, &&
    \featureB_{\rm 6} = \frac{2}{\nSamples} \sum_{\sampleI = 1}^\nSamples \realPart{ \PrecInDAdjXSample \odot \PrecOutDFwdZSample + \PrecInDAdjZSample \odot \PrecOutDFwdXSample } .
\end{align*}

\bmsection{Mathematical proofs}
\label{sec:appendix-proofs}

\begin{proof}[Proof of \cref{spectrum-lemma}]
By exploiting the property of the unitary matrix $\trafoD$ and by virtue of a similarity transformation, it follows
\begin{align*}
	\eig{\PrecD_\weights} = \eig{\trafoDHer \PrecSparse \, \trafoD} = \eig{\PrecSparse} \,.
\end{align*}
Due to the special structure of~$\PrecSparse$ in \eqref{eq:prec-sparse-def}, there exists a permutation matrix~$\ull{R}_\pi$ that transforms $\PrecSparse$ into the block-diagonal matrix
\begin{align*}
    \widetilde{\ull{Q}} = \ull{R}_\pi \ull{Q} \, \T{\ull{R}}_\pi = \mathrm{blockdiag}\left( \fnocgFundSolMat\modeIdx \right) = \begin{bmatrix} \fnocgFundSolMat\modeIdxFirst & \ull{0} & \cdots & \ull{0} \\ \ull{0} & \fnocgFundSolMat^{\langle 2 \rangle} & \cdots & \ull{0} \\ \vdots & \vdots & \ddots & \vdots \\ \ull{0} & \ull{0} & \cdots & \fnocgFundSolMat^{\langle n \rangle} \end{bmatrix} \in \ffR^{\nComp n \times \nComp n} \,, && \fnocgFundSolMat\modeIdx \in \ffR^{\nComp \times \nComp} \,, && 1 \leq \modeI \leq n \,.
\end{align*}
Since permutation matrices are orthogonal, this is also a similarity transformation. The eigenvalues of $\widetilde{\ull{Q}}$ are directly available:
\begin{align*}
	\eig{\PrecSparse} = \eig{\T{\ull{R}}_\pi \widetilde{\ull{Q}} \, \ull{R}_\pi} = \eig{\widetilde{\ull{Q}}} = {\color{revision}\mathrm{eig}\left(\mathrm{blockdiag}\left( \fnocgFundSolMat\modeIdx \right) \right)} = \bigcup\limits_{i=1}^{n} \eig{\fnocgFundSolMat\modeIdx} \,.
\end{align*}
\end{proof}

\clearpage

\bmsection{Information about training and test data}
\label{ssec:hyperparameters}

\begin{table}[ht!]
    \caption{Training and test data considered for the different problems. All data is available in our data repository \cite{UNOCGdarus2025}, where also additional details about data generation can be found.
    }
\begin{adjustbox}{center}
    \centering \begin{tabular}{ccccccc}
        \toprule 
        Problem & Training data & Test data \\
        \midrule 
        Thermal 2D (periodic BC) & \makecell{$\solDParam, \rhsDParam$ for 1000 2D microstructures ($400^2$)\\and $\tempGradMacro \in \{\unitVec_{(1)},\unitVec_{(2)}\}$ each $\Rightarrow \nSamples=2000$} & \makecell{$\solDParam, \rhsDParam$ for 400 2D microstructures ($400^2$)\\and $\tempGradMacro \in \{\unitVec_{(1)},\unitVec_{(2)}\}$ each $\Rightarrow \nSamples=800$} \\
        \midrule 
        \makecell{Thermal 3D (periodic BC) \\Thermal 3D (Dirichlet BC)} & \makecell{$\solDParam, \rhsDParam$ for 1000 3D microstructures ($192^3$)\\and $\tempGradMacro \in \{\unitVec_{(1)},\unitVec_{(2)},\unitVec_{(3)}\}$ each $\Rightarrow \nSamples=3000$} & \makecell{$\solDParam, \rhsDParam$ for 400 3D microstructures ($192^3$)\\and $\tempGradMacro \in \{\unitVec_{(1)},\unitVec_{(2)},\unitVec_{(3)}\}$ each $\Rightarrow \nSamples=1200$} \\
        \midrule 
        \makecell{Mechanical 2D (periodic BC)\\Mechanical 2D (Dirichlet BC)\\Mechanical 2D (mixed BC)} & \makecell{$\solDParam, \rhsDParam$ for 1000 2D microstructures ($400^2$)\\and $\strainMandelMacro \in \{\unitVec_{(1)},\unitVec_{(2)},\unitVec_{(3)}\}$ each $\Rightarrow \nSamples=3000$} & \makecell{$\solDParam, \rhsDParam$ for 400 microstructures ($400^2$)\\and $\strainMandelMacro \in \{\unitVec_{(1)},\unitVec_{(2)},\unitVec_{(3)}\}$ each $\Rightarrow \nSamples=1200$} \\
        \midrule 
        \makecell{Mechanical 3D (periodic BC)\\Mechanical 3D (Dirichlet BC)\\Mechanical 3D (mixed BC)} & \makecell{$\solDParam, \rhsDParam$ for 1000 3D microstructures ($192^3$)\\and $\strainMandelMacro \in \{\unitVec_{(1)},\dots,\unitVec_{(6)}\}$ each $\Rightarrow \nSamples=6000$} & \makecell{$\solDParam, \rhsDParam$ for 400 3D microstructures ($192^3$)\\and $\strainMandelMacro \in \{\unitVec_{(1)},\dots,\unitVec_{(6)}\}$ each $\Rightarrow \nSamples=2400$} \\
        \bottomrule
    \end{tabular}
\end{adjustbox}
    \label{tab:hyperparameters}
\end{table}

\bmsection{{\color{revision}Solver runtimes}}
\label{ssec:runtimes}
\begin{table}[ht!]
    \caption{\color{revision}Wall-clock runtimes required to reach a relative residual tolerance of $10^{-6}$
 for all considered solvers in \cref{sec:results}. Measurements were performed on an NVIDIA H100 GPU using identical matrix-free operator implementations. The timings include all solver components (matrix-free operators, preconditioner applications, vector updates), providing a fair comparison of the practical efficiency of \UNOCG{} relative to Jacobi-preconditioned and unpreconditioned CG as well as \FANS{} (for periodic BC only). Runtimes naturally depend on hardware and software implementations and are provided for reference.}
\begin{adjustbox}{center}
{\color{revision}
    \centering \begin{tabular}{crrrr}
        \toprule 
        Problem & Unpreconditioned \CG{} & \JacCG{} & \UNOCG{} & \FANS{} \\
        \midrule 
        Thermal 2D (periodic BC) & 466.3 ms & 406.8 ms & 11.6 ms & 8.4 ms \\
        \midrule
        Thermal 3D (periodic BC) & 10,808.9 ms & 10,304.2 ms & 872.5 ms & 547.6 ms \\
        Thermal 3D (Dirichlet BC) & 11,317.6 ms & 9,854.1 ms & 968.3 ms & --- \\
        \midrule 
        Mechanical 2D (periodic BC) & 1,630.3 ms & 1,217.8 ms & 41.1 ms & 18.8 ms \\
        Mechanical 2D (Dirichlet BC) & 1,549.9 ms & 1,222.0 ms & 73.7 ms & --- \\
        Mechanical 2D (mixed BC) & 2,093.5 ms & 1,706.5 ms & 61.2 ms & --- \\
        \midrule 
        Mechanical 3D (periodic BC) & 68,129.1 ms & 50,796.1 ms & 3,222.4 ms & 1,300.8 ms \\
        Mechanical 3D (Dirichlet BC) & 60,814.4 ms & 42,908.6 ms & 3,627.3 ms & --- \\
        Mechanical 3D (mixed BC) & 98,785.7 ms & 61,316.5 ms & 3,349.6 ms & --- \\
        \bottomrule
    \end{tabular}
}
\end{adjustbox}
    \label{tab:hyperparameters}
\end{table}

\end{document}